\documentstyle[eepic,12pt]{article}
\font\goth=eufm10 scaled 1200

\newtheorem{example}{Example}[section]

\newtheorem{theorem}[example]{Theorem}
\newtheorem{corollary}[example]{Corollary}

\def\boxit#1#2{\setbox1=\hbox{\kern#1{#2}\kern#1}%
\dimen1=\ht1 \advance\dimen1 by #1 \dimen2=\dp1 \advance\dimen2 by #1
\setbox1=\hbox{\vrule height\dimen1 depth\dimen2\box1\vrule}%
\setbox1=\vbox{\hrule\box1\hrule}%
\advance\dimen1 by .4pt \ht1=\dimen1
\advance\dimen2 by .4pt \dp1=\dimen2 \box1\relax}
\def\bo#1{\boxit{1pt}{$#1$}}

\def\Proof{\medskip\noindent {\it Proof: }}
\def\cqfd{\hfill $\Box$ \bigskip}
\def\adots{\mathinner{\mkern2mu\raise1pt\hbox{.}
\mkern3mu\raise4pt\hbox{.}\mkern1mu\raise7pt\hbox{.}}}
\def\<{\langle}
\def\>{\rangle}

\def\cf{{\it cf.$\ $}}
\def\resp{{\it resp.$\ $}}
\def\ie{{\it i.e. }}
\def\mat{{\rm Mat}}

\newfont{\bb}{cmbx10}

\def\Mat{{\rm Mat\, }}
\def\mod{{\rm mod\ }}

\def\L{{\cal L}}
\def\B{{\cal B}}

\def\Sl{\hbox{\goth sl}\hskip 2.5pt}

\def\gl{\hbox{\goth gl}\hskip 2.5pt}

\def\F{{\cal F}}
\def\A{{\cal A}}

\def\F{{\cal F}}
\def\dim{{\rm dim\,}}
\def\det{{\rm det\,}}

\def\F{{\cal F}}

\title{The Robinson-Schensted correspondence
\\ as the quantum straightening at $q=0$ \thanks{Partially supported
by PRC Math-Info and EEC grant n$^0$ ERBCHRXCT930400}}

\author{Bernard \sc   Leclerc\thanks{L.I.T.P., Universit\'e
Paris 7, 2 place Jussieu, 75251 Paris cedex 05, France}
\rm and Jean-Yves \sc Thibon\thanks{Institut Gaspard Monge, Universit\'e de
Marne-la-Vall\'ee, 2 rue de la Butte-Verte, 93166 Noisy-le-Grand cedex,
France}}
\date{\footnotesize March 1995.}

\begin{document}

\maketitle

\centerline{\it Dedicated to Dominique Foata}

\begin{abstract}
We show that the quantum straightening algorithm for Young tableaux
and Young bitableaux reduces in the crystal limit $q \mapsto 0$ to the
Robinson-Schensted algorithm.
\end{abstract}


\section{Introduction}\label{INT}

Let $K$ be a field of characteristic $0$, and $X= (x_{ij})_{1\le i,\,j \le n}$
a matrix of commutative indeterminates. The ring $K[x_{ij}]$ may be regarded
as the algebra $\F[\Mat_n]$ of polynomial functions on the space of $n\times n$
matrices over $K$. A linear basis of this algebra
is given by the bitableaux of D\'esarm\'enien, Kung, Rota \cite{DKR}, which
are defined in the following way.

Given two (semistandard) Young tableaux $\tau$ and $\tau'$ of the same shape,
with
columns $c_1,\ldots ,c_k$ and $c'_1,\ldots ,c'_k$, the Young
{\it bitableau} $(\tau\,|\,\tau')$ is the product of the $k$ minors of $X$
whose
row indices belong to $c_i$ and column indices to $c_i'$, $i=1,\ldots ,k$.
For example,
\begin{center}
\setlength{\unitlength}{0.01in}
\begin{picture}(202,100)(0,-10)
\put(116.000,40.000){\arc{170.000}{5.7932}{6.7731}}
\path(161,50)(181,50)(181,10)
	(161,10)(161,50)
\path(141,50)(161,50)(161,10)
	(141,10)(141,50)
\path(121,70)(141,70)(141,10)
	(121,10)(121,70)
\path(61,50)(81,50)(81,10)
	(61,10)(61,50)
\path(41,50)(61,50)(61,10)
	(41,10)(41,50)
\path(21,70)(41,70)(41,10)
	(21,10)(21,70)
\path(101,85)(101,0)
\put(171,15){\makebox(0,0)[b]{\raisebox{0pt}[0pt][0pt]{\shortstack{{5}}}}}
\put(171,35){\makebox(0,0)[b]{\raisebox{0pt}[0pt][0pt]{\shortstack{{7}}}}}
\put(86.000,40.000){\arc{170.000}{2.6516}{3.6315}}
\put(151,15){\makebox(0,0)[b]{\raisebox{0pt}[0pt][0pt]{\shortstack{{3}}}}}
\put(31,55){\makebox(0,0)[b]{\raisebox{0pt}[0pt][0pt]{\shortstack{{4}}}}}
\put(151,35){\makebox(0,0)[b]{\raisebox{0pt}[0pt][0pt]{\shortstack{{4}}}}}
\put(131,15){\makebox(0,0)[b]{\raisebox{0pt}[0pt][0pt]{\shortstack{{1}}}}}
\put(131,35){\makebox(0,0)[b]{\raisebox{0pt}[0pt][0pt]{\shortstack{{2}}}}}
\put(131,55){\makebox(0,0)[b]{\raisebox{0pt}[0pt][0pt]{\shortstack{{5}}}}}
\put(71,15){\makebox(0,0)[b]{\raisebox{0pt}[0pt][0pt]{\shortstack{{2}}}}}
\put(71,35){\makebox(0,0)[b]{\raisebox{0pt}[0pt][0pt]{\shortstack{{6}}}}}
\put(51,15){\makebox(0,0)[b]{\raisebox{0pt}[0pt][0pt]{\shortstack{{1}}}}}
\put(51,35){\makebox(0,0)[b]{\raisebox{0pt}[0pt][0pt]{\shortstack{{5}}}}}
\put(31,15){\makebox(0,0)[b]{\raisebox{0pt}[0pt][0pt]{\shortstack{{1}}}}}
\put(31,35){\makebox(0,0)[b]{\raisebox{0pt}[0pt][0pt]{\shortstack{{3}}}}}
\end{picture}
\end{center}
$$
:=
\left|\matrix{
x_{11} & x_{12} & x_{15}\cr
x_{31} & x_{32} & x_{35}\cr
x_{41} & x_{42} & x_{45}\cr
}\right|\times
\left|\matrix{
x_{13} & x_{14} \cr
x_{53} & x_{54} \cr
}\right|\times
\left|\matrix{
x_{25} & x_{27} \cr
x_{65} & x_{67} \cr
}\right|\,.
$$
More generally, we shall call {\it tabloid} a sequence of column-shaped Young
tableaux,
and we shall associate to each pair $\delta ,\, \delta'$ of tabloids of the
same
shape a {\it bitabloid} $(\delta \,|\,\delta')$ defined as the product of
minors
indexed by the columns of $\delta$ and $\delta'$.

There exists an algorithm due to D\'esarm\'enien \cite{D} for expanding any
polynomial
in $K[x_{ij}]$ on the basis of bitableaux. This is the so-called {\it
straightening
algorithm} (for bitableaux). In particular, the monomials
$x_{i_1j_1}\cdots x_{i_kj_k}$, which obviously form another linear basis of
$K[x_{ij}]$,
can be expressed in a unique way as linear combinations of bitableaux. Thus,
the straightening of $x_{23}\,x_{11}\,x_{32}$ reads
\begin{center}
\setlength{\unitlength}{0.01in}
\begingroup\makeatletter\ifx\SetFigFont\undefined
\def\x#1#2#3#4#5#6#7\relax{\def\x{#1#2#3#4#5#6}}%
\expandafter\x\fmtname xxxxxx\relax \def\y{splain}%
\ifx\x\y   
\gdef\SetFigFont#1#2#3{%
  \ifnum #1<17\tiny\else \ifnum #1<20\small\else
  \ifnum #1<24\normalsize\else \ifnum #1<29\large\else
  \ifnum #1<34\Large\else \ifnum #1<41\LARGE\else
     \huge\fi\fi\fi\fi\fi\fi
  \csname #3\endcsname}%
\else
\gdef\SetFigFont#1#2#3{\begingroup
  \count@#1\relax \ifnum 25<\count@\count@25\fi
  \def\x{\endgroup\@setsize\SetFigFont{#2pt}}%
  \expandafter\x
    \csname \romannumeral\the\count@ pt\expandafter\endcsname
    \csname @\romannumeral\the\count@ pt\endcsname
  \csname #3\endcsname}%
\fi
\fi\endgroup
\begin{picture}(596,175)(0,-10)
\put(535.000,130.000){\arc{100.000}{5.6397}{6.9267}}
\put(495.000,130.000){\arc{100.000}{2.4981}{3.7851}}
\path(475,110)(475,150)(455,150)
	(455,110)(475,110)
\path(495,110)(495,130)(475,130)
	(475,110)(495,110)
\path(555,110)(555,150)(535,150)
	(535,110)(555,110)
\path(575,110)(575,130)(555,130)
	(555,110)(575,110)
\path(515,160)(515,100)(515,105)
\put(460,115){\makebox(0,0)[lb]{\smash{{{\SetFigFont{12}{14.4}{rm}1}}}}}
\put(460,135){\makebox(0,0)[lb]{\smash{{{\SetFigFont{12}{14.4}{rm}3}}}}}
\put(480,115){\makebox(0,0)[lb]{\smash{{{\SetFigFont{12}{14.4}{rm}2}}}}}
\put(540,115){\makebox(0,0)[lb]{\smash{{{\SetFigFont{12}{14.4}{rm}1}}}}}
\put(540,135){\makebox(0,0)[lb]{\smash{{{\SetFigFont{12}{14.4}{rm}3}}}}}
\put(560,115){\makebox(0,0)[lb]{\smash{{{\SetFigFont{12}{14.4}{rm}2}}}}}
\put(260.000,120.000){\arc{50.000}{2.2143}{4.0689}}
\put(390.000,120.000){\arc{50.000}{5.3559}{7.2105}}
\path(265,110)(265,130)(245,130)
	(245,110)(265,110)
\path(285,110)(285,130)(265,130)
	(265,110)(285,110)
\path(305,110)(305,130)(285,130)
	(285,110)(305,110)
\path(365,110)(365,130)(345,130)
	(345,110)(365,110)
\path(385,110)(385,130)(365,130)
	(365,110)(385,110)
\path(405,110)(405,130)(385,130)
	(385,110)(405,110)
\path(325,140)(325,100)
\put(250,115){\makebox(0,0)[lb]{\smash{{{\SetFigFont{12}{14.4}{rm}1}}}}}
\put(270,115){\makebox(0,0)[lb]{\smash{{{\SetFigFont{12}{14.4}{rm}2}}}}}
\put(290,115){\makebox(0,0)[lb]{\smash{{{\SetFigFont{12}{14.4}{rm}3}}}}}
\put(350,115){\makebox(0,0)[lb]{\smash{{{\SetFigFont{12}{14.4}{rm}1}}}}}
\put(370,115){\makebox(0,0)[lb]{\smash{{{\SetFigFont{12}{14.4}{rm}2}}}}}
\put(390,115){\makebox(0,0)[lb]{\smash{{{\SetFigFont{12}{14.4}{rm}3}}}}}
\put(50.000,120.000){\arc{50.000}{2.2143}{4.0689}}
\put(180.000,120.000){\arc{50.000}{5.3559}{7.2105}}
\path(55,110)(55,130)(35,130)
	(35,110)(55,110)
\path(75,110)(75,130)(55,130)
	(55,110)(75,110)
\path(95,110)(95,130)(75,130)
	(75,110)(95,110)
\path(155,110)(155,130)(135,130)
	(135,110)(155,110)
\path(175,110)(175,130)(155,130)
	(155,110)(175,110)
\path(195,110)(195,130)(175,130)
	(175,110)(195,110)
\path(115,140)(115,100)
\put(40,115){\makebox(0,0)[lb]{\smash{{{\SetFigFont{12}{14.4}{rm}2}}}}}
\put(60,115){\makebox(0,0)[lb]{\smash{{{\SetFigFont{12}{14.4}{rm}1}}}}}
\put(80,115){\makebox(0,0)[lb]{\smash{{{\SetFigFont{12}{14.4}{rm}3}}}}}
\put(140,115){\makebox(0,0)[lb]{\smash{{{\SetFigFont{12}{14.4}{rm}3}}}}}
\put(160,115){\makebox(0,0)[lb]{\smash{{{\SetFigFont{12}{14.4}{rm}1}}}}}
\put(180,115){\makebox(0,0)[lb]{\smash{{{\SetFigFont{12}{14.4}{rm}2}}}}}
\put(580.000,40.000){\arc{170.000}{2.6516}{3.6315}}
\put(510.000,40.000){\arc{170.000}{5.7932}{6.7731}}
\put(100.000,30.000){\arc{100.000}{5.6397}{6.9267}}
\put(60.000,30.000){\arc{100.000}{2.4981}{3.7851}}
\put(260.000,30.000){\arc{100.000}{5.6397}{6.9267}}
\put(220.000,30.000){\arc{100.000}{2.4981}{3.7851}}
\put(420.000,30.000){\arc{100.000}{5.6397}{6.9267}}
\put(380.000,30.000){\arc{100.000}{2.4981}{3.7851}}
\path(575,10)(575,70)(555,70)
	(555,10)(575,10)
\path(535,10)(535,70)(515,70)
	(515,10)(535,10)
\path(545,80)(545,0)
\path(40,10)(40,50)(20,50)
	(20,10)(40,10)
\path(60,10)(60,30)(40,30)
	(40,10)(60,10)
\path(120,10)(120,50)(100,50)
	(100,10)(120,10)
\path(140,10)(140,30)(120,30)
	(120,10)(140,10)
\path(80,60)(80,0)(80,5)
\path(200,10)(200,50)(180,50)
	(180,10)(200,10)
\path(220,10)(220,30)(200,30)
	(200,10)(220,10)
\path(280,10)(280,50)(260,50)
	(260,10)(280,10)
\path(300,10)(300,30)(280,30)
	(280,10)(300,10)
\path(240,60)(240,0)(240,5)
\path(360,10)(360,50)(340,50)
	(340,10)(360,10)
\path(380,10)(380,30)(360,30)
	(360,10)(380,10)
\path(440,10)(440,50)(420,50)
	(420,10)(440,10)
\path(460,10)(460,30)(440,30)
	(440,10)(460,10)
\path(400,60)(400,0)(400,5)
\put(425,115){\makebox(0,0)[lb]{\smash{{{\SetFigFont{12}{14.4}{rm}$-$}}}}}
\put(216,115){\makebox(0,0)[lb]{\smash{{{\SetFigFont{12}{14.4}{rm}$=$}}}}}
\put(520,55){\makebox(0,0)[lb]{\smash{{{\SetFigFont{12}{14.4}{rm}3}}}}}
\put(520,35){\makebox(0,0)[lb]{\smash{{{\SetFigFont{12}{14.4}{rm}2}}}}}
\put(520,15){\makebox(0,0)[lb]{\smash{{{\SetFigFont{12}{14.4}{rm}1}}}}}
\put(560,55){\makebox(0,0)[lb]{\smash{{{\SetFigFont{12}{14.4}{rm}3}}}}}
\put(560,35){\makebox(0,0)[lb]{\smash{{{\SetFigFont{12}{14.4}{rm}2}}}}}
\put(560,15){\makebox(0,0)[lb]{\smash{{{\SetFigFont{12}{14.4}{rm}1}}}}}
\put(-5,15){\makebox(0,0)[lb]{\smash{{{\SetFigFont{12}{14.4}{rm}$+$}}}}}
\put(155,15){\makebox(0,0)[lb]{\smash{{{\SetFigFont{12}{14.4}{rm}$+$}}}}}
\put(315,15){\makebox(0,0)[lb]{\smash{{{\SetFigFont{12}{14.4}{rm}$-$}}}}}
\put(480,15){\makebox(0,0)[lb]{\smash{{{\SetFigFont{12}{14.4}{rm}$+$}}}}}
\put(25,15){\makebox(0,0)[lb]{\smash{{{\SetFigFont{12}{14.4}{rm}1}}}}}
\put(105,15){\makebox(0,0)[lb]{\smash{{{\SetFigFont{12}{14.4}{rm}1}}}}}
\put(105,35){\makebox(0,0)[lb]{\smash{{{\SetFigFont{12}{14.4}{rm}3}}}}}
\put(125,15){\makebox(0,0)[lb]{\smash{{{\SetFigFont{12}{14.4}{rm}2}}}}}
\put(185,15){\makebox(0,0)[lb]{\smash{{{\SetFigFont{12}{14.4}{rm}1}}}}}
\put(185,35){\makebox(0,0)[lb]{\smash{{{\SetFigFont{12}{14.4}{rm}3}}}}}
\put(205,15){\makebox(0,0)[lb]{\smash{{{\SetFigFont{12}{14.4}{rm}2}}}}}
\put(265,15){\makebox(0,0)[lb]{\smash{{{\SetFigFont{12}{14.4}{rm}1}}}}}
\put(345,15){\makebox(0,0)[lb]{\smash{{{\SetFigFont{12}{14.4}{rm}1}}}}}
\put(425,15){\makebox(0,0)[lb]{\smash{{{\SetFigFont{12}{14.4}{rm}1}}}}}
\put(25,35){\makebox(0,0)[lb]{\smash{{{\SetFigFont{12}{14.4}{rm}2}}}}}
\put(45,15){\makebox(0,0)[lb]{\smash{{{\SetFigFont{12}{14.4}{rm}3}}}}}
\put(265,35){\makebox(0,0)[lb]{\smash{{{\SetFigFont{12}{14.4}{rm}2}}}}}
\put(285,15){\makebox(0,0)[lb]{\smash{{{\SetFigFont{12}{14.4}{rm}3}}}}}
\put(345,35){\makebox(0,0)[lb]{\smash{{{\SetFigFont{12}{14.4}{rm}2}}}}}
\put(365,15){\makebox(0,0)[lb]{\smash{{{\SetFigFont{12}{14.4}{rm}3}}}}}
\put(425,35){\makebox(0,0)[lb]{\smash{{{\SetFigFont{12}{14.4}{rm}2}}}}}
\put(445,15){\makebox(0,0)[lb]{\smash{{{\SetFigFont{12}{14.4}{rm}3}}}}}
\end{picture}
\end{center}

On the other hand, the Robinson-Schensted correspondence \cite{Ro,Sch,S,Kn}
associates to any word $w$ on the alphabet of symbols $\{1,\,\ldots\,,n\}$
a pair $(P(w),\,Q(w))$ of Young tableaux of the same shape. For example,
the image of the word $w = 2\,1\,4\,3\,5\,1\,2$ under this correspondence
is the pair
\begin{center}
\setlength{\unitlength}{0.01in}
\begingroup\makeatletter\ifx\SetFigFont\undefined
\def\x#1#2#3#4#5#6#7\relax{\def\x{#1#2#3#4#5#6}}%
\expandafter\x\fmtname xxxxxx\relax \def\y{splain}%
\ifx\x\y   
\gdef\SetFigFont#1#2#3{%
  \ifnum #1<17\tiny\else \ifnum #1<20\small\else
  \ifnum #1<24\normalsize\else \ifnum #1<29\large\else
  \ifnum #1<34\Large\else \ifnum #1<41\LARGE\else
     \huge\fi\fi\fi\fi\fi\fi
  \csname #3\endcsname}%
\else
\gdef\SetFigFont#1#2#3{\begingroup
  \count@#1\relax \ifnum 25<\count@\count@25\fi
  \def\x{\endgroup\@setsize\SetFigFont{#2pt}}%
  \expandafter\x
    \csname \romannumeral\the\count@ pt\expandafter\endcsname
    \csname @\romannumeral\the\count@ pt\endcsname
  \csname #3\endcsname}%
\fi
\fi\endgroup
\begin{picture}(192,95)(0,-10)
\put(86.000,40.000){\arc{170.000}{2.6516}{3.6315}}
\put(106.000,40.000){\arc{170.000}{5.7932}{6.7731}}
\path(36,10)(36,70)(16,70)
	(16,10)(36,10)
\path(56,10)(56,50)(36,50)
	(36,10)(56,10)
\path(76,10)(76,50)(56,50)
	(56,10)(76,10)
\path(136,10)(136,70)(116,70)
	(116,10)(136,10)
\path(156,10)(156,50)(136,50)
	(136,10)(156,10)
\path(176,10)(176,50)(156,50)
	(156,10)(176,10)
\put(91,10){\makebox(0,0)[lb]{\smash{{{\SetFigFont{12}{14.4}{rm},}}}}}
\put(21,55){\makebox(0,0)[lb]{\smash{{{\SetFigFont{12}{14.4}{rm}4}}}}}
\put(21,35){\makebox(0,0)[lb]{\smash{{{\SetFigFont{12}{14.4}{rm}2}}}}}
\put(21,15){\makebox(0,0)[lb]{\smash{{{\SetFigFont{12}{14.4}{rm}1}}}}}
\put(41,35){\makebox(0,0)[lb]{\smash{{{\SetFigFont{12}{14.4}{rm}3}}}}}
\put(41,15){\makebox(0,0)[lb]{\smash{{{\SetFigFont{12}{14.4}{rm}1}}}}}
\put(61,35){\makebox(0,0)[lb]{\smash{{{\SetFigFont{12}{14.4}{rm}5}}}}}
\put(61,15){\makebox(0,0)[lb]{\smash{{{\SetFigFont{12}{14.4}{rm}2}}}}}
\put(121,55){\makebox(0,0)[lb]{\smash{{{\SetFigFont{12}{14.4}{rm}6}}}}}
\put(121,35){\makebox(0,0)[lb]{\smash{{{\SetFigFont{12}{14.4}{rm}2}}}}}
\put(121,15){\makebox(0,0)[lb]{\smash{{{\SetFigFont{12}{14.4}{rm}1}}}}}
\put(141,35){\makebox(0,0)[lb]{\smash{{{\SetFigFont{12}{14.4}{rm}4}}}}}
\put(141,15){\makebox(0,0)[lb]{\smash{{{\SetFigFont{12}{14.4}{rm}3}}}}}
\put(161,15){\makebox(0,0)[lb]{\smash{{{\SetFigFont{12}{14.4}{rm}5}}}}}
\put(161,35){\makebox(0,0)[lb]{\smash{{{\SetFigFont{12}{14.4}{rm}7}}}}}
\end{picture}
\end{center}

Both the straightening algorithm and the Robinson-Schensted algorithm are
strongly connected with the representation theory of $GL_n$. Indeed, the
straightening algorithm allows to compute the action of $GL_n$ in
(polynomial) irreducible representations, while the Robinson-Schensted
correspondence was devised by Robinson to obtain a proof of the
Littlewood-Ri\-chard\-son rule for decomposing into irreducibles the tensor
product of two irreducible representations.
The problem that we want to investigate is whether there exists any relation
between the straightening algorithm and the Schensted algorithm. It turns out
that to answer this question, one has to replace the algebra $\F[\Mat_n]$ by
its quantum analogue $\F_q[\Mat_n]$ \cite{RTF}.
This is the associative algebra over $K(q)$
generated by $n^2$ letters $t_{ij},\, i,\,j = 1,\,\ldots ,\, n$ subject to the
relations
\begin{eqnarray}
t_{ik}\,t_{il} & = & q^{-1}\,t_{il}\,t_{ik}\,,  \label{qGL1} \\
t_{ik}\,t_{jk} & = & q^{-1}\,t_{jk}\,t_{ik}\,,  \label{qGL2} \\
t_{il}\,t_{jk} & = & t_{jk}\,t_{il}\,,  \label{qGL3} \\
t_{ik}\,t_{jl}-t_{jl}\,t_{ik} & = & (q^{-1}\!-\! q)\, t_{il}\,t_{jk}\,,
\label{qGL4}
\end{eqnarray}
for $1\le i<j\le n$ and $1\le k<l\le n$. The {\it quantum determinant}
of $T=(t_{ij})$ is the element
$$
\det_q = \det_q\, T :=
\sum_{w\in{S_n}} \,
(-q)^{-\ell(w)} \, t_{1{w_1}}\,\ldots\, t_{n{w_n}}
$$
of $\F_q[\mat_n]$. Here $S_n$ is the symmetric group on $\{1,\ldots ,n\}$ and
$\ell(w)$ denotes the length of the permutation $w$.
The quantum determinant of $T$ belongs to the center of $\F_q[\mat_n]$.
More generally, for $I = ( i_1,\ldots ,i_k)$ and $J=(j_1,\ldots ,j_k)$
one defines the {\it quantum minor} of the submatrix $T_{IJ}$ as
$$
\det_q\, T_{IJ} :=
\sum_{w\in S_k} \,
(-q)^{-\ell(w)} \, t_{i_1j_{w_1}}\,\ldots\, t_{i_kj_{w_k}}
\ ,
$$
and the {\it quantum bitableau} $(\tau |\tau')$ (resp. {\it quantum bitabloid}
$(\delta \,|\,\delta')$)
as the product of the quantum minors indexed by the columns
of the Young tableaux $\tau$ and $\tau'$ (resp. of the tabloids
$\delta$ and $\delta'$). As proved by Huang, Zhang \cite{HZ}, the set of
quantum
bitableaux is again a linear basis of $\F_q[\Mat_n]$. For example, the
expansion of the monomial $t_{23}\,t_{11}\,t_{32}$ on this basis is

\begin{center}
\setlength{\unitlength}{0.01in}
\begingroup\makeatletter\ifx\SetFigFont\undefined
\def\x#1#2#3#4#5#6#7\relax{\def\x{#1#2#3#4#5#6}}%
\expandafter\x\fmtname xxxxxx\relax \def\y{splain}%
\ifx\x\y   
\gdef\SetFigFont#1#2#3{%
  \ifnum #1<17\tiny\else \ifnum #1<20\small\else
  \ifnum #1<24\normalsize\else \ifnum #1<29\large\else
  \ifnum #1<34\Large\else \ifnum #1<41\LARGE\else
     \huge\fi\fi\fi\fi\fi\fi
  \csname #3\endcsname}%
\else
\gdef\SetFigFont#1#2#3{\begingroup
  \count@#1\relax \ifnum 25<\count@\count@25\fi
  \def\x{\endgroup\@setsize\SetFigFont{#2pt}}%
  \expandafter\x
    \csname \romannumeral\the\count@ pt\expandafter\endcsname
    \csname @\romannumeral\the\count@ pt\endcsname
  \csname #3\endcsname}%
\fi
\fi\endgroup
\begin{picture}(631,255)(0,-10)
\put(50.000,200.000){\arc{50.000}{2.2143}{4.0689}}
\put(180.000,200.000){\arc{50.000}{5.3559}{7.2105}}
\path(55,190)(55,210)(35,210)
	(35,190)(55,190)
\path(75,190)(75,210)(55,210)
	(55,190)(75,190)
\path(95,190)(95,210)(75,210)
	(75,190)(95,190)
\path(155,190)(155,210)(135,210)
	(135,190)(155,190)
\path(175,190)(175,210)(155,210)
	(155,190)(175,190)
\path(195,190)(195,210)(175,210)
	(175,190)(195,190)
\path(115,220)(115,180)
\put(40,195){\makebox(0,0)[lb]{\smash{{{\SetFigFont{12}{14.4}{rm}2}}}}}
\put(60,195){\makebox(0,0)[lb]{\smash{{{\SetFigFont{12}{14.4}{rm}1}}}}}
\put(80,195){\makebox(0,0)[lb]{\smash{{{\SetFigFont{12}{14.4}{rm}3}}}}}
\put(140,195){\makebox(0,0)[lb]{\smash{{{\SetFigFont{12}{14.4}{rm}3}}}}}
\put(160,195){\makebox(0,0)[lb]{\smash{{{\SetFigFont{12}{14.4}{rm}1}}}}}
\put(180,195){\makebox(0,0)[lb]{\smash{{{\SetFigFont{12}{14.4}{rm}2}}}}}
\put(290.000,200.000){\arc{50.000}{2.2143}{4.0689}}
\put(420.000,200.000){\arc{50.000}{5.3559}{7.2105}}
\path(295,190)(295,210)(275,210)
	(275,190)(295,190)
\path(315,190)(315,210)(295,210)
	(295,190)(315,190)
\path(335,190)(335,210)(315,210)
	(315,190)(335,190)
\path(395,190)(395,210)(375,210)
	(375,190)(395,190)
\path(415,190)(415,210)(395,210)
	(395,190)(415,190)
\path(435,190)(435,210)(415,210)
	(415,190)(435,190)
\path(355,220)(355,180)
\put(280,195){\makebox(0,0)[lb]{\smash{{{\SetFigFont{12}{14.4}{rm}1}}}}}
\put(300,195){\makebox(0,0)[lb]{\smash{{{\SetFigFont{12}{14.4}{rm}2}}}}}
\put(320,195){\makebox(0,0)[lb]{\smash{{{\SetFigFont{12}{14.4}{rm}3}}}}}
\put(380,195){\makebox(0,0)[lb]{\smash{{{\SetFigFont{12}{14.4}{rm}1}}}}}
\put(400,195){\makebox(0,0)[lb]{\smash{{{\SetFigFont{12}{14.4}{rm}2}}}}}
\put(420,195){\makebox(0,0)[lb]{\smash{{{\SetFigFont{12}{14.4}{rm}3}}}}}
\put(575.000,130.000){\arc{100.000}{5.6397}{6.9267}}
\put(535.000,130.000){\arc{100.000}{2.4981}{3.7851}}
\path(515,110)(515,150)(495,150)
	(495,110)(515,110)
\path(535,110)(535,130)(515,130)
	(515,110)(535,110)
\path(595,110)(595,150)(575,150)
	(575,110)(595,110)
\path(615,110)(615,130)(595,130)
	(595,110)(615,110)
\path(555,160)(555,100)(555,105)
\put(500,115){\makebox(0,0)[lb]{\smash{{{\SetFigFont{12}{14.4}{rm}1}}}}}
\put(580,115){\makebox(0,0)[lb]{\smash{{{\SetFigFont{12}{14.4}{rm}1}}}}}
\put(500,135){\makebox(0,0)[lb]{\smash{{{\SetFigFont{12}{14.4}{rm}2}}}}}
\put(520,115){\makebox(0,0)[lb]{\smash{{{\SetFigFont{12}{14.4}{rm}3}}}}}
\put(580,135){\makebox(0,0)[lb]{\smash{{{\SetFigFont{12}{14.4}{rm}2}}}}}
\put(600,115){\makebox(0,0)[lb]{\smash{{{\SetFigFont{12}{14.4}{rm}3}}}}}
\put(380.000,40.000){\arc{170.000}{2.6516}{3.6315}}
\put(310.000,40.000){\arc{170.000}{5.7932}{6.7731}}
\path(375,10)(375,70)(355,70)
	(355,10)(375,10)
\path(335,10)(335,70)(315,70)
	(315,10)(335,10)
\path(345,80)(345,0)
\put(320,55){\makebox(0,0)[lb]{\smash{{{\SetFigFont{12}{14.4}{rm}3}}}}}
\put(320,35){\makebox(0,0)[lb]{\smash{{{\SetFigFont{12}{14.4}{rm}2}}}}}
\put(320,15){\makebox(0,0)[lb]{\smash{{{\SetFigFont{12}{14.4}{rm}1}}}}}
\put(360,55){\makebox(0,0)[lb]{\smash{{{\SetFigFont{12}{14.4}{rm}3}}}}}
\put(360,35){\makebox(0,0)[lb]{\smash{{{\SetFigFont{12}{14.4}{rm}2}}}}}
\put(360,15){\makebox(0,0)[lb]{\smash{{{\SetFigFont{12}{14.4}{rm}1}}}}}
\put(375.000,130.000){\arc{100.000}{5.6397}{6.9267}}
\put(335.000,130.000){\arc{100.000}{2.4981}{3.7851}}
\path(315,110)(315,150)(295,150)
	(295,110)(315,110)
\path(335,110)(335,130)(315,130)
	(315,110)(335,110)
\path(395,110)(395,150)(375,150)
	(375,110)(395,110)
\path(415,110)(415,130)(395,130)
	(395,110)(415,110)
\path(355,160)(355,100)(355,105)
\put(300,115){\makebox(0,0)[lb]{\smash{{{\SetFigFont{12}{14.4}{rm}1}}}}}
\put(300,135){\makebox(0,0)[lb]{\smash{{{\SetFigFont{12}{14.4}{rm}3}}}}}
\put(320,115){\makebox(0,0)[lb]{\smash{{{\SetFigFont{12}{14.4}{rm}2}}}}}
\put(380,115){\makebox(0,0)[lb]{\smash{{{\SetFigFont{12}{14.4}{rm}1}}}}}
\put(380,135){\makebox(0,0)[lb]{\smash{{{\SetFigFont{12}{14.4}{rm}2}}}}}
\put(400,115){\makebox(0,0)[lb]{\smash{{{\SetFigFont{12}{14.4}{rm}3}}}}}
\put(195.000,130.000){\arc{100.000}{5.6397}{6.9267}}
\put(155.000,130.000){\arc{100.000}{2.4981}{3.7851}}
\path(135,110)(135,150)(115,150)
	(115,110)(135,110)
\path(155,110)(155,130)(135,130)
	(135,110)(155,110)
\path(215,110)(215,150)(195,150)
	(195,110)(215,110)
\path(235,110)(235,130)(215,130)
	(215,110)(235,110)
\path(175,160)(175,100)(175,105)
\put(120,115){\makebox(0,0)[lb]{\smash{{{\SetFigFont{12}{14.4}{rm}1}}}}}
\put(200,115){\makebox(0,0)[lb]{\smash{{{\SetFigFont{12}{14.4}{rm}1}}}}}
\put(200,135){\makebox(0,0)[lb]{\smash{{{\SetFigFont{12}{14.4}{rm}3}}}}}
\put(220,115){\makebox(0,0)[lb]{\smash{{{\SetFigFont{12}{14.4}{rm}2}}}}}
\put(120,135){\makebox(0,0)[lb]{\smash{{{\SetFigFont{12}{14.4}{rm}2}}}}}
\put(140,115){\makebox(0,0)[lb]{\smash{{{\SetFigFont{12}{14.4}{rm}3}}}}}
\put(580.000,210.000){\arc{100.000}{5.6397}{6.9267}}
\put(540.000,210.000){\arc{100.000}{2.4981}{3.7851}}
\path(520,190)(520,230)(500,230)
	(500,190)(520,190)
\path(540,190)(540,210)(520,210)
	(520,190)(540,190)
\path(600,190)(600,230)(580,230)
	(580,190)(600,190)
\path(620,190)(620,210)(600,210)
	(600,190)(620,190)
\path(560,240)(560,180)(560,185)
\put(505,195){\makebox(0,0)[lb]{\smash{{{\SetFigFont{12}{14.4}{rm}1}}}}}
\put(505,215){\makebox(0,0)[lb]{\smash{{{\SetFigFont{12}{14.4}{rm}3}}}}}
\put(525,195){\makebox(0,0)[lb]{\smash{{{\SetFigFont{12}{14.4}{rm}2}}}}}
\put(585,195){\makebox(0,0)[lb]{\smash{{{\SetFigFont{12}{14.4}{rm}1}}}}}
\put(585,215){\makebox(0,0)[lb]{\smash{{{\SetFigFont{12}{14.4}{rm}3}}}}}
\put(605,195){\makebox(0,0)[lb]{\smash{{{\SetFigFont{12}{14.4}{rm}2}}}}}
\put(220,195){\makebox(0,0)[lb]{\smash{{{\SetFigFont{12}{14.4}{rm}$=$}}}}}
\put(460,195){\makebox(0,0)[lb]{\smash{{{\SetFigFont{12}{14.4}{rm}$-$}}}}}
\put(0,115){\makebox(0,0)[lb]{\smash{{{\SetFigFont{12}{14.4}{rm}$+$}}}}}
\put(440,115){\makebox(0,0)[lb]{\smash{{{\SetFigFont{12}{14.4}{rm}$-$}}}}}
\put(245,195){\makebox(0,0)[lb]{\smash{{{\SetFigFont{12}{14.4}{rm}$q^3$}}}}}
\put(475,195){\makebox(0,0)[lb]{\smash{{{\SetFigFont{12}{14.4}{rm}$q^3$}}}}}
\put(15,115){\makebox(0,0)[lb]{\smash{{{\SetFigFont{12}{14.4}{rm}$(1-q^2+q^4)$}}}}}
\put(270,115){\makebox(0,0)[lb]{\smash{{{\SetFigFont{12}{14.4}{rm}$q^4$}}}}}
\put(275,35){\makebox(0,0)[lb]{\smash{{{\SetFigFont{12}{14.4}{rm}$q^5$}}}}}
\put(255,35){\makebox(0,0)[lb]{\smash{{{\SetFigFont{12}{14.4}{rm}$+$}}}}}
\put(465,115){\makebox(0,0)[lb]{\smash{{{\SetFigFont{12}{14.4}{rm}$q^5$}}}}}
\put(255,115){\makebox(0,0)[lb]{\smash{{{\SetFigFont{12}{14.4}{rm}$+$}}}}}
\end{picture}
\end{center}
In the present example the coefficients of the expansion are
polynomials in $q$, and only one of them has a nonzero constant term. In other
words, denote by $\B$ the linear basis of quantum bitableaux in $\F_q[\Mat_n]$,
and let $\L$ be the lattice generated over $K[q]$ by the elements of $\B$. Then
\begin{center}
\setlength{\unitlength}{0.01in}
\begingroup\makeatletter\ifx\SetFigFont\undefined
\def\x#1#2#3#4#5#6#7\relax{\def\x{#1#2#3#4#5#6}}%
\expandafter\x\fmtname xxxxxx\relax \def\y{splain}%
\ifx\x\y   
\gdef\SetFigFont#1#2#3{%
  \ifnum #1<17\tiny\else \ifnum #1<20\small\else
  \ifnum #1<24\normalsize\else \ifnum #1<29\large\else
  \ifnum #1<34\Large\else \ifnum #1<41\LARGE\else
     \huge\fi\fi\fi\fi\fi\fi
  \csname #3\endcsname}%
\else
\gdef\SetFigFont#1#2#3{\begingroup
  \count@#1\relax \ifnum 25<\count@\count@25\fi
  \def\x{\endgroup\@setsize\SetFigFont{#2pt}}%
  \expandafter\x
    \csname \romannumeral\the\count@ pt\expandafter\endcsname
    \csname @\romannumeral\the\count@ pt\endcsname
  \csname #3\endcsname}%
\fi
\fi\endgroup
\begin{picture}(393,75)(0,-10)
\put(26.000,20.000){\arc{50.000}{2.2143}{4.0689}}
\put(156.000,20.000){\arc{50.000}{5.3559}{7.2105}}
\path(31,10)(31,30)(11,30)
	(11,10)(31,10)
\path(51,10)(51,30)(31,30)
	(31,10)(51,10)
\path(71,10)(71,30)(51,30)
	(51,10)(71,10)
\path(131,10)(131,30)(111,30)
	(111,10)(131,10)
\path(151,10)(151,30)(131,30)
	(131,10)(151,10)
\path(171,10)(171,30)(151,30)
	(151,10)(171,10)
\path(91,40)(91,0)
\put(16,15){\makebox(0,0)[lb]{\smash{{{\SetFigFont{12}{14.4}{rm}2}}}}}
\put(36,15){\makebox(0,0)[lb]{\smash{{{\SetFigFont{12}{14.4}{rm}1}}}}}
\put(56,15){\makebox(0,0)[lb]{\smash{{{\SetFigFont{12}{14.4}{rm}3}}}}}
\put(116,15){\makebox(0,0)[lb]{\smash{{{\SetFigFont{12}{14.4}{rm}3}}}}}
\put(136,15){\makebox(0,0)[lb]{\smash{{{\SetFigFont{12}{14.4}{rm}1}}}}}
\put(156,15){\makebox(0,0)[lb]{\smash{{{\SetFigFont{12}{14.4}{rm}2}}}}}
\put(291.000,30.000){\arc{100.000}{5.6397}{6.9267}}
\put(251.000,30.000){\arc{100.000}{2.4981}{3.7851}}
\path(231,10)(231,50)(211,50)
	(211,10)(231,10)
\path(251,10)(251,30)(231,30)
	(231,10)(251,10)
\path(311,10)(311,50)(291,50)
	(291,10)(311,10)
\path(331,10)(331,30)(311,30)
	(311,10)(331,10)
\path(271,60)(271,0)(271,5)
\put(216,15){\makebox(0,0)[lb]{\smash{{{\SetFigFont{12}{14.4}{rm}1}}}}}
\put(296,15){\makebox(0,0)[lb]{\smash{{{\SetFigFont{12}{14.4}{rm}1}}}}}
\put(296,35){\makebox(0,0)[lb]{\smash{{{\SetFigFont{12}{14.4}{rm}3}}}}}
\put(316,15){\makebox(0,0)[lb]{\smash{{{\SetFigFont{12}{14.4}{rm}2}}}}}
\put(216,35){\makebox(0,0)[lb]{\smash{{{\SetFigFont{12}{14.4}{rm}2}}}}}
\put(236,15){\makebox(0,0)[lb]{\smash{{{\SetFigFont{12}{14.4}{rm}3}}}}}
\put(186,15){\makebox(0,0)[lb]{\smash{{{\SetFigFont{12}{14.4}{rm}$\equiv$}}}}}
\put(356,15){\makebox(0,0)[lb]{\smash{{{\SetFigFont{12}{14.4}{rm}$\mod
q\L$}}}}}
\end{picture}
\end{center}
This is an illustration of the following
\begin{theorem} \label{TH1}
Let $w=i_1\cdots i_k$ and $u=j_1\cdots j_k$ be two words
on $\{1,\ldots ,n\}$, and denote by $(P(w),\,Q(w))$, $(P(u),\,Q(u))$
their images under the Robinson-Schensted correspondence. Then,
$$
t_{i_1j_1}\cdots t_{i_kj_k} \equiv
\left\{\matrix{
(P(w)|P(u)) &\mod q\L &{\it if}\  Q(w)=Q(u) \cr
0           &\mod q\L &{\it otherwise}} \right.\,.
$$
\end{theorem}

As a corollary we obtain an unexpected characterization of the {\it plactic
congruence} on words \cite{Kn,LS1}, defined by
$$
w \sim u \quad \Longleftrightarrow \quad P(w) = P(u)\,.
$$
\begin{corollary} With the same notations as above,
$$
w \sim u \quad \Longleftrightarrow \quad t_{i_1i_1}\cdots t_{i_ki_k} \equiv
t_{j_1j_1}\cdots t_{j_kj_k}\ \mod q\L\,.
$$
\end{corollary}

At this point, it is important to recall that the existence of a connection
between the Robinson-Schensted correspondence and the representation
theory of the quantized enveloping algebra $U_q(\gl_n)$ at $q=0$ was
first discovered by Date, Jimbo, Miwa \cite{DJM}. The results presented here
are in fact of the same kind as those of \cite{DJM}, namely, it is shown
in \cite{DJM} that if $V_{(1)}$ denotes the basic representation of
$U_q(\gl_n)$,
the transition matrix in $V_{(1)}^{\otimes k}$ from the basis of
monomial tensors to the Gelfand-Zetlin basis specializes when $q=0$
to a permutation matrix given by the Robinson-Schensted map.
Similarly, Theorem~\ref{TH1} states that the transition matrix in the
$U_q(\gl_n)$-module $\F_q[\Mat_n]$ from the basis of monomials
$$
B_\tau =\{t_{i_1j_1}\cdots t_{i_kj_k}\ | \ Q(i_1\cdots i_k) = Q(j_1\cdots j_k)
=
\tau_\nu \ {\rm for \ some}\ \nu\}
$$
to the basis of quantum bitableaux is equal at $q=0$ to a permutation matrix
also computed from the Robinson-Schensted algorithm.
Here, $\tau_\nu$ denotes for each partition $\nu$ a fixed standard Young
tableau
of shape $\nu$.

The work of  Date, Jimbo, Miwa, provided the
starting point from which Kashiwara developed
his theory of {\it crystal bases} for quantized enveloping algebras
\cite{Ka1,Ka2}.
We shall use crystal bases as the main tool for proving Theorem~\ref{TH1}.

The paper is organized as follows. In Section~\ref{Uq}, we collect the
necessary
material about $U_q(\gl_n)$, $U_q(\Sl_n)$ and their representation theory.
In Section~\ref{BC}, we review the definition and basic
properties of Kashiwara's crystal bases at $q=0$. In Section~\ref{qG/B},
we formulate and prove a version of Theorem~\ref{TH1} for the case of
(single) tableaux, that is, we work in the subring of $\F_q[\Mat_n]$
generated by the quantum minors taken on the initial rows of $T$. Finally in
Section~\ref{PR} we prove Theorem~\ref{TH1}, as well as a slightly more
general statement.


\section{$U_q(\gl_n)$ and $U_q(\Sl_n)$} \label{Uq}

A general reference for this Section and the following one is the excellent
exposition \cite{ChPr}. We first recall the definition of the quantized
enveloping algebras $U_q(\gl_n)$ \cite{Ji2} and $U_q(\Sl_n)$ \cite{Ji1,Dr}.
$U_q(\gl_n)$ is the associative algebra over $K(q)$ generated by the
$4n-2$ symbols $e_i,\,f_i,\,i=1,\ldots , n-1$ and $q^{\epsilon_i},\,
q^{-\epsilon_i},\,i=1,\ldots , n$, subject to the relations
\begin{equation}
q^{\epsilon_i} q^{-\epsilon_i} = q^{-\epsilon_i} q^{\epsilon_i} =1\,,
\quad [q^{\epsilon_i},q^{\epsilon_j}] = 0\,,
\end{equation}
\begin{equation}\label{Uq2}
q^{\epsilon_i} e_j q^{-\epsilon_i} =
\left\{\matrix{q e_j     &{\rm for} \ i=j \cr
	       q^{-1} e_j &{\rm for} \ i=j+1 \cr
               e_j       &{\rm otherwise}  \cr
	       }\right.\,,
\end{equation}
\begin{equation}\label{Uq3}
q^{\epsilon_i} f_j q^{-\epsilon_i} =
\left\{\matrix{q^{-1} f_j     &{\rm for} \ i=j \cr
	       q f_j &{\rm for} \ i=j+1 \cr
	      f_j       &{\rm otherwise}  \cr
	     }\right.\,,
\end{equation}
\begin{equation}\label{Uq4}
[e_i,f_j]= \delta_{ij} {q^{\epsilon_i}q^{-\epsilon_{i+1}}-
q^{-\epsilon_i}q^{\epsilon_{i+1}}\over q-q^{-1}}\,,
\end{equation}
\begin{equation}
[e_i,e_j] = [f_i,f_j] = 0 \ \ {\rm for} \ |i-j|>1\,,
\end{equation}
\begin{equation}
e_je_i^2 -(q+q^{-1})e_ie_je_i +e_i^2e_j =
f_jf_i^2 -(q+q^{-1})f_if_jf_i +f_i^2f_j =0 \ {\rm for} \ |i-j|=1  \,.
\end{equation}
The subalgebra of $U_q(\gl_n)$ generated by
$e_i,\,f_i$, and
\begin{equation}\label{Uq7}
q^{h_i}=q^{\epsilon_i}q^{-\epsilon_{i+1}},\quad
q^{-h_i}=q^{-\epsilon_i}q^{\epsilon_{i+1}},\quad i=1,\ldots , n-1,
\end{equation}
is denoted by $U_q(\Sl_n)$.

The representation theories of $U_q(\gl_n)$ and $U_q(\Sl_n)$ are closely
parallel to those of their classical counterparts $U(\gl_n)$ and $U(\Sl_n)$.
Let $M$ be a $U_q(\gl_n)$-module and $\mu = (\mu_1,\ldots ,\mu_n)$ be a
$n$-tuple of nonnegative integers. The subspace
$$
M_\mu=\{v\in M\ |\ q^{\epsilon_i} v=q^{\mu_i}v\,,\quad i=1,\ldots , n\}
$$
is called a {\it weight space} and its elements are called {\it weight vectors}
(of weight $\mu$). Relations (\ref{Uq2}) (\ref{Uq3}) show that
$$
e_i M_\mu \subset M_{\mu^+}\,, \quad \ f_i M_\mu \subset M_{\mu^-} \,,
$$
where $\mu^+ = (\mu_1,\ldots ,\mu_i+1,\mu_{i+1}-1,\ldots ,\mu_n)$
and $\mu^- = (\mu_1,\ldots ,\mu_i-1,\mu_{i+1}+1,\ldots ,\mu_n)$. Thus,
ordering the weights in the usual way by setting
$$
\mu \le \lambda \quad \Longleftrightarrow \quad
\sum_{i=1}^k \mu_i \le \sum_{i=1}^k \lambda_i\,,
\quad \ k=1,\ldots ,n,
$$
we see that the $e_i$'s act as raising operators and the $f_i$'s as
lowering operators. A weight vector is said to be a {\it highest weight vector}
if it is annihilated by the $e_i$'s.
$M$ is called a {\it highest weight module} if it contains a highest weight
vector
$v$ such that $M = U_q(\gl_n)\,v$. If $v$ is of weight $\lambda$, it
follows that $\dim M_{\lambda} = 1$ and $M=\oplus_{\mu \le \lambda} M_{\mu}$.
One then shows that there exists for each partition $\lambda$ of length $\le n$
a unique highest weight finite-dimensional irreducible $U_q(\gl_n)$-module
$V_{\lambda}$, with highest weight $\lambda$.

\begin{example}\label{BAS} {\rm
The basic representation $V=V_{(1)}$ of $U_q(\gl_n)$ is the $n$-dimensional
vector space over $K(q)$ with basis $\{v_i,\ 1\le i \le n\}$, on which the
action of $U_q(\gl_n)$ is as follows:
$$
q^{\epsilon_i} v_j  =  q^{\delta_{ij}} v_j ,\quad
e_i   v_j  =  \delta_{i+1\,j} v_i ,\quad
f_i   v_j  =  \delta_{ij} v_{i+1} \,.
$$
}
\end{example}

\begin{example} \label{PEXT} {\rm
More generally,
the $U_q(\gl_n)$-module $V_{(1^k)}$ is a
${n\choose k}$-dimensional vector space with basis $\{v_c\}$ labelled by the
subsets
$c$ of $\{1,\ldots ,n\}$ with $k$ elements (\ie by the Young tableaux of shape
$(1^k)$ over $\{1,\ldots ,n\}$). The action of $U_q(\gl_n)$ on this basis is
given by
$$
q^{\epsilon_i} v_c =
\left\{\matrix{v_c  &{\rm if}\  i\not\in c\cr
	       q v_c &{\rm otherwise}  \cr
}\right.\,,
$$
$$
e_i v_c =
\left\{\matrix{0  &{\rm if}\  i+1\not\in c \ {\rm or}\ i\in c\cr
	        v_d &{\rm otherwise,\ where}\ d = (c\setminus \{i+1\}) \cup \{i\}
\cr
}\right.\,,
$$
$$
f_i v_c =
\left\{\matrix{0  &{\rm if}\  i+1\in c \ {\rm or}\ i\not\in c\cr
		v_d &{\rm otherwise,\ where}\ d = (c\setminus \{i\}) \cup \{i+1\}   \cr
}\right.\,.
$$
We see that the action of the lowering operators $f_i$ does not depend on
$q$, and can be recorded on a colored graph whose vertices are the
column-shaped
Young tableaux $c$ and whose arrows are given by:
$$
c \stackrel{i}{\longrightarrow} d \quad \Longleftrightarrow \quad f_i v_c = v_d
\,.
$$
Thus for $k=2, \ n=4$, one has the following graph:
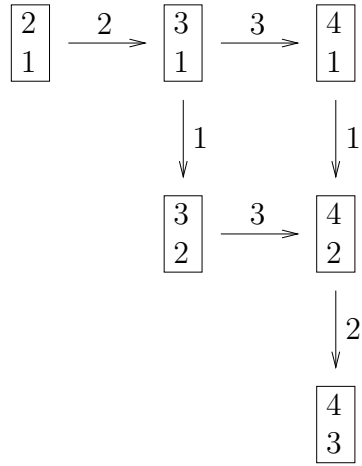
\begin{center}
\setlength{\unitlength}{0.01in}
\begingroup\makeatletter\ifx\SetFigFont\undefined
\def\x#1#2#3#4#5#6#7\relax{\def\x{#1#2#3#4#5#6}}%
\expandafter\x\fmtname xxxxxx\relax \def\y{splain}%
\ifx\x\y   
\gdef\SetFigFont#1#2#3{%
  \ifnum #1<17\tiny\else \ifnum #1<20\small\else
  \ifnum #1<24\normalsize\else \ifnum #1<29\large\else
  \ifnum #1<34\Large\else \ifnum #1<41\LARGE\else
     \huge\fi\fi\fi\fi\fi\fi
  \csname #3\endcsname}%
\else
\gdef\SetFigFont#1#2#3{\begingroup
  \count@#1\relax \ifnum 25<\count@\count@25\fi
  \def\x{\endgroup\@setsize\SetFigFont{#2pt}}%
  \expandafter\x
    \csname \romannumeral\the\count@ pt\expandafter\endcsname
    \csname @\romannumeral\the\count@ pt\endcsname
  \csname #3\endcsname}%
\fi
\fi\endgroup
\begin{picture}(181,255)(0,-10)
\path(20,200)(20,240)(0,240)
	(0,200)(20,200)
\path(100,200)(100,240)(80,240)
	(80,200)(100,200)
\path(180,200)(180,240)(160,240)
	(160,200)(180,200)
\path(100,100)(100,140)(80,140)
	(80,100)(100,100)
\path(180,100)(180,140)(160,140)
	(160,100)(180,100)
\path(180,0)(180,40)(160,40)
	(160,0)(180,0)
\path(30,220)(70,220)
\path(62.000,218.000)(70.000,220.000)(62.000,222.000)
\path(110,220)(150,220)
\path(142.000,218.000)(150.000,220.000)(142.000,222.000)
\path(90,190)(90,150)
\path(88.000,158.000)(90.000,150.000)(92.000,158.000)
\path(170,190)(170,150)
\path(168.000,158.000)(170.000,150.000)(172.000,158.000)
\path(110,120)(150,120)
\path(142.000,118.000)(150.000,120.000)(142.000,122.000)
\path(170,90)(170,50)
\path(168.000,58.000)(170.000,50.000)(172.000,58.000)
\put(5,205){\makebox(0,0)[lb]{\smash{{{\SetFigFont{12}{14.4}{rm}1}}}}}
\put(5,225){\makebox(0,0)[lb]{\smash{{{\SetFigFont{12}{14.4}{rm}2}}}}}
\put(85,205){\makebox(0,0)[lb]{\smash{{{\SetFigFont{12}{14.4}{rm}1}}}}}
\put(85,225){\makebox(0,0)[lb]{\smash{{{\SetFigFont{12}{14.4}{rm}3}}}}}
\put(165,205){\makebox(0,0)[lb]{\smash{{{\SetFigFont{12}{14.4}{rm}1}}}}}
\put(165,225){\makebox(0,0)[lb]{\smash{{{\SetFigFont{12}{14.4}{rm}4}}}}}
\put(85,105){\makebox(0,0)[lb]{\smash{{{\SetFigFont{12}{14.4}{rm}2}}}}}
\put(85,125){\makebox(0,0)[lb]{\smash{{{\SetFigFont{12}{14.4}{rm}3}}}}}
\put(165,105){\makebox(0,0)[lb]{\smash{{{\SetFigFont{12}{14.4}{rm}2}}}}}
\put(165,125){\makebox(0,0)[lb]{\smash{{{\SetFigFont{12}{14.4}{rm}4}}}}}
\put(165,5){\makebox(0,0)[lb]{\smash{{{\SetFigFont{12}{14.4}{rm}3}}}}}
\put(165,25){\makebox(0,0)[lb]{\smash{{{\SetFigFont{12}{14.4}{rm}4}}}}}
\put(45,225){\makebox(0,0)[lb]{\smash{{{\SetFigFont{12}{14.4}{rm}2}}}}}
\put(125,225){\makebox(0,0)[lb]{\smash{{{\SetFigFont{12}{14.4}{rm}3}}}}}
\put(95,165){\makebox(0,0)[lb]{\smash{{{\SetFigFont{12}{14.4}{rm}1}}}}}
\put(175,165){\makebox(0,0)[lb]{\smash{{{\SetFigFont{12}{14.4}{rm}1}}}}}
\put(125,125){\makebox(0,0)[lb]{\smash{{{\SetFigFont{12}{14.4}{rm}3}}}}}
\put(175,65){\makebox(0,0)[lb]{\smash{{{\SetFigFont{12}{14.4}{rm}2}}}}}
\end{picture}
\end{center}
This is one of the simplest examples of crystal graphs (\cf Section~\ref{BC}).
}
\end{example}

In order to construct more interesting $U_q(\gl_n)$-modules, we use the tensor
product operation. Given two $U_q(\gl_n)$-modules $M,\ N$, we can define a
structure of $U_q(\gl_n)$-module on $M\otimes N$ by putting
\begin{eqnarray}
q^{\epsilon_i} (u\otimes v) & = & q^{\epsilon_i} u \otimes q^{\epsilon_i} v,
\label{ATP1}\\
e_i (u\otimes v) & = & e_i u \otimes v + q^{-h_i} u \otimes e_i v,
\label{ATP2}\\
f_i (u\otimes v) & = & f_i u \otimes q^{h_i} v + u \otimes f_i v. \label{ATP3}
\end{eqnarray}
Indeed, the formulas
$$
\Delta q^{\epsilon_i}  =  q^{\epsilon_i} \otimes q^{\epsilon_i} , \quad
\Delta e_i   =  e_i  \otimes 1 + q^{-h_i} \otimes e_i , \quad
\Delta f_i   =  f_i  \otimes q^{h_i}  + 1 \otimes f_i ,
$$
define a comultiplication on $U_q(\gl_n)$. One shows that the decomposition
into irreducible components of the tensor product of two irreducible
$U_q(\gl_n)$-modules is given by
\begin{equation}\label{LR}
V_\lambda \otimes V_\mu \simeq \bigoplus_\nu\, c_{\lambda\,\mu}^\nu\,V_\nu \,,
\end{equation}
where the $c_{\lambda\,\mu}^\nu$ are the classical Littlewood-Richardson
numbers. In particular, it follows that
\begin{equation}\label{TS}
V^{\otimes k} \simeq \bigoplus_{\nu \vdash k} f_\nu\,V_\nu \,,
\end{equation}
where $f_\nu$ denotes the number of standard Young tableaux of shape $\nu$.

\begin{example}{\rm
The $n^2$-dimensional $U_q(\gl_n)$-module $V^{\otimes 2}$ decomposes into
the $q$-symmetric squa\-re $V_{(2)}$ and the $q$-alternating square
$V_{(1,1)}$.
For $n=2$ this decomposition is described by the following diagram:
$$
0 \stackrel{e_1}{\longleftarrow} v_1\otimes v_1
\stackrel{f_1}{\longrightarrow}
v_1\otimes v_2 + q v_2\otimes v_1  \stackrel{f_1}{\longrightarrow}
(q+q^{-1})\, v_2\otimes v_2 \stackrel{f_1}{\longrightarrow} 0 \ \ \simeq
V_{(2)}
$$
$$
0 \stackrel{e_1}{\longleftarrow} v_2\otimes v_1 - q v_1\otimes v_2
\stackrel{f_1}{\longrightarrow} 0 \ \ \simeq V_{(1,1)}
$$
}
\end{example}

The algebra $\F_q[\Mat_n]$ defined in Section~\ref{INT} is also endowed with
a natural structure of $U_q(\gl_n)$-module via the action defined by
\begin{equation}
q^{\epsilon_i}\, t_{kl}  =  q^{\delta_{il}}\, t_{kl} \label{ACTD} ,\quad
e_i\,  t_{kl}   =  \delta_{i+1\,l}\,  t_{k\,l-1}   ,\quad
f_i\,  t_{kl}  =  \delta_{il}\,  t_{k\,l+1} ,
\end{equation}
and the Leibniz formulas
\begin{eqnarray}
q^{\epsilon_i} (PQ) & = & (q^{\epsilon_i} P)\, .\,( q^{\epsilon_i} Q) ,
\label{L1} \\
e_i (PQ) & = & (e_i P)\, .\, Q  + (q^{-h_i} P)\, .\, (e_i Q) , \label{L2}\\
f_i (PQ) & = & (f_i P)\, .\,( q^{h_i} Q) + P\, .\,( f_i Q) , \label{L3}
\end{eqnarray}
for $P,\ Q$ in $\F_q[\Mat_n]$. This provides a very convenient realization
of the irreducible modules $V_\lambda$ as natural subspaces of $\F_q[\Mat_n]$.
To describe it, we introduce some notations. We shall write $y_\lambda$ for
the unique Young tableau of shape and weight $\lambda$. This is the so-called
{\it Yamanouchi tableau} of shape $\lambda$. Let $\tau$ be any Young tableau of
shape $\lambda$. The quantum bitableau $(y_\lambda \,|\, \tau)$ will be
simply denoted by $(\tau )$ and will be called a {\it quantum tableau}. This is
a product of quantum minors taken on the first rows of the matrix $T$.
{\it Quantum tabloids} are defined similarly.
Finally, denote by ${\cal T}_\lambda$ the subspace of $\F_q[\Mat_n]$ spanned
by quantum tableaux $(\tau)$ of shape $\lambda$. Then one can show
\cite{LR,NYM}
the following $q$-analogue of a classical result of Deruyts (see \cite{Gr}).
\begin{theorem}\label{qDERUYTS}
The subspace ${\cal T}_\lambda$ is invariant under the action of $U_q(\gl_n)$
on $\F_q[\Mat_n]$, and is isomorphic as a $U_q(\gl_n)$-module to the simple
module $V_\lambda$.
\end{theorem}
The action of $U_q(\gl_n)$ on ${\cal T}_\lambda$ is computed by means of
the $q$-straightening formula. Namely, for column-shaped quantum tableaux
one checks easily that the action coincides with the one previously
described in Example~\ref{PEXT}.  For general quantum
tableaux we use Leibniz formulas (\ref{L1}) (\ref{L2}) (\ref{L3}), and when
necessary we use the $q$-straightening algorithm (\cf Section~\ref{qG/B})
for converting the quantum
tabloids of the right-hand side into a linear combination of quantum tableaux.
\begin{example}{\rm
We choose $n=3$ and $\lambda = (2,1)$. ${\cal T}_{(2,1)}$ is $8$-dimensional
and one has for instance
\begin{center}
\setlength{\unitlength}{0.01in}
\begingroup\makeatletter\ifx\SetFigFont\undefined
\def\x#1#2#3#4#5#6#7\relax{\def\x{#1#2#3#4#5#6}}%
\expandafter\x\fmtname xxxxxx\relax \def\y{splain}%
\ifx\x\y   
\gdef\SetFigFont#1#2#3{%
  \ifnum #1<17\tiny\else \ifnum #1<20\small\else
  \ifnum #1<24\normalsize\else \ifnum #1<29\large\else
  \ifnum #1<34\Large\else \ifnum #1<41\LARGE\else
     \huge\fi\fi\fi\fi\fi\fi
  \csname #3\endcsname}%
\else
\gdef\SetFigFont#1#2#3{\begingroup
  \count@#1\relax \ifnum 25<\count@\count@25\fi
  \def\x{\endgroup\@setsize\SetFigFont{#2pt}}%
  \expandafter\x
    \csname \romannumeral\the\count@ pt\expandafter\endcsname
    \csname @\romannumeral\the\count@ pt\endcsname
  \csname #3\endcsname}%
\fi
\fi\endgroup
\begin{picture}(415,55)(0,-10)
\path(35,0)(35,40)(15,40)
	(15,0)(35,0)
\path(55,0)(55,20)(35,20)
	(35,0)(55,0)
\path(115,0)(115,40)(95,40)
	(95,0)(115,0)
\path(135,0)(135,20)(115,20)
	(115,0)(135,0)
\path(195,0)(195,40)(175,40)
	(175,0)(195,0)
\path(215,0)(215,20)(195,20)
	(195,0)(215,0)
\path(315,0)(315,40)(295,40)
	(295,0)(315,0)
\path(335,0)(335,20)(315,20)
	(315,0)(335,0)
\path(395,0)(395,40)(375,40)
	(375,0)(395,0)
\path(415,0)(415,20)(395,20)
	(395,0)(415,0)
\put(-3,5){\makebox(0,0)[lb]{\smash{{{\SetFigFont{12}{14.4}{rm}$f_1$}}}}}
\put(70,5){\makebox(0,0)[lb]{\smash{{{\SetFigFont{12}{14.4}{rm}$=$}}}}}
\put(145,5){\makebox(0,0)[lb]{\smash{{{\SetFigFont{12}{14.4}{rm}$+\ q$}}}}}
\put(223,5){\makebox(0,0)[lb]{\smash{{{\SetFigFont{12}{14.4}{rm}$=(1+q^2)$}}}}}
\put(342,5){\makebox(0,0)[lb]{\smash{{{\SetFigFont{12}{14.4}{rm}$-q^3$}}}}}
\put(20,25){\makebox(0,0)[lb]{\smash{{{\SetFigFont{12}{14.4}{rm}3}}}}}
\put(20,5){\makebox(0,0)[lb]{\smash{{{\SetFigFont{12}{14.4}{rm}1}}}}}
\put(40,5){\makebox(0,0)[lb]{\smash{{{\SetFigFont{12}{14.4}{rm}1}}}}}
\put(100,25){\makebox(0,0)[lb]{\smash{{{\SetFigFont{12}{14.4}{rm}3}}}}}
\put(100,5){\makebox(0,0)[lb]{\smash{{{\SetFigFont{12}{14.4}{rm}1}}}}}
\put(120,5){\makebox(0,0)[lb]{\smash{{{\SetFigFont{12}{14.4}{rm}2}}}}}
\put(180,25){\makebox(0,0)[lb]{\smash{{{\SetFigFont{12}{14.4}{rm}3}}}}}
\put(180,5){\makebox(0,0)[lb]{\smash{{{\SetFigFont{12}{14.4}{rm}2}}}}}
\put(200,5){\makebox(0,0)[lb]{\smash{{{\SetFigFont{12}{14.4}{rm}1}}}}}
\put(300,25){\makebox(0,0)[lb]{\smash{{{\SetFigFont{12}{14.4}{rm}3}}}}}
\put(300,5){\makebox(0,0)[lb]{\smash{{{\SetFigFont{12}{14.4}{rm}1}}}}}
\put(320,5){\makebox(0,0)[lb]{\smash{{{\SetFigFont{12}{14.4}{rm}2}}}}}
\put(380,25){\makebox(0,0)[lb]{\smash{{{\SetFigFont{12}{14.4}{rm}2}}}}}
\put(380,5){\makebox(0,0)[lb]{\smash{{{\SetFigFont{12}{14.4}{rm}1}}}}}
\put(400,5){\makebox(0,0)[lb]{\smash{{{\SetFigFont{12}{14.4}{rm}3}}}}}
\end{picture}
\end{center}
In this example, the quantum tableaux $(\tau)$ have been written for short
$\tau$
(without brackets). This small abuse of notation will be used freely in the
sequel.
}\end{example}

This realization of $V_\lambda$ is, up to some minor changes
of convention, the same as the
one described in \cite{BKW} via the $q$-Young symmetrizers of the Hecke algebra
of type $A$.

We denote by $\F_q[GL_n/B]$ the subspace of $\F_q[\Mat_n]$ spanned by the
quantum tableaux. It follows from the $q$-straightening formula that this is
in fact a subalgebra of $\F_q[\Mat_n]$. It coincides for $q=1$ with the
ring of polynomial functions on the flag variety $GL_n/B$, hence the notation.
The quantum deformation $\F_q[GL_n/B]$ has been studied by Lakshmibai,
Reshetikhin
\cite{LR} and Taft, Towber \cite{TT}. As a $U_q(\gl_n)$-module it decomposes
into:
\begin{equation}
\F_q[GL_n/B] \simeq \bigoplus_{\ell (\lambda) \le n} V_\lambda \,.
\end{equation}

Returning to $\F_q[\Mat_n]$ and its linear basis formed by quantum bitableaux,
we note that the action of $U_q(\gl_n)$ defined by (\ref{ACTD})
involves only the column indices of the variables $t_{ij}$. The
defining relations (\ref{qGL1}) (\ref{qGL2}) (\ref{qGL3}) (\ref{qGL4}) being
invariant under transposition of the matrix $T$, we see that we have
another action of $U_q(\gl_n)$ given by
\begin{equation}
(q^{\epsilon_i})^\dagger t_{kl}  =  q^{\delta_{ik}}\, t_{kl} , \quad
\label{ACTG}
e_i^\dagger  t_{kl}  =  \delta_{i+1\,k}\,  t_{k-1\,l} , \quad
f_i^\dagger  t_{kl}  =  \delta_{ik}\,  t_{k+1\,l} ,
\end{equation}
(where the symbol $\dagger$ has been added to distinguish this action from
the previous one), and the Leibniz formulas (\ref{L1}) (\ref{L2}) (\ref{L3}).
These two actions obviously commute with each other, so that $\F_q[\Mat_n]$
is now endowed with the structure of a (left) bimodule over $U_q(\gl_n)$. The
quantum version of the Peter-Weyl theorem provides the decomposition \cite{NYM}
\begin{equation}
\F_q[\Mat_n] \simeq \bigoplus_{\ell (\lambda) \le n} V_\lambda \otimes
V_\lambda\,.
\end{equation}
Here, the irreducible bimodule $V_\lambda \otimes V_\lambda$ is generated by
applying
all possible products of lowering operators $f_i^\dagger,\,f_j$ to the highest
weight vector $(y_\lambda | y_\lambda)$.

We end this Section by noting that every $U_q(\gl_n)$-module $M$ can be
regarded
by restriction as a $U_q(\Sl_n)$-module (that we still denote by $M$).
In particular, the $V_\lambda$ are also irreducible under $U_q(\Sl_n)$. However
we point out that, as $U_q(\Sl_n)$-modules,
$$
V_\lambda \simeq V_\mu \quad \Longleftrightarrow
\quad \lambda_i -\lambda_{i+1} = \mu_i -\mu_{i+1}\,, \ i=1,\ldots ,n-1 \,.
$$
\begin{example} \label{Vl}{\rm
The $U_q(\Sl_2)$-modules $V_{(l)}$ will be very important in the sequel and
we describe them precisely. For $l\ge 0$, $V_{(l)}$ is a $(l+1)$-dimensional
vector space over $K(q)$ with basis $\{u_k,\ 0\le k \le l\}$, on which the
action
of $U_q(\Sl_2)$ is as follows:
$$
q^{h_1}\,u_k = q^{l-2k}\,u_k\,,\quad
e_1\,u_k = [l-k+1]\,u_{k-1}\,, \quad
f_1\,u_k = [k+1]\,u_{k+1}\,.
$$
In these formulas $[m]$ denotes the $q$-integer $(q^m-q^{-m})/(q-q^{-1})$,
and we understand $u_{-1}=u_{l+1}=0$. Setting $[m]!=[m][m-1]\cdots [1]$ and
$f_1^{(m)} = f_1^m/[m]!$, we see that the basis $\{u_k\}$ is characterized by
$u_k = f_1^{(k)}\,u_0$.
Also, we note that the weight spaces being one-dimensional, there is up to
normalization a unique basis of $V_{(l)}$ whose elements are weight vectors.
The basis $\{u_k\}$ may therefore be regarded as canonical. This
will provide the starting point for defining the crystal basis of a
$U_q(\gl_n)$-module.
}\end{example}

\section{Crystal bases} \label{BC}

It follows from relations (\ref{Uq2}) (\ref{Uq3}) (\ref{Uq4}) (\ref{Uq7}) that
for any $i = 1,\ldots ,n-1$, the subalgebra $U_i$ generated by $e_i,\ f_i,\
q^{h_i},\ q^{-h_i}$ is isomorphic to $U_q(\Sl_2)$. Hence a $U_q(\gl_n)$-module
$M$ can be regarded by restriction to $U_i$ as a $U_q(\Sl_2)$-module. We
shall assume from now on that the weight spaces $M_\mu$ are finite-dimensional,
that $M = \oplus_\mu M_\mu$, and that for any $i$, $M$ decomposes into a
direct sum of finite-dimensional $U_i$-modules. Such modules $M$ are said
to be {\it integrable}. It follows from the representation theory of
$U_q(\Sl_2)$
that for any $i$, the integrable module $M$ is a direct sum of irreducible
$U_i$-modules $V_{(l)}$.
\begin{example}\label{EXA1}{\rm Let $M$ denote the $U_q(\gl_3)$-module
$V_{(2,1)}$
in the realization given by Theorem~\ref{qDERUYTS}. As
a $U_1$-module, $M$ decomposes into $4$ irreducible components, as shown
by the following diagram:
\begin{center}
\setlength{\unitlength}{0.0098in}
\begingroup\makeatletter\ifx\SetFigFont\undefined
\def\x#1#2#3#4#5#6#7\relax{\def\x{#1#2#3#4#5#6}}%
\expandafter\x\fmtname xxxxxx\relax \def\y{splain}%
\ifx\x\y   
\gdef\SetFigFont#1#2#3{%
  \ifnum #1<17\tiny\else \ifnum #1<20\small\else
  \ifnum #1<24\normalsize\else \ifnum #1<29\large\else
  \ifnum #1<34\Large\else \ifnum #1<41\LARGE\else
     \huge\fi\fi\fi\fi\fi\fi
  \csname #3\endcsname}%
\else
\gdef\SetFigFont#1#2#3{\begingroup
  \count@#1\relax \ifnum 25<\count@\count@25\fi
  \def\x{\endgroup\@setsize\SetFigFont{#2pt}}%
  \expandafter\x
    \csname \romannumeral\the\count@ pt\expandafter\endcsname
    \csname @\romannumeral\the\count@ pt\endcsname
  \csname #3\endcsname}%
\fi
\fi\endgroup
\begin{picture}(671,295)(0,-10)
\path(115,240)(115,280)(95,280)
	(95,240)(115,240)
\path(135,240)(135,260)(115,260)
	(115,240)(135,240)
\path(255,240)(255,280)(235,280)
	(235,240)(255,240)
\path(275,240)(275,260)(255,260)
	(255,240)(275,240)
\path(115,160)(115,200)(95,200)
	(95,160)(115,160)
\path(135,160)(135,180)(115,180)
	(115,160)(135,160)
\path(315,160)(315,200)(295,200)
	(295,160)(315,160)
\path(335,160)(335,180)(315,180)
	(315,160)(335,160)
\path(415,160)(415,200)(395,200)
	(395,160)(415,160)
\path(435,160)(435,180)(415,180)
	(415,160)(435,160)
\path(555,160)(555,200)(535,200)
	(535,160)(555,160)
\path(575,160)(575,180)(555,180)
	(555,160)(575,160)
\path(115,80)(115,120)(95,120)
	(95,80)(115,80)
\path(135,80)(135,100)(115,100)
	(115,80)(135,80)
\path(115,0)(115,40)(95,40)
	(95,0)(115,0)
\path(135,0)(135,20)(115,20)
	(115,0)(135,0)
\path(255,0)(255,40)(235,40)
	(235,0)(255,0)
\path(275,0)(275,20)(255,20)
	(255,0)(275,0)
\path(155,250)(215,250)
\path(207.000,248.000)(215.000,250.000)(207.000,252.000)
\path(155,170)(215,170)
\path(207.000,168.000)(215.000,170.000)(207.000,172.000)
\path(155,90)(215,90)
\path(207.000,88.000)(215.000,90.000)(207.000,92.000)
\path(155,10)(215,10)
\path(207.000,8.000)(215.000,10.000)(207.000,12.000)
\path(455,170)(515,170)
\path(507.000,168.000)(515.000,170.000)(507.000,172.000)
\path(295,250)(355,250)
\path(347.000,248.000)(355.000,250.000)(347.000,252.000)
\path(575,170)(635,170)
\path(627.000,168.000)(635.000,170.000)(627.000,172.000)
\path(295,10)(355,10)
\path(347.000,8.000)(355.000,10.000)(347.000,12.000)
\path(75,250)(15,250)
\path(23.000,252.000)(15.000,250.000)(23.000,248.000)
\path(75,170)(15,170)
\path(23.000,172.000)(15.000,170.000)(23.000,168.000)
\path(75,90)(15,90)
\path(23.000,92.000)(15.000,90.000)(23.000,88.000)
\path(75,10)(15,10)
\path(23.000,12.000)(15.000,10.000)(23.000,8.000)
\put(100,245){\makebox(0,0)[lb]{\smash{{{\SetFigFont{12}{14.4}{rm}1}}}}}
\put(100,265){\makebox(0,0)[lb]{\smash{{{\SetFigFont{12}{14.4}{rm}2}}}}}
\put(120,245){\makebox(0,0)[lb]{\smash{{{\SetFigFont{12}{14.4}{rm}1}}}}}
\put(240,245){\makebox(0,0)[lb]{\smash{{{\SetFigFont{12}{14.4}{rm}1}}}}}
\put(240,265){\makebox(0,0)[lb]{\smash{{{\SetFigFont{12}{14.4}{rm}2}}}}}
\put(260,245){\makebox(0,0)[lb]{\smash{{{\SetFigFont{12}{14.4}{rm}2}}}}}
\put(100,165){\makebox(0,0)[lb]{\smash{{{\SetFigFont{12}{14.4}{rm}1}}}}}
\put(100,185){\makebox(0,0)[lb]{\smash{{{\SetFigFont{12}{14.4}{rm}3}}}}}
\put(120,165){\makebox(0,0)[lb]{\smash{{{\SetFigFont{12}{14.4}{rm}1}}}}}
\put(300,165){\makebox(0,0)[lb]{\smash{{{\SetFigFont{12}{14.4}{rm}1}}}}}
\put(300,185){\makebox(0,0)[lb]{\smash{{{\SetFigFont{12}{14.4}{rm}3}}}}}
\put(320,165){\makebox(0,0)[lb]{\smash{{{\SetFigFont{12}{14.4}{rm}2}}}}}
\put(400,165){\makebox(0,0)[lb]{\smash{{{\SetFigFont{12}{14.4}{rm}1}}}}}
\put(400,185){\makebox(0,0)[lb]{\smash{{{\SetFigFont{12}{14.4}{rm}2}}}}}
\put(420,165){\makebox(0,0)[lb]{\smash{{{\SetFigFont{12}{14.4}{rm}3}}}}}
\put(540,165){\makebox(0,0)[lb]{\smash{{{\SetFigFont{12}{14.4}{rm}2}}}}}
\put(560,165){\makebox(0,0)[lb]{\smash{{{\SetFigFont{12}{14.4}{rm}2}}}}}
\put(540,185){\makebox(0,0)[lb]{\smash{{{\SetFigFont{12}{14.4}{rm}3}}}}}
\put(100,85){\makebox(0,0)[lb]{\smash{{{\SetFigFont{12}{14.4}{rm}1}}}}}
\put(100,105){\makebox(0,0)[lb]{\smash{{{\SetFigFont{12}{14.4}{rm}2}}}}}
\put(120,85){\makebox(0,0)[lb]{\smash{{{\SetFigFont{12}{14.4}{rm}3}}}}}
\put(100,5){\makebox(0,0)[lb]{\smash{{{\SetFigFont{12}{14.4}{rm}1}}}}}
\put(100,25){\makebox(0,0)[lb]{\smash{{{\SetFigFont{12}{14.4}{rm}3}}}}}
\put(120,5){\makebox(0,0)[lb]{\smash{{{\SetFigFont{12}{14.4}{rm}3}}}}}
\put(240,5){\makebox(0,0)[lb]{\smash{{{\SetFigFont{12}{14.4}{rm}2}}}}}
\put(240,25){\makebox(0,0)[lb]{\smash{{{\SetFigFont{12}{14.4}{rm}3}}}}}
\put(260,5){\makebox(0,0)[lb]{\smash{{{\SetFigFont{12}{14.4}{rm}3}}}}}
\put(40,255){\makebox(0,0)[lb]{\smash{{{\SetFigFont{12}{14.4}{rm}$e_1$}}}}}
\put(40,175){\makebox(0,0)[lb]{\smash{{{\SetFigFont{12}{14.4}{rm}$e_1$}}}}}
\put(40,95){\makebox(0,0)[lb]{\smash{{{\SetFigFont{12}{14.4}{rm}$e_1$}}}}}
\put(40,15){\makebox(0,0)[lb]{\smash{{{\SetFigFont{12}{14.4}{rm}$e_1$}}}}}
\put(180,255){\makebox(0,0)[lb]{\smash{{{\SetFigFont{12}{14.4}{rm}$f_1$}}}}}
\put(320,255){\makebox(0,0)[lb]{\smash{{{\SetFigFont{12}{14.4}{rm}$f_1$}}}}}
\put(180,175){\makebox(0,0)[lb]{\smash{{{\SetFigFont{12}{14.4}{rm}$f_1$}}}}}
\put(462,175){\makebox(0,0)[lb]{\smash{{{\SetFigFont{12}{14.4}{rm}$f_1/[2]$}}}}}
\put(600,175){\makebox(0,0)[lb]{\smash{{{\SetFigFont{12}{14.4}{rm}$f_1$}}}}}
\put(180,95){\makebox(0,0)[lb]{\smash{{{\SetFigFont{12}{14.4}{rm}$f_1$}}}}}
\put(180,15){\makebox(0,0)[lb]{\smash{{{\SetFigFont{12}{14.4}{rm}$f_1$}}}}}
\put(320,15){\makebox(0,0)[lb]{\smash{{{\SetFigFont{12}{14.4}{rm}$f_1$}}}}}
\put(0,5){\makebox(0,0)[lb]{\smash{{{\SetFigFont{12}{14.4}{rm}0}}}}}
\put(0,85){\makebox(0,0)[lb]{\smash{{{\SetFigFont{12}{14.4}{rm}0}}}}}
\put(0,165){\makebox(0,0)[lb]{\smash{{{\SetFigFont{12}{14.4}{rm}0}}}}}
\put(0,245){\makebox(0,0)[lb]{\smash{{{\SetFigFont{12}{14.4}{rm}0}}}}}
\put(360,245){\makebox(0,0)[lb]{\smash{{{\SetFigFont{12}{14.4}{rm}0}}}}}
\put(640,165){\makebox(0,0)[lb]{\smash{{{\SetFigFont{12}{14.4}{rm}0}}}}}
\put(220,85){\makebox(0,0)[lb]{\smash{{{\SetFigFont{12}{14.4}{rm}0}}}}}
\put(360,5){\makebox(0,0)[lb]{\smash{{{\SetFigFont{12}{14.4}{rm}0}}}}}
\put(230,165){\makebox(0,0)[lb]{\smash{{{\SetFigFont{12}{14.4}{rm}$(1+q^2)$}}}}}
\put(350,165){\makebox(0,0)[lb]{\smash{{{\SetFigFont{12}{14.4}{rm}$-\ q^3$}}}}}
\end{picture}
\end{center}
On the other hand, as a $U_2$-module, $M$ decomposes into:
\begin{center}
\setlength{\unitlength}{0.0098in}
\begingroup\makeatletter\ifx\SetFigFont\undefined
\def\x#1#2#3#4#5#6#7\relax{\def\x{#1#2#3#4#5#6}}%
\expandafter\x\fmtname xxxxxx\relax \def\y{splain}%
\ifx\x\y   
\gdef\SetFigFont#1#2#3{%
  \ifnum #1<17\tiny\else \ifnum #1<20\small\else
  \ifnum #1<24\normalsize\else \ifnum #1<29\large\else
  \ifnum #1<34\Large\else \ifnum #1<41\LARGE\else
     \huge\fi\fi\fi\fi\fi\fi
  \csname #3\endcsname}%
\else
\gdef\SetFigFont#1#2#3{\begingroup
  \count@#1\relax \ifnum 25<\count@\count@25\fi
  \def\x{\endgroup\@setsize\SetFigFont{#2pt}}%
  \expandafter\x
    \csname \romannumeral\the\count@ pt\expandafter\endcsname
    \csname @\romannumeral\the\count@ pt\endcsname
  \csname #3\endcsname}%
\fi
\fi\endgroup
\begin{picture}(605,295)(0,-10)
\path(115,240)(115,280)(95,280)
	(95,240)(115,240)
\path(135,240)(135,260)(115,260)
	(115,240)(135,240)
\path(255,240)(255,280)(235,280)
	(235,240)(255,240)
\path(275,240)(275,260)(255,260)
	(255,240)(275,240)
\path(115,160)(115,200)(95,200)
	(95,160)(115,160)
\path(135,160)(135,180)(115,180)
	(115,160)(135,160)
\path(115,80)(115,120)(95,120)
	(95,80)(115,80)
\path(135,80)(135,100)(115,100)
	(115,80)(135,80)
\path(115,0)(115,40)(95,40)
	(95,0)(115,0)
\path(135,0)(135,20)(115,20)
	(115,0)(135,0)
\path(255,0)(255,40)(235,40)
	(235,0)(255,0)
\path(275,0)(275,20)(255,20)
	(255,0)(275,0)
\path(255,160)(255,200)(235,200)
	(235,160)(255,160)
\path(275,160)(275,180)(255,180)
	(255,160)(275,160)
\path(355,160)(355,200)(335,200)
	(335,160)(355,160)
\path(375,160)(375,180)(355,180)
	(355,160)(375,160)
\path(495,160)(495,200)(475,200)
	(475,160)(495,160)
\path(515,160)(515,180)(495,180)
	(495,160)(515,160)
\path(395,170)(455,170)
\path(447.000,168.000)(455.000,170.000)(447.000,172.000)
\path(540,170)(600,170)
\path(592.000,168.000)(600.000,170.000)(592.000,172.000)
\put(240,165){\makebox(0,0)[lb]{\smash{{{\SetFigFont{12}{14.4}{rm}1}}}}}
\put(340,165){\makebox(0,0)[lb]{\smash{{{\SetFigFont{12}{14.4}{rm}1}}}}}
\put(480,185){\makebox(0,0)[lb]{\smash{{{\SetFigFont{12}{14.4}{rm}3}}}}}
\put(402,175){\makebox(0,0)[lb]{\smash{{{\SetFigFont{12}{14.4}{rm}$f_2/[2]$}}}}}
\put(560,175){\makebox(0,0)[lb]{\smash{{{\SetFigFont{12}{14.4}{rm}$f_2$}}}}}
\put(605,165){\makebox(0,0)[lb]{\smash{{{\SetFigFont{12}{14.4}{rm}0}}}}}
\put(240,185){\makebox(0,0)[lb]{\smash{{{\SetFigFont{12}{14.4}{rm}2}}}}}
\put(260,165){\makebox(0,0)[lb]{\smash{{{\SetFigFont{12}{14.4}{rm}3}}}}}
\put(340,185){\makebox(0,0)[lb]{\smash{{{\SetFigFont{12}{14.4}{rm}3}}}}}
\put(360,165){\makebox(0,0)[lb]{\smash{{{\SetFigFont{12}{14.4}{rm}2}}}}}
\put(480,165){\makebox(0,0)[lb]{\smash{{{\SetFigFont{12}{14.4}{rm}1}}}}}
\put(500,165){\makebox(0,0)[lb]{\smash{{{\SetFigFont{12}{14.4}{rm}3}}}}}
\put(290,165){\makebox(0,0)[lb]{\smash{{{\SetFigFont{12}{14.4}{rm}$+\ q$}}}}}
\path(255,90)(315,90)
\path(307.000,88.000)(315.000,90.000)(307.000,92.000)
\put(280,95){\makebox(0,0)[lb]{\smash{{{\SetFigFont{12}{14.4}{rm}$f_2$}}}}}
\put(320,85){\makebox(0,0)[lb]{\smash{{{\SetFigFont{12}{14.4}{rm}0}}}}}
\path(155,250)(215,250)
\path(207.000,248.000)(215.000,250.000)(207.000,252.000)
\path(155,170)(215,170)
\path(207.000,168.000)(215.000,170.000)(207.000,172.000)
\path(155,10)(215,10)
\path(207.000,8.000)(215.000,10.000)(207.000,12.000)
\path(295,250)(355,250)
\path(347.000,248.000)(355.000,250.000)(347.000,252.000)
\path(295,10)(355,10)
\path(347.000,8.000)(355.000,10.000)(347.000,12.000)
\path(75,250)(15,250)
\path(23.000,252.000)(15.000,250.000)(23.000,248.000)
\path(75,170)(15,170)
\path(23.000,172.000)(15.000,170.000)(23.000,168.000)
\path(75,90)(15,90)
\path(23.000,92.000)(15.000,90.000)(23.000,88.000)
\path(75,10)(15,10)
\path(23.000,12.000)(15.000,10.000)(23.000,8.000)
\path(215,80)(215,120)(195,120)
	(195,80)(215,80)
\path(235,80)(235,100)(215,100)
	(215,80)(235,80)
\put(100,245){\makebox(0,0)[lb]{\smash{{{\SetFigFont{12}{14.4}{rm}1}}}}}
\put(100,265){\makebox(0,0)[lb]{\smash{{{\SetFigFont{12}{14.4}{rm}2}}}}}
\put(120,245){\makebox(0,0)[lb]{\smash{{{\SetFigFont{12}{14.4}{rm}1}}}}}
\put(240,245){\makebox(0,0)[lb]{\smash{{{\SetFigFont{12}{14.4}{rm}1}}}}}
\put(100,165){\makebox(0,0)[lb]{\smash{{{\SetFigFont{12}{14.4}{rm}1}}}}}
\put(100,85){\makebox(0,0)[lb]{\smash{{{\SetFigFont{12}{14.4}{rm}1}}}}}
\put(100,25){\makebox(0,0)[lb]{\smash{{{\SetFigFont{12}{14.4}{rm}3}}}}}
\put(240,5){\makebox(0,0)[lb]{\smash{{{\SetFigFont{12}{14.4}{rm}2}}}}}
\put(240,25){\makebox(0,0)[lb]{\smash{{{\SetFigFont{12}{14.4}{rm}3}}}}}
\put(260,5){\makebox(0,0)[lb]{\smash{{{\SetFigFont{12}{14.4}{rm}3}}}}}
\put(40,255){\makebox(0,0)[lb]{\smash{{{\SetFigFont{12}{14.4}{rm}$e_2$}}}}}
\put(40,175){\makebox(0,0)[lb]{\smash{{{\SetFigFont{12}{14.4}{rm}$e_2$}}}}}
\put(40,95){\makebox(0,0)[lb]{\smash{{{\SetFigFont{12}{14.4}{rm}$e_2$}}}}}
\put(40,15){\makebox(0,0)[lb]{\smash{{{\SetFigFont{12}{14.4}{rm}$e_2$}}}}}
\put(180,255){\makebox(0,0)[lb]{\smash{{{\SetFigFont{12}{14.4}{rm}$f_2$}}}}}
\put(320,255){\makebox(0,0)[lb]{\smash{{{\SetFigFont{12}{14.4}{rm}$f_2$}}}}}
\put(180,175){\makebox(0,0)[lb]{\smash{{{\SetFigFont{12}{14.4}{rm}$f_2$}}}}}
\put(180,15){\makebox(0,0)[lb]{\smash{{{\SetFigFont{12}{14.4}{rm}$f_2$}}}}}
\put(320,15){\makebox(0,0)[lb]{\smash{{{\SetFigFont{12}{14.4}{rm}$f_2$}}}}}
\put(0,5){\makebox(0,0)[lb]{\smash{{{\SetFigFont{12}{14.4}{rm}0}}}}}
\put(0,85){\makebox(0,0)[lb]{\smash{{{\SetFigFont{12}{14.4}{rm}0}}}}}
\put(0,165){\makebox(0,0)[lb]{\smash{{{\SetFigFont{12}{14.4}{rm}0}}}}}
\put(0,245){\makebox(0,0)[lb]{\smash{{{\SetFigFont{12}{14.4}{rm}0}}}}}
\put(360,245){\makebox(0,0)[lb]{\smash{{{\SetFigFont{12}{14.4}{rm}0}}}}}
\put(360,5){\makebox(0,0)[lb]{\smash{{{\SetFigFont{12}{14.4}{rm}0}}}}}
\put(240,265){\makebox(0,0)[lb]{\smash{{{\SetFigFont{12}{14.4}{rm}3}}}}}
\put(100,185){\makebox(0,0)[lb]{\smash{{{\SetFigFont{12}{14.4}{rm}2}}}}}
\put(120,165){\makebox(0,0)[lb]{\smash{{{\SetFigFont{12}{14.4}{rm}2}}}}}
\put(100,5){\makebox(0,0)[lb]{\smash{{{\SetFigFont{12}{14.4}{rm}2}}}}}
\put(120,5){\makebox(0,0)[lb]{\smash{{{\SetFigFont{12}{14.4}{rm}2}}}}}
\put(200,85){\makebox(0,0)[lb]{\smash{{{\SetFigFont{12}{14.4}{rm}1}}}}}
\put(150,85){\makebox(0,0)[lb]{\smash{{{\SetFigFont{12}{14.4}{rm}$-\ q$}}}}}
\put(260,245){\makebox(0,0)[lb]{\smash{{{\SetFigFont{12}{14.4}{rm}1}}}}}
\put(100,105){\makebox(0,0)[lb]{\smash{{{\SetFigFont{12}{14.4}{rm}3}}}}}
\put(120,85){\makebox(0,0)[lb]{\smash{{{\SetFigFont{12}{14.4}{rm}2}}}}}
\put(200,105){\makebox(0,0)[lb]{\smash{{{\SetFigFont{12}{14.4}{rm}2}}}}}
\put(220,85){\makebox(0,0)[lb]{\smash{{{\SetFigFont{12}{14.4}{rm}3}}}}}
\end{picture}
\end{center}
We observe that the $U_1$-decomposition leads to the basis
\begin{center}
\setlength{\unitlength}{0.0098in}
\begingroup\makeatletter\ifx\SetFigFont\undefined
\def\x#1#2#3#4#5#6#7\relax{\def\x{#1#2#3#4#5#6}}%
\expandafter\x\fmtname xxxxxx\relax \def\y{splain}%
\ifx\x\y   
\gdef\SetFigFont#1#2#3{%
  \ifnum #1<17\tiny\else \ifnum #1<20\small\else
  \ifnum #1<24\normalsize\else \ifnum #1<29\large\else
  \ifnum #1<34\Large\else \ifnum #1<41\LARGE\else
     \huge\fi\fi\fi\fi\fi\fi
  \csname #3\endcsname}%
\else
\gdef\SetFigFont#1#2#3{\begingroup
  \count@#1\relax \ifnum 25<\count@\count@25\fi
  \def\x{\endgroup\@setsize\SetFigFont{#2pt}}%
  \expandafter\x
    \csname \romannumeral\the\count@ pt\expandafter\endcsname
    \csname @\romannumeral\the\count@ pt\endcsname
  \csname #3\endcsname}%
\fi
\fi\endgroup
\begin{picture}(375,55)(0,-10)
\path(60,0)(60,40)(40,40)
	(40,0)(60,0)
\path(80,0)(80,20)(60,20)
	(60,0)(80,0)
\path(160,0)(160,40)(140,40)
	(140,0)(160,0)
\path(180,0)(180,20)(160,20)
	(160,0)(180,0)
\put(45,5){\makebox(0,0)[lb]{\smash{{{\SetFigFont{12}{14.4}{rm}1}}}}}
\put(145,5){\makebox(0,0)[lb]{\smash{{{\SetFigFont{12}{14.4}{rm}1}}}}}
\put(45,25){\makebox(0,0)[lb]{\smash{{{\SetFigFont{12}{14.4}{rm}2}}}}}
\put(65,5){\makebox(0,0)[lb]{\smash{{{\SetFigFont{12}{14.4}{rm}3}}}}}
\put(145,25){\makebox(0,0)[lb]{\smash{{{\SetFigFont{12}{14.4}{rm}3}}}}}
\put(165,5){\makebox(0,0)[lb]{\smash{{{\SetFigFont{12}{14.4}{rm}2}}}}}
\put(95,5){\makebox(0,0)[lb]{\smash{{{\SetFigFont{12}{14.4}{rm}$+\ q$}}}}}
\path(240,0)(240,40)(220,40)
	(220,0)(240,0)
\path(260,0)(260,20)(240,20)
	(240,0)(260,0)
\path(340,0)(340,40)(320,40)
	(320,0)(340,0)
\path(360,0)(360,20)(340,20)
	(340,0)(360,0)
\put(225,5){\makebox(0,0)[lb]{\smash{{{\SetFigFont{12}{14.4}{rm}1}}}}}
\put(325,5){\makebox(0,0)[lb]{\smash{{{\SetFigFont{12}{14.4}{rm}1}}}}}
\put(275,5){\makebox(0,0)[lb]{\smash{{{\SetFigFont{12}{14.4}{rm}$-\ q$}}}}}
\put(225,25){\makebox(0,0)[lb]{\smash{{{\SetFigFont{12}{14.4}{rm}3}}}}}
\put(245,5){\makebox(0,0)[lb]{\smash{{{\SetFigFont{12}{14.4}{rm}2}}}}}
\put(325,25){\makebox(0,0)[lb]{\smash{{{\SetFigFont{12}{14.4}{rm}2}}}}}
\put(345,5){\makebox(0,0)[lb]{\smash{{{\SetFigFont{12}{14.4}{rm}3}}}}}
\put(30,5){\makebox(0,0)[lb]{\smash{{{\SetFigFont{12}{14.4}{rm}$($}}}}}
\put(200,5){\makebox(0,0)[lb]{\smash{{{\SetFigFont{12}{14.4}{rm}$;$}}}}}
\put(365,5){\makebox(0,0)[lb]{\smash{{{\SetFigFont{12}{14.4}{rm}$)$}}}}}
\put(-10,5){\makebox(0,0)[lb]{\smash{{{\SetFigFont{12}{14.4}{rm}$B_2 =$}}}}}
\end{picture}
\end{center}
for the 2-dimensional weight space $M_{(1,1,1)}$, while the $U_2$-decomposition
leads to
\begin{center}
\setlength{\unitlength}{0.0098in}
\begingroup\makeatletter\ifx\SetFigFont\undefined
\def\x#1#2#3#4#5#6#7\relax{\def\x{#1#2#3#4#5#6}}%
\expandafter\x\fmtname xxxxxx\relax \def\y{splain}%
\ifx\x\y   
\gdef\SetFigFont#1#2#3{%
  \ifnum #1<17\tiny\else \ifnum #1<20\small\else
  \ifnum #1<24\normalsize\else \ifnum #1<29\large\else
  \ifnum #1<34\Large\else \ifnum #1<41\LARGE\else
     \huge\fi\fi\fi\fi\fi\fi
  \csname #3\endcsname}%
\else
\gdef\SetFigFont#1#2#3{\begingroup
  \count@#1\relax \ifnum 25<\count@\count@25\fi
  \def\x{\endgroup\@setsize\SetFigFont{#2pt}}%
  \expandafter\x
    \csname \romannumeral\the\count@ pt\expandafter\endcsname
    \csname @\romannumeral\the\count@ pt\endcsname
  \csname #3\endcsname}%
\fi
\fi\endgroup
\begin{picture}(375,55)(0,-10)
\path(60,0)(60,40)(40,40)
	(40,0)(60,0)
\path(80,0)(80,20)(60,20)
	(60,0)(80,0)
\path(160,0)(160,40)(140,40)
	(140,0)(160,0)
\path(180,0)(180,20)(160,20)
	(160,0)(180,0)
\put(45,5){\makebox(0,0)[lb]{\smash{{{\SetFigFont{12}{14.4}{rm}1}}}}}
\put(145,5){\makebox(0,0)[lb]{\smash{{{\SetFigFont{12}{14.4}{rm}1}}}}}
\put(45,25){\makebox(0,0)[lb]{\smash{{{\SetFigFont{12}{14.4}{rm}2}}}}}
\put(65,5){\makebox(0,0)[lb]{\smash{{{\SetFigFont{12}{14.4}{rm}3}}}}}
\put(145,25){\makebox(0,0)[lb]{\smash{{{\SetFigFont{12}{14.4}{rm}3}}}}}
\put(165,5){\makebox(0,0)[lb]{\smash{{{\SetFigFont{12}{14.4}{rm}2}}}}}
\put(95,5){\makebox(0,0)[lb]{\smash{{{\SetFigFont{12}{14.4}{rm}$+\ q$}}}}}
\path(240,0)(240,40)(220,40)
	(220,0)(240,0)
\path(260,0)(260,20)(240,20)
	(240,0)(260,0)
\path(340,0)(340,40)(320,40)
	(320,0)(340,0)
\path(360,0)(360,20)(340,20)
	(340,0)(360,0)
\put(225,5){\makebox(0,0)[lb]{\smash{{{\SetFigFont{12}{14.4}{rm}1}}}}}
\put(325,5){\makebox(0,0)[lb]{\smash{{{\SetFigFont{12}{14.4}{rm}1}}}}}
\put(275,5){\makebox(0,0)[lb]{\smash{{{\SetFigFont{12}{14.4}{rm}$-\ q$}}}}}
\put(225,25){\makebox(0,0)[lb]{\smash{{{\SetFigFont{12}{14.4}{rm}3}}}}}
\put(245,5){\makebox(0,0)[lb]{\smash{{{\SetFigFont{12}{14.4}{rm}2}}}}}
\put(325,25){\makebox(0,0)[lb]{\smash{{{\SetFigFont{12}{14.4}{rm}2}}}}}
\put(345,5){\makebox(0,0)[lb]{\smash{{{\SetFigFont{12}{14.4}{rm}3}}}}}
\put(30,5){\makebox(0,0)[lb]{\smash{{{\SetFigFont{12}{14.4}{rm}$($}}}}}
\put(200,5){\makebox(0,0)[lb]{\smash{{{\SetFigFont{12}{14.4}{rm}$;$}}}}}
\put(365,5){\makebox(0,0)[lb]{\smash{{{\SetFigFont{12}{14.4}{rm}$)$}}}}}
\put(-10,5){\makebox(0,0)[lb]{\smash{{{\SetFigFont{12}{14.4}{rm}$B_2 =$}}}}}
\end{picture}
\end{center}
These two bases are different and therefore one cannot find a basis $B$ of
weight vectors in $M$ compatible with both decompositions. However, as noted
by Kashiwara, `at $q=0$' the bases $B_1$ and $B_2$ coincide. The next
definitions
will allow us to state this in a more formal way.
}\end{example}

Consider the simple $U_q(\Sl_2)$-module $V_{(l)}$ with basis $\{u_k\}$
(\cf Example~\ref{Vl}).
Kashiwara \cite{Ka1,Ka2} introduces
the endomorphisms $\tilde e,\ \tilde f$ of $V_{(l)}$ defined by
$$
\tilde e \, u_k = u_{k-1} \, , \quad \tilde f \, u_k = u_{k+1} \, , \quad
k=0,\ldots ,l,
$$
where $u_{-1}= u_{l+1} = 0$. More generally, if $M$ is a direct sum of modules
$V_{(l)}$, that is, if there exists an isomorphism of $U_q(\Sl_2)$-modules
$\phi : M \stackrel{\sim}{\longrightarrow} \oplus V_{(l)}^{\oplus \alpha_l}$,
one defines endomorphisms $\tilde e,\ \tilde f$ of $M$ by means of $\phi$ in
the obvious way, and one checks easily that they do not depend on the choice
of $\phi$. In particular, if $M$ is an integrable $U_q(\gl_n)$-module,
regarding
$M$ as a $U_i$-module we define operators $\tilde e_i,\ \tilde f_i$ on $M$ for
$i = 1,\ldots ,n-1$.
\begin{example}\label{EXA2}{\rm We keep the notations of Example~\ref{EXA1}. We
have
\begin{center}
\setlength{\unitlength}{0.0098in}
\begingroup\makeatletter\ifx\SetFigFont\undefined
\def\x#1#2#3#4#5#6#7\relax{\def\x{#1#2#3#4#5#6}}%
\expandafter\x\fmtname xxxxxx\relax \def\y{splain}%
\ifx\x\y   
\gdef\SetFigFont#1#2#3{%
  \ifnum #1<17\tiny\else \ifnum #1<20\small\else
  \ifnum #1<24\normalsize\else \ifnum #1<29\large\else
  \ifnum #1<34\Large\else \ifnum #1<41\LARGE\else
     \huge\fi\fi\fi\fi\fi\fi
  \csname #3\endcsname}%
\else
\gdef\SetFigFont#1#2#3{\begingroup
  \count@#1\relax \ifnum 25<\count@\count@25\fi
  \def\x{\endgroup\@setsize\SetFigFont{#2pt}}%
  \expandafter\x
    \csname \romannumeral\the\count@ pt\expandafter\endcsname
    \csname @\romannumeral\the\count@ pt\endcsname
  \csname #3\endcsname}%
\fi
\fi\endgroup
\begin{picture}(265,115)(0,-10)
\path(20,60)(20,100)(0,100)
	(0,60)(20,60)
\path(40,60)(40,80)(20,80)
	(20,60)(40,60)
\path(105,60)(105,100)(85,100)
	(85,60)(105,60)
\path(125,60)(125,80)(105,80)
	(105,60)(125,60)
\path(245,60)(245,100)(225,100)
	(225,60)(245,60)
\path(265,60)(265,80)(245,80)
	(245,60)(265,60)
\path(145,70)(205,70)
\path(197.000,68.000)(205.000,70.000)(197.000,72.000)
\put(90,65){\makebox(0,0)[lb]{\smash{{{\SetFigFont{12}{14.4}{rm}1}}}}}
\put(230,85){\makebox(0,0)[lb]{\smash{{{\SetFigFont{12}{14.4}{rm}3}}}}}
\put(160,79){\makebox(0,0)[lb]{\smash{{{\SetFigFont{12}{14.4}{rm}$\tilde
f_2$}}}}}
\put(90,85){\makebox(0,0)[lb]{\smash{{{\SetFigFont{12}{14.4}{rm}3}}}}}
\put(110,65){\makebox(0,0)[lb]{\smash{{{\SetFigFont{12}{14.4}{rm}2}}}}}
\put(230,65){\makebox(0,0)[lb]{\smash{{{\SetFigFont{12}{14.4}{rm}1}}}}}
\put(250,65){\makebox(0,0)[lb]{\smash{{{\SetFigFont{12}{14.4}{rm}3}}}}}
\put(5,65){\makebox(0,0)[lb]{\smash{{{\SetFigFont{12}{14.4}{rm}1}}}}}
\put(5,85){\makebox(0,0)[lb]{\smash{{{\SetFigFont{12}{14.4}{rm}2}}}}}
\put(25,65){\makebox(0,0)[lb]{\smash{{{\SetFigFont{12}{14.4}{rm}3}}}}}
\put(55,65){\makebox(0,0)[lb]{\smash{{{\SetFigFont{12}{14.4}{rm}$+\ q$}}}}}
\path(150,10)(210,10)
\path(202.000,8.000)(210.000,10.000)(202.000,12.000)
\put(175,15){\makebox(0,0)[lb]{\smash{{{\SetFigFont{12}{14.4}{rm}$\tilde
f_2$}}}}}
\put(215,5){\makebox(0,0)[lb]{\smash{{{\SetFigFont{12}{14.4}{rm}0}}}}}
\path(110,0)(110,40)(90,40)
	(90,0)(110,0)
\path(130,0)(130,20)(110,20)
	(110,0)(130,0)
\put(95,5){\makebox(0,0)[lb]{\smash{{{\SetFigFont{12}{14.4}{rm}1}}}}}
\put(95,25){\makebox(0,0)[lb]{\smash{{{\SetFigFont{12}{14.4}{rm}2}}}}}
\put(115,5){\makebox(0,0)[lb]{\smash{{{\SetFigFont{12}{14.4}{rm}3}}}}}
\path(20,0)(20,40)(0,40)
	(0,0)(20,0)
\path(40,0)(40,20)(20,20)
	(20,0)(40,0)
\put(5,5){\makebox(0,0)[lb]{\smash{{{\SetFigFont{12}{14.4}{rm}1}}}}}
\put(55,5){\makebox(0,0)[lb]{\smash{{{\SetFigFont{12}{14.4}{rm}$-\ q$}}}}}
\put(5,25){\makebox(0,0)[lb]{\smash{{{\SetFigFont{12}{14.4}{rm}3}}}}}
\put(25,5){\makebox(0,0)[lb]{\smash{{{\SetFigFont{12}{14.4}{rm}2}}}}}
\end{picture}
\end{center}
Hence,
\begin{center}
\setlength{\unitlength}{0.0098in}
\begingroup\makeatletter\ifx\SetFigFont\undefined
\def\x#1#2#3#4#5#6#7\relax{\def\x{#1#2#3#4#5#6}}%
\expandafter\x\fmtname xxxxxx\relax \def\y{splain}%
\ifx\x\y   
\gdef\SetFigFont#1#2#3{%
  \ifnum #1<17\tiny\else \ifnum #1<20\small\else
  \ifnum #1<24\normalsize\else \ifnum #1<29\large\else
  \ifnum #1<34\Large\else \ifnum #1<41\LARGE\else
     \huge\fi\fi\fi\fi\fi\fi
  \csname #3\endcsname}%
\else
\gdef\SetFigFont#1#2#3{\begingroup
  \count@#1\relax \ifnum 25<\count@\count@25\fi
  \def\x{\endgroup\@setsize\SetFigFont{#2pt}}%
  \expandafter\x
    \csname \romannumeral\the\count@ pt\expandafter\endcsname
    \csname @\romannumeral\the\count@ pt\endcsname
  \csname #3\endcsname}%
\fi
\fi\endgroup
\begin{picture}(220,115)(0,-10)
\path(20,60)(20,100)(0,100)
	(0,60)(20,60)
\path(40,60)(40,80)(20,80)
	(20,60)(40,60)
\path(200,60)(200,100)(180,100)
	(180,60)(200,60)
\path(220,60)(220,80)(200,80)
	(200,60)(220,60)
\put(185,85){\makebox(0,0)[lb]{\smash{{{\SetFigFont{12}{14.4}{rm}3}}}}}
\put(185,65){\makebox(0,0)[lb]{\smash{{{\SetFigFont{12}{14.4}{rm}1}}}}}
\put(205,65){\makebox(0,0)[lb]{\smash{{{\SetFigFont{12}{14.4}{rm}3}}}}}
\path(200,0)(200,40)(180,40)
	(180,0)(200,0)
\path(220,0)(220,20)(200,20)
	(200,0)(220,0)
\put(185,25){\makebox(0,0)[lb]{\smash{{{\SetFigFont{12}{14.4}{rm}3}}}}}
\put(185,5){\makebox(0,0)[lb]{\smash{{{\SetFigFont{12}{14.4}{rm}1}}}}}
\put(205,5){\makebox(0,0)[lb]{\smash{{{\SetFigFont{12}{14.4}{rm}3}}}}}
\path(20,0)(20,40)(0,40)
	(0,0)(20,0)
\path(40,0)(40,20)(20,20)
	(20,0)(40,0)
\path(60,10)(120,10)
\path(112.000,8.000)(120.000,10.000)(112.000,12.000)
\path(60,70)(120,70)
\path(112.000,68.000)(120.000,70.000)(112.000,72.000)
\put(5,65){\makebox(0,0)[lb]{\smash{{{\SetFigFont{12}{14.4}{rm}1}}}}}
\put(5,85){\makebox(0,0)[lb]{\smash{{{\SetFigFont{12}{14.4}{rm}2}}}}}
\put(25,65){\makebox(0,0)[lb]{\smash{{{\SetFigFont{12}{14.4}{rm}3}}}}}
\put(5,5){\makebox(0,0)[lb]{\smash{{{\SetFigFont{12}{14.4}{rm}1}}}}}
\put(5,25){\makebox(0,0)[lb]{\smash{{{\SetFigFont{12}{14.4}{rm}3}}}}}
\put(25,5){\makebox(0,0)[lb]{\smash{{{\SetFigFont{12}{14.4}{rm}2}}}}}
\put(85,15){\makebox(0,0)[lb]{\smash{{{\SetFigFont{12}{14.4}{rm}$\tilde
f_2$}}}}}
\put(85,75){\makebox(0,0)[lb]{\smash{{{\SetFigFont{12}{14.4}{rm}$\tilde
f_2$}}}}}
\put(130,65){\makebox(0,0)[lb]{\smash{{{\SetFigFont{12}{14.4}{rm}$\displaystyle{1\over1+q^2}$}}}}}
\put(130,5){\makebox(0,0)[lb]{\smash{{{\SetFigFont{12}{14.4}{rm}$\displaystyle{q\over 1+q^2}$}}}}}
\end{picture}
\end{center}
}\end{example}

Since we want to let $q$ tend to 0, we introduce the subring $\A$ of $K(q)$
consisting of rational functions without pole at $q=0$. A {\it crystal lattice}
of $M$ is a free $\A$-module $L$ such that $M=K(q)\otimes_\A L,\
L = \oplus_\mu L_\mu$ where $L_\mu = L \cap M_\mu$, and
\begin{equation}
\tilde e_i L \subset L\,, \quad
\tilde f_i L \subset L\,, \quad i = 1,\ldots ,n-1\,.
\end{equation}
In other words, $L$ spans $M$ over $K(q)$, $L$ is compatible with the weight
space decomposition of $M$ and is stable under the operators $\tilde e_i$
and $\tilde f_i$. It follows that $\tilde e_i, \, \tilde f_i$ induce
endomorphisms
of the $K$-vector space $L/qL$ that we shall still denote by $\tilde e_i, \,
\tilde f_i$.
Now Kashiwara defines a {\it crystal basis} of $M$ (at $q=0$) to be a pair
$(L,B)$ where $L$ is a crystal lattice in $M$ and $B$ is a basis of $L/qL$ such
that $B=\sqcup B_\mu$ where $B_\mu = B \cap (L_\mu/qL_\mu)$, and
\begin{equation}
\tilde e_i B \subset B\sqcup \{0\}\,, \quad
\tilde f_i B \subset B\sqcup \{0\}\,, \quad i = 1,\ldots ,n-1\,,
\end{equation}
\begin{equation}
\tilde e_i v = u \quad \Longleftrightarrow \quad \tilde f_i u = v ,\, \quad u,v
\in B,\,
\quad i = 1,\ldots ,n-1\,.
\end{equation}
\begin{example}\label{EXA3}{\rm
We continue Examples~\ref{EXA1} and \ref{EXA2}. Denote by ${\cal B}$ the basis
of
quantum tableaux in $M$, and let $L$ be the ${\cal A}$-lattice in $M$ spanned
by
the elements of ${\cal B}$. Let $B$ be the projection of ${\cal B}$ in $L/qL$.
Then, Example~\ref{EXA1} shows that $(L,B)$ is a crystal basis of $M$.
}\end{example}

Kashiwara has proven the following existence and uniqueness result for crystal
bases \cite{Ka1,Ka2}.
\begin{theorem}\label{EXIST}
Any integrable $U_q(\gl_n)$-module $M$ has a crystal basis $(L,B)$. Moreover,
if $(L',B')$ is another crystal basis of $M$, then there exists a
$U_q(\gl_n)$-automorphism of $M$ sending $L$ on $L'$ and inducing an
isomorphism
of vector spaces from $L/qL$ to $L'/qL'$ which sends $B$ on $B'$. In
particular,
if $M= V_\lambda$ is irreducible, its crystal basis $(L(\lambda),B(\lambda))$
is
unique up to an overall scalar multiple. It is given by
\begin{equation}
L(\lambda)= \sum_{1\le i_1,i_2,\ldots,i_r\le n-1} {\cal A}
\, \tilde f_{i_1}\tilde f_{i_2}\cdots \tilde f_{i_r}\, u_\lambda ,
\end{equation}
\begin{equation}
B(\lambda) = \{
\tilde f_{i_1}\tilde f_{i_2}\cdots \tilde f_{i_r}\, u_\lambda\ \mod
qL(\lambda)\, |\,
1\le i_1,\ldots,i_r\le n\}\backslash\{0\} ,
\end{equation}
where $u_\lambda$ is a highest weight vector of $V_\lambda$.
\end{theorem}

It follows that one can associate to each integrable $U_q(\gl_n)$-module $M$
a well-defined colored graph $\Gamma(M)$ whose vertices are labelled by the
elements of $B$ and whose edges describe the action of the operators $\tilde
f_i$ :
$$
u \stackrel{i}{\longrightarrow} v \ \Longleftrightarrow \ \tilde f_iu = v\,.
$$
$\Gamma(M)$ is called the {\it crystal graph} of $M$.
\begin{example}\label{EXA4}{\rm The crystal graph of the  $U_q(\gl_3)$-module
$V_{(2,1)}$ is
readily deduced from Examples~\ref{EXA1}, \ref{EXA2} and \ref{EXA3}. It is
shown in
Figure~\ref{V21}.
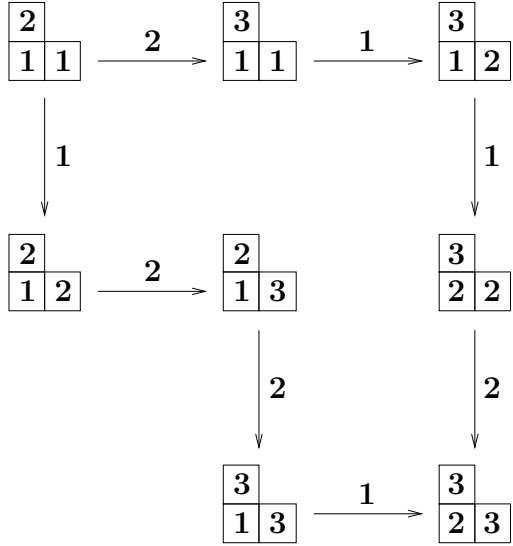
\begin{figure}[t]
\begin{center}
\setlength{\unitlength}{0.0085in}
\begingroup\makeatletter\ifx\SetFigFont\undefined
\def\x#1#2#3#4#5#6#7\relax{\def\x{#1#2#3#4#5#6}}%
\expandafter\x\fmtname xxxxxx\relax \def\y{splain}%
\ifx\x\y   
\gdef\SetFigFont#1#2#3{%
  \ifnum #1<17\tiny\else \ifnum #1<20\small\else
  \ifnum #1<24\normalsize\else \ifnum #1<29\large\else
  \ifnum #1<34\Large\else \ifnum #1<41\LARGE\else
     \huge\fi\fi\fi\fi\fi\fi
  \csname #3\endcsname}%
\else
\gdef\SetFigFont#1#2#3{\begingroup
  \count@#1\relax \ifnum 25<\count@\count@25\fi
  \def\x{\endgroup\@setsize\SetFigFont{#2pt}}%
  \expandafter\x
    \csname \romannumeral\the\count@ pt\expandafter\endcsname
    \csname @\romannumeral\the\count@ pt\endcsname
  \csname #3\endcsname}%
\fi
\fi\endgroup
\begin{picture}(309,349)(0,-10)
\path(22,309)(22,333)(0,333)
	(0,309)(22,309)
\path(22,285)(22,309)(0,309)
	(0,285)(22,285)
\path(44,285)(44,309)(22,309)
	(22,285)(44,285)
\path(22,167)(22,190)(0,190)
	(0,167)(22,167)
\path(44,143)(44,167)(22,167)
	(22,143)(44,143)
\path(22,143)(22,167)(0,167)
	(0,143)(22,143)
\path(154,309)(154,333)(132,333)
	(132,309)(154,309)
\path(154,285)(154,309)(132,309)
	(132,285)(154,285)
\path(177,285)(177,309)(154,309)
	(154,285)(177,285)
\path(287,309)(287,333)(265,333)
	(265,309)(287,309)
\path(287,285)(287,309)(265,309)
	(265,285)(287,285)
\path(309,285)(309,309)(287,309)
	(287,285)(309,285)
\path(154,167)(154,190)(132,190)
	(132,167)(154,167)
\path(154,143)(154,167)(132,167)
	(132,143)(154,143)
\path(177,143)(177,167)(154,167)
	(154,143)(177,143)
\path(287,167)(287,190)(265,190)
	(265,167)(287,167)
\path(287,143)(287,167)(265,167)
	(265,143)(287,143)
\path(309,143)(309,167)(287,167)
	(287,143)(309,143)
\path(154,24)(154,48)(132,48)
	(132,24)(154,24)
\path(154,0)(154,24)(132,24)
	(132,0)(154,0)
\path(177,0)(177,24)(154,24)
	(154,0)(177,0)
\path(287,24)(287,48)(265,48)
	(265,24)(287,24)
\path(287,0)(287,24)(265,24)
	(265,0)(287,0)
\path(309,0)(309,24)(287,24)
	(287,0)(309,0)
\path(55,297)(121,297)
\path(113.000,295.000)(121.000,297.000)(113.000,299.000)
\path(188,297)(254,297)
\path(246.000,295.000)(254.000,297.000)(246.000,299.000)
\path(287,274)(287,202)
\path(285.000,210.000)(287.000,202.000)(289.000,210.000)
\path(287,131)(287,59)
\path(285.000,67.000)(287.000,59.000)(289.000,67.000)
\path(188,18)(254,18)
\path(246.000,16.000)(254.000,18.000)(246.000,20.000)
\path(154,131)(154,59)
\path(152.000,67.000)(154.000,59.000)(156.000,67.000)
\path(55,155)(121,155)
\path(113.000,153.000)(121.000,155.000)(113.000,157.000)
\path(22,274)(22,202)
\path(20.000,210.000)(22.000,202.000)(24.000,210.000)
\put(6,315){\makebox(0,0)[lb]{\smash{{{\SetFigFont{12}{14.4}{bf}2}}}}}
\put(6,291){\makebox(0,0)[lb]{\smash{{{\SetFigFont{12}{14.4}{bf}1}}}}}
\put(28,291){\makebox(0,0)[lb]{\smash{{{\SetFigFont{12}{14.4}{bf}1}}}}}
\put(138,315){\makebox(0,0)[lb]{\smash{{{\SetFigFont{12}{14.4}{bf}3}}}}}
\put(138,291){\makebox(0,0)[lb]{\smash{{{\SetFigFont{12}{14.4}{bf}1}}}}}
\put(160,291){\makebox(0,0)[lb]{\smash{{{\SetFigFont{12}{14.4}{bf}1}}}}}
\put(270,315){\makebox(0,0)[lb]{\smash{{{\SetFigFont{12}{14.4}{bf}3}}}}}
\put(270,291){\makebox(0,0)[lb]{\smash{{{\SetFigFont{12}{14.4}{bf}1}}}}}
\put(292,291){\makebox(0,0)[lb]{\smash{{{\SetFigFont{12}{14.4}{bf}2}}}}}
\put(6,172){\makebox(0,0)[lb]{\smash{{{\SetFigFont{12}{14.4}{bf}2}}}}}
\put(6,149){\makebox(0,0)[lb]{\smash{{{\SetFigFont{12}{14.4}{bf}1}}}}}
\put(28,149){\makebox(0,0)[lb]{\smash{{{\SetFigFont{12}{14.4}{bf}2}}}}}
\put(138,149){\makebox(0,0)[lb]{\smash{{{\SetFigFont{12}{14.4}{bf}1}}}}}
\put(160,149){\makebox(0,0)[lb]{\smash{{{\SetFigFont{12}{14.4}{bf}3}}}}}
\put(270,172){\makebox(0,0)[lb]{\smash{{{\SetFigFont{12}{14.4}{bf}3}}}}}
\put(270,149){\makebox(0,0)[lb]{\smash{{{\SetFigFont{12}{14.4}{bf}2}}}}}
\put(292,149){\makebox(0,0)[lb]{\smash{{{\SetFigFont{12}{14.4}{bf}2}}}}}
\put(270,30){\makebox(0,0)[lb]{\smash{{{\SetFigFont{12}{14.4}{bf}3}}}}}
\put(270,6){\makebox(0,0)[lb]{\smash{{{\SetFigFont{12}{14.4}{bf}2}}}}}
\put(292,6){\makebox(0,0)[lb]{\smash{{{\SetFigFont{12}{14.4}{bf}3}}}}}
\put(138,30){\makebox(0,0)[lb]{\smash{{{\SetFigFont{12}{14.4}{bf}3}}}}}
\put(138,6){\makebox(0,0)[lb]{\smash{{{\SetFigFont{12}{14.4}{bf}1}}}}}
\put(160,6){\makebox(0,0)[lb]{\smash{{{\SetFigFont{12}{14.4}{bf}3}}}}}
\put(83,303){\makebox(0,0)[lb]{\smash{{{\SetFigFont{12}{14.4}{bf}2}}}}}
\put(215,303){\makebox(0,0)[lb]{\smash{{{\SetFigFont{12}{14.4}{bf}1}}}}}
\put(292,232){\makebox(0,0)[lb]{\smash{{{\SetFigFont{12}{14.4}{bf}1}}}}}
\put(292,89){\makebox(0,0)[lb]{\smash{{{\SetFigFont{12}{14.4}{bf}2}}}}}
\put(215,24){\makebox(0,0)[lb]{\smash{{{\SetFigFont{12}{14.4}{bf}1}}}}}
\put(160,89){\makebox(0,0)[lb]{\smash{{{\SetFigFont{12}{14.4}{bf}2}}}}}
\put(83,161){\makebox(0,0)[lb]{\smash{{{\SetFigFont{12}{14.4}{bf}2}}}}}
\put(28,232){\makebox(0,0)[lb]{\smash{{{\SetFigFont{12}{14.4}{bf}1}}}}}
\put(138,172){\makebox(0,0)[lb]{\smash{{{\SetFigFont{12}{14.4}{bf}2}}}}}
\end{picture}
\end{center}
\caption{Crystal graph of the $U_q(\gl_3)$-module $V_{(2,\,1)}$\label{V21}}
\end{figure}
}\end{example}

It is clear from this definition that the crystal graph $\Gamma(M)$ of the
direct
sum $M=M_1\oplus M_2$ of two $U_q(\gl_n)$-modules is the disjoint union of
$\Gamma(M_1)$ and $\Gamma(M_2)$. It follows from the complete reducibility
of $M$ that
the connected components of $\Gamma(M)$ are the crystal graphs of the
irreducible
components of $M$.
Note also that if one restricts the crystal graph $\Gamma(M)$ to its edges of
colour $i$, one obtains a decomposition of this graph into {\it strings} of
colour $i$ corresponding to the $U_i$-decomposition of $M$. For a vertex $v$
of $\Gamma(M)$, we shall denote by $\epsilon_i(v)$ (\resp $\phi_i(v)$)
the distance of $v$ to the origin (\resp end) of its string of colour $i$,
that is,
$$
\epsilon_i(v) = {\rm max}\{k \,|\, \tilde e_i^k v \not = 0\} \,, \quad
\phi_i(v) = {\rm max}\{k \,|\, \tilde f_i^k v \not = 0\} \,.
$$
The integers $\epsilon_i(v)$, $\phi_i(v)$ give information on the geometry of
the graph $\Gamma(M)$ around the vertex $v$. They seem to be very significant
from a representation theoretical point of view. Indeed, both the
Littlewood-Richardson multiplicities $c_{\lambda\,\mu}^\nu$ and the $q$-weight
multiplicities $K_{\lambda\,\mu}(q)$ can be computed in a simple way from
the integers $\epsilon_i(v)$, $\phi_i(v)$ attached to the crystal graph
$\Gamma(V_\lambda)$
\cite{BZ,LLT}.

One of the nicest properties of crystal bases is that they behave well under
the tensor product operation. We shall first consider an example, from which
Kashiwara deduces by induction the general description of the crystal basis
of a tensor product \cite{Ka1}.
\begin{example}{\rm
We slightly modify the notations of Example~\ref{Vl} and write
$\{u_k^{(l)},\,  k = 0,\ldots , l\}$ for the canonical basis of the
$U_q(\Sl_2)$-module
$V_{(l)}$. We recall that for convenience we set
$u_{-1}^{(l)} = u_{l+1}^{(l)} = 0$.
A first basis of the  tensor product $V_{(1)} \otimes V_{(l)}$ is
$\{u_j^{(1)}\otimes u_k^{(l)},\, j =0,1,\ k=0,\ldots , l\}$. We define
\begin{eqnarray}
w_k & = & u_0^{(1)} \otimes u_k^{(l)} + q^{l-k+1} u_1^{(1)}\otimes
u_{k-1}^{(l)}\,,
\quad k = 0,\ldots , l+1 \,, \label{wk}\\
z_k & = & q^{l-k-1}\, [l-k]\, u_1^{(1)} \otimes u_k^{(l)} -
q^l\, [k+1]\, u_0^{(1)} \otimes u_{k+1}^{(l)} \,,
\quad k = 0,\ldots , l-1\,. \label{zk}
\end{eqnarray}
It is straightforward to check that $e_1 w_0 = e_1 z_0 = 0,\ f_1 w_{l+1} = f_1
z_{l-1} = 0$,
and
$$
f_1 w_k = [k+1]\, w_{k+1}\,, \quad f_1 z_k = [k+1]\, z_{k+1} \,.
$$
Hence $\{w_k\}$ and $\{z_k\}$ span submodules of $V_{(1)} \otimes V_{(l)}$
isomorphic
to $V_{(l+1)}$ and $V_{(l-1)}$ respectively, and are the canonical bases of
these
irreducible representations of $U_q(\Sl_2)$. Therefore, the ${\cal A}$-lattice
$L$
spanned by $\{w_k\}$, $\{z_k\}$ is a crystal lattice and if we define
$B= \{w_k\, \mod qL\} \sqcup  \{z_k\, \mod qL\}$, then $(B,L)$ is a crystal
basis of
$V_{(1)} \otimes V_{(l)}$. On the other hand, equations (\ref{wk}) (\ref{zk})
show
that $L$ coincides with the tensor product $L(1)\otimes L(l)$ of the crystal
lattices of $V_{(1)}$ and $V_{(l)}$, and
that
$$
w_k \equiv u_0^{(1)}\otimes u_k^{(l)}\ \mod qL\,, \ k=0,\ldots ,l \,, \quad
w_{l+1} \equiv u_1^{(1)}\otimes u_l^{(l)}\ \mod qL \,.
$$
$$
z_k \equiv u_1^{(1)} \otimes u_k^{(l)}\ \mod qL\,,\quad k=0,\ldots ,l-1 \,.
$$
Thus, $(L,B)$ coincides with $(L(1)\otimes L(l), B(1)\otimes B(l))$, where we
have denoted
by $B(1)\otimes B(l)$ the basis $\{u_j^{(1)}\otimes u_k^{(l)},\,\mod qL\,,\
j =0,1,\ k=0,\ldots , l\}$. Finally, the action of $\tilde f_1$ on $B(1)\otimes
B(l)$
is described by the crystal graph:
$$
\matrix{
u_0^{(1)}\otimes u_0^{(l)} & \rightarrow & u_0^{(1)}\otimes u_1^{(l)} &
\rightarrow &
\cdots & \rightarrow & u_0^{(1)}\otimes u_{l-1}^{(l)} & \rightarrow &
u_0^{(1)}\otimes u_l^{(l)} \cr
& & & & & & & & \downarrow \cr
u_1^{(1)}\otimes u_0^{(l)} & \rightarrow & u_1^{(1)}\otimes u_1^{(l)} &
\rightarrow &
\cdots & \rightarrow & u_1^{(1)}\otimes u_{l-1}^{(l)} &  &
u_1^{(1)}\otimes u_l^{(l)} \cr
}
$$
}\end{example}

More generally, we have the following property \cite{Ka1}.
\begin{theorem}\label{TP}
Let $(L_1,B_1)$ and $(L_2,B_2)$ be crystal bases of integrable
$U_q(\gl_n)$-modules
$M_1$ and $M_2$. Let $B_1\otimes B_2$ denote the basis $\{ u\otimes v,\ u\in
B_1,\, v\in B_2\}$
of $(L_1/qL_1 )\otimes (L_2/qL_2)$ (which is isomorphic to
$(L_1\otimes L_2)/q(L_1\otimes L_2)$). Then, $(L_1\otimes L_2,\,B_1\otimes
B_2)$ is a
crystal basis of $M_1\otimes M_2$, the action of $\tilde e_i$, $\tilde f_i$ on
$B_1\otimes B_2$ being given by
\begin{eqnarray}
\tilde f_i (u\otimes v) & = &
\left\{\matrix{
u\otimes \tilde f_i v  & {\it if}& \epsilon_i(u) < \phi_i(v) \cr
\tilde f_i u\otimes v  & {\it if}& \epsilon_i(u) \ge \phi_i(v)
}\right. \,, \\
\tilde e_i (u\otimes v) & = &
\left\{\matrix{
\tilde e_i u\otimes v  & {\it if}& \epsilon_i(u) > \phi_i(v) \cr
u\otimes \tilde e_iv  & {\it if}& \epsilon_i(u) \le \phi_i(v)
}\right. \,.
\end{eqnarray}
\end{theorem}

Theorem~\ref{TP} enables one to describe the crystal graph of $V_{(1)}^{\otimes
m}$
for any $m$, and to deduce from that the description of the crystal graph
$\Gamma_\lambda$
of the simple $U_q(\gl_n)$-module $V_\lambda$. For convenience, we shall
identify
the tensor algebra $T(V_{(1)})$ with the free associative algebra $K(q)\<A\>$
over
the alphabet $A=\{1,\ldots , n\}$ via the isomorphism $v_i \mapsto i$, where
$\{v_i\}$
is the canonical basis of $V_{(1)}$ defined in Example~\ref{BAS}. Accordingly,
the crystal lattice $L$ spanned by the monomials in the $v_i$ is identified
with
${\cal A}\<A\>$ and the $K$-vector space $L/qL$ with $K\<A\>$. Define linear
operators $\hat e_i$, $\hat f_i$ on $K\<A\>$ in the following way.
Consider first the case of a two-letter alphabet $A=\{i,i{+}1\}$. Let
$w = x_1\cdots x_m$ be
a word on $A$. Delete every factor $((i{+}1)\,i)$ of $w$. The remaining letters
constitute a subword $w_1$ of $w$. Then, delete again every factor
$((i{+}1)\,i)$ of $w_1$. There remains a subword $w_2$. Continue this procedure
until it stops, leaving a word $w_k=x_{j_1}\cdots x_{j_{r+s}}$ of the form
$w_k = i^r (i{+}1)^s$.
The image of $w_k$ under $\hat e_i,\, \hat f_i$ is the word $y_{j_1}\cdots
y_{j_{r+s}}$
given by
$$
\hat e_i(i^r (i{+}1)^s)
=\left\{ \matrix{ i^{r+1} (i{+}1)^{s-1} & (s\ge 1) \cr
0               & (s=0)}\right. ,\ \quad
\hat f_i(i^r (i{+}1)^s)=
\left\{ \matrix{ i^{r-1} (i{+}1)^{s+1} & (r\ge 1) \cr
0               & (r=0)}\right. \ .
$$
The image of the initial word $w$ is then $w' = y_1\cdots y_m$, where
$y_k=x_k$ for $k \not \in \{j_1\,,\ldots ,\,j_{r+s}\}$.
For instance, if
$$w=(2\ 1)\ 1\ 1\ 2\ (2\ 1)\ 1\ 1\ 1\ 2\,,$$
we shall have
$$w_1=.\ .\ 1\ 1\ (2\ . \ . \ 1)\ 1\ 1\ 2\,,$$
$$w_2=.\ .\ \ 1\ 1\ .\ .\ .\ .\ \ 1\ 1\ 2\,.$$
Thus, $$\hat e_1(w)=\ 2\ 1{\bf \ 1\ 1}\ 2\ 2\ 1\ 1{\bf \ 1\ 1\ 1}\,,$$
$$\hat f_1(w)=2\ 1{\bf \ 1\ 1}\ 2\ 2\ 1\ 1\ {\bf 1\ 2\ 2}\,,$$
where the letters printed in bold type are those of the image of the subword
$w_2$. Finally, in the general case, the action of the operators
$\hat e_i,\, \hat f_i$ on $w$ is defined by the previous
rules applied to the subword consisting of the letters $i,\,i{+}1$,
the remaining letters being unchanged.

The operators $\hat e_i$, $\hat f_i$ have been considered in \cite{LS2}, where
they were used as building blocks for defining noncommutative analogues of
Demazure symmetrization operators. It follows from Theorem~\ref{TP} that they
coincide in the above identification with the endomorphisms $\tilde e_i$,
$\tilde f_i$ on the $K$-space $L/qL$ \cite{KN}.

It is straightforward to deduce from the definition of $\hat e_i$, $\hat f_i$
the
following compatibility properties with the Robinson-Schensted correspondence:
\begin{description}
\item[{\it (a)}] for any word $w$ on $A$ such that $\hat e_i w \not = 0$,
we have $Q(\hat e_i w) = Q(w)$,
\item[{\it (b)}] for any pair of words $w$, $u$ such that $P(w) = P(u)$ and
$\hat e_i w \not = 0$, we have $P(\hat e_i w) = P(\hat e_i u)$.
\end{description}
In other words, the operators $\hat e_i$, $\hat f_i$ do not change the
insertion
tableau and are compatible with the plactic equivalence. Moreover, the words
$y$
such that $\hat e_i y = 0$ for any $i$ are characterized by the Yamanouchi
property:
each right factor of $y$ contains at least as many letters $i$ than $i+1$, and
this
for any $i$. It is well-known that for any standard Young tableau $\tau$ there
exists
a unique Yamanouchi word $y$ such that $Q(y)=\tau$. This yields the following
cristallization of (\ref{TS}).
\begin{theorem}\label{TCG}
The crystal graph of the $U_q(\gl_n)$-module $V_{(1)}^{\otimes m}$ is the
colored graph
whose vertices are the words of length $m$ over $A$, and whose edges are given
by:
$$
w \stackrel{i}{\longrightarrow} u \quad \Longleftrightarrow \quad \hat f_iw =
u\,.
$$
The connected components $\Gamma_\tau$ of $\Gamma(V_{(1)}^{\otimes m})$ are
parametrized by the set of Young tableaux $\tau$ of weight $(1^m)$. The
vertices
of $\Gamma_\tau$ are those words $w$ which satisfy $Q(w)= \tau$. Moreover, if
$\lambda$ denotes the shape of $\tau$, then $\Gamma_\tau$ is isomorphic to
the crystal graph $\Gamma(V_\lambda)$.
\end{theorem}

Replacing the vertices $w$ of $\Gamma_\tau$ by their associated Young tableaux
$P(w)$, we obtain a labelling of $\Gamma(V_\lambda)$ by the set of Young
tableaux
of shape $\lambda$, as shown in Example~\ref{EXA4}. It follows from
property~{\it (b)}
above that this labelling does not depend on the particular choice of $\tau$
among the standard Young tableaux of shape $\lambda$. We end this section by
showing
a less trivial example of crystal graph, which is computed easily using
the previous description of $\hat f_i$.

\begin{example}{\rm The crystal graph of the $U_q(\gl_4)$-module $V_{(2,2)}$ is
shown
in Figure~\ref{V22}.
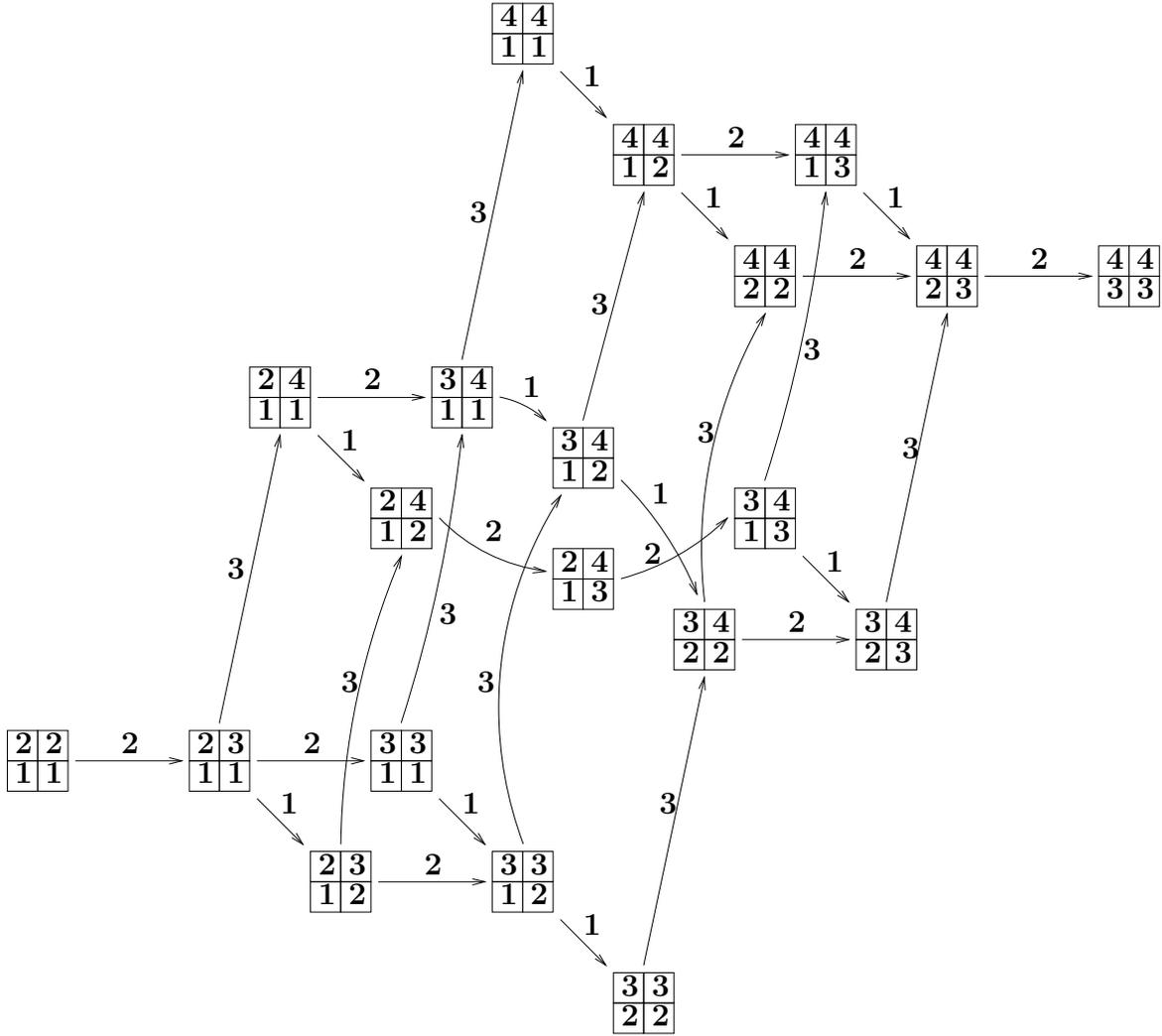
\begin{figure}[t]
\begin{center}
\setlength{\unitlength}{0.008in}
\begingroup\makeatletter\ifx\SetFigFont\undefined
\def\x#1#2#3#4#5#6#7\relax{\def\x{#1#2#3#4#5#6}}%
\expandafter\x\fmtname xxxxxx\relax \def\y{splain}%
\ifx\x\y   
\gdef\SetFigFont#1#2#3{%
  \ifnum #1<17\tiny\else \ifnum #1<20\small\else
  \ifnum #1<24\normalsize\else \ifnum #1<29\large\else
  \ifnum #1<34\Large\else \ifnum #1<41\LARGE\else
     \huge\fi\fi\fi\fi\fi\fi
  \csname #3\endcsname}%
\else
\gdef\SetFigFont#1#2#3{\begingroup
  \count@#1\relax \ifnum 25<\count@\count@25\fi
  \def\x{\endgroup\@setsize\SetFigFont{#2pt}}%
  \expandafter\x
    \csname \romannumeral\the\count@ pt\expandafter\endcsname
    \csname @\romannumeral\the\count@ pt\endcsname
  \csname #3\endcsname}%
\fi
\fi\endgroup
\begin{picture}(760,697)(0,-10)
\path(180,420)(180,440)(160,440)
	(160,420)(180,420)
\path(200,420)(200,440)(180,440)
	(180,420)(200,420)
\path(200,400)(200,420)(180,420)
	(180,400)(200,400)
\path(180,400)(180,420)(160,420)
	(160,400)(180,400)
\path(260,340)(260,360)(240,360)
	(240,340)(260,340)
\path(280,340)(280,360)(260,360)
	(260,340)(280,340)
\path(280,320)(280,340)(260,340)
	(260,320)(280,320)
\path(260,320)(260,340)(240,340)
	(240,320)(260,320)
\path(300,420)(300,440)(280,440)
	(280,420)(300,420)
\path(320,420)(320,440)(300,440)
	(300,420)(320,420)
\path(320,400)(320,420)(300,420)
	(300,400)(320,400)
\path(300,400)(300,420)(280,420)
	(280,400)(300,400)
\path(460,260)(460,280)(440,280)
	(440,260)(460,260)
\path(480,260)(480,280)(460,280)
	(460,260)(480,260)
\path(480,240)(480,260)(460,260)
	(460,240)(480,240)
\path(460,240)(460,260)(440,260)
	(440,240)(460,240)
\path(500,340)(500,360)(480,360)
	(480,340)(500,340)
\path(520,340)(520,360)(500,360)
	(500,340)(520,340)
\path(520,320)(520,340)(500,340)
	(500,320)(520,320)
\path(500,320)(500,340)(480,340)
	(480,320)(500,320)
\path(580,260)(580,280)(560,280)
	(560,260)(580,260)
\path(600,260)(600,280)(580,280)
	(580,260)(600,260)
\path(600,240)(600,260)(580,260)
	(580,240)(600,240)
\path(580,240)(580,260)(560,260)
	(560,240)(580,240)
\path(380,380)(380,400)(360,400)
	(360,380)(380,380)
\path(400,380)(400,400)(380,400)
	(380,380)(400,380)
\path(400,360)(400,380)(380,380)
	(380,360)(400,360)
\path(380,360)(380,380)(360,380)
	(360,360)(380,360)
\path(380,300)(380,320)(360,320)
	(360,300)(380,300)
\path(400,300)(400,320)(380,320)
	(380,300)(400,300)
\path(400,280)(400,300)(380,300)
	(380,280)(400,280)
\path(380,280)(380,300)(360,300)
	(360,280)(380,280)
\path(420,20)(420,40)(400,40)
	(400,20)(420,20)
\path(440,20)(440,40)(420,40)
	(420,20)(440,20)
\path(420,0)(420,20)(400,20)
	(400,0)(420,0)
\path(440,0)(440,20)(420,20)
	(420,0)(440,0)
\path(340,100)(340,120)(320,120)
	(320,100)(340,100)
\path(360,100)(360,120)(340,120)
	(340,100)(360,100)
\path(360,80)(360,100)(340,100)
	(340,80)(360,80)
\path(340,80)(340,100)(320,100)
	(320,80)(340,80)
\path(260,180)(260,200)(240,200)
	(240,180)(260,180)
\path(280,180)(280,200)(260,200)
	(260,180)(280,180)
\path(280,160)(280,180)(260,180)
	(260,160)(280,160)
\path(260,160)(260,180)(240,180)
	(240,160)(260,160)
\path(220,100)(220,120)(200,120)
	(200,100)(220,100)
\path(240,100)(240,120)(220,120)
	(220,100)(240,100)
\path(240,80)(240,100)(220,100)
	(220,80)(240,80)
\path(220,80)(220,100)(200,100)
	(200,80)(220,80)
\path(140,180)(140,200)(120,200)
	(120,180)(140,180)
\path(160,180)(160,200)(140,200)
	(140,180)(160,180)
\path(160,160)(160,180)(140,180)
	(140,160)(160,160)
\path(140,160)(140,180)(120,180)
	(120,160)(140,160)
\path(20,180)(20,200)(0,200)
	(0,180)(20,180)
\path(40,180)(40,200)(20,200)
	(20,180)(40,180)
\path(40,160)(40,180)(20,180)
	(20,160)(40,160)
\path(20,160)(20,180)(0,180)
	(0,160)(20,160)
\path(45,180)(115,180)
\path(107.000,178.000)(115.000,180.000)(107.000,182.000)
\path(165,180)(235,180)
\path(227.000,178.000)(235.000,180.000)(227.000,182.000)
\path(245,100)(315,100)
\path(307.000,98.000)(315.000,100.000)(307.000,102.000)
\path(285,155)(315,125)
\path(307.929,129.243)(315.000,125.000)(310.757,132.071)
\path(165,155)(195,125)
\path(187.929,129.243)(195.000,125.000)(190.757,132.071)
\path(365,75)(395,45)
\path(387.929,49.243)(395.000,45.000)(390.757,52.071)
\path(340,660)(340,680)(320,680)
	(320,660)(340,660)
\path(360,660)(360,680)(340,680)
	(340,660)(360,660)
\path(360,640)(360,660)(340,660)
	(340,640)(360,640)
\path(340,640)(340,660)(320,660)
	(320,640)(340,640)
\path(420,580)(420,600)(400,600)
	(400,580)(420,580)
\path(440,580)(440,600)(420,600)
	(420,580)(440,580)
\path(440,560)(440,580)(420,580)
	(420,560)(440,560)
\path(420,560)(420,580)(400,580)
	(400,560)(420,560)
\path(540,580)(540,600)(520,600)
	(520,580)(540,580)
\path(560,580)(560,600)(540,600)
	(540,580)(560,580)
\path(560,560)(560,580)(540,580)
	(540,560)(560,560)
\path(540,560)(540,580)(520,580)
	(520,560)(540,560)
\path(500,500)(500,520)(480,520)
	(480,500)(500,500)
\path(520,500)(520,520)(500,520)
	(500,500)(520,500)
\path(520,480)(520,500)(500,500)
	(500,480)(520,480)
\path(500,480)(500,500)(480,500)
	(480,480)(500,480)
\path(620,500)(620,520)(600,520)
	(600,500)(620,500)
\path(640,500)(640,520)(620,520)
	(620,500)(640,500)
\path(640,480)(640,500)(620,500)
	(620,480)(640,480)
\path(620,480)(620,500)(600,500)
	(600,480)(620,480)
\path(740,500)(740,520)(720,520)
	(720,500)(740,500)
\path(760,500)(760,520)(740,520)
	(740,500)(760,500)
\path(760,480)(760,500)(740,500)
	(740,480)(760,480)
\path(740,480)(740,500)(720,500)
	(720,480)(740,480)
\path(445,580)(515,580)
\path(507.000,578.000)(515.000,580.000)(507.000,582.000)
\path(525,500)(595,500)
\path(587.000,498.000)(595.000,500.000)(587.000,502.000)
\path(645,500)(715,500)
\path(707.000,498.000)(715.000,500.000)(707.000,502.000)
\path(445,555)(475,525)
\path(467.929,529.243)(475.000,525.000)(470.757,532.071)
\path(565,555)(595,525)
\path(587.929,529.243)(595.000,525.000)(590.757,532.071)
\path(365,635)(395,605)
\path(387.929,609.243)(395.000,605.000)(390.757,612.071)
\path(450.105,296.636)(455.000,290.000)(453.803,298.159)
\put(242.500,202.500){\arc{459.619}{5.4978}{5.8926}}
\path(347.704,408.844)(355.000,405.000)(350.372,411.825)
\put(312.500,357.500){\arc{127.475}{4.9098}{5.4423}}
\path(346.816,303.986)(355.000,305.000)(347.302,307.956)
\put(369.167,420.833){\arc{233.393}{1.6925}{2.3764}}
\path(470.689,332.971)(475.000,340.000)(467.888,335.826)
\put(363.750,453.438){\arc{317.772}{0.7951}{1.3082}}
\path(362.352,347.191)(365.000,355.000)(358.988,349.356)
\put(582.500,215.000){\arc{517.325}{2.7862}{3.7135}}
\path(301.155,386.835)(300.000,395.000)(297.177,387.252)
\put(-638.333,493.333){\arc{1886.943}{0.1044}{0.3106}}
\path(258.730,306.852)(260.000,315.000)(255.046,308.408)
\put(728.571,117.143){\arc{1017.263}{3.1570}{3.5411}}
\path(541.155,546.835)(540.000,555.000)(537.177,547.252)
\put(-398.333,653.333){\arc{1886.943}{0.1044}{0.3106}}
\path(497.687,467.085)(500.000,475.000)(494.235,469.104)
\put(765.000,320.000){\arc{614.003}{3.0273}{3.6708}}
\path(205,420)(275,420)
\path(267.000,418.000)(275.000,420.000)(267.000,422.000)
\path(485,260)(555,260)
\path(547.000,258.000)(555.000,260.000)(547.000,262.000)
\path(205,395)(235,365)
\path(227.929,369.243)(235.000,365.000)(230.757,372.071)
\path(525,315)(555,285)
\path(547.929,289.243)(555.000,285.000)(550.757,292.071)
\path(140,205)(180,395)
\path(180.309,386.760)(180.000,395.000)(176.395,387.584)
\path(300,445)(340,635)
\path(340.309,626.760)(340.000,635.000)(336.395,627.584)
\path(580,285)(620,475)
\path(620.309,466.760)(620.000,475.000)(616.395,467.584)
\path(420,45)(460,235)
\path(460.309,226.760)(460.000,235.000)(456.395,227.584)
\path(380,405)(420,555)
\path(419.871,546.755)(420.000,555.000)(416.006,547.785)
\put(75,185){\makebox(0,0)[lb]{\smash{{{\SetFigFont{12}{14.4}{bf}2}}}}}
\put(195,185){\makebox(0,0)[lb]{\smash{{{\SetFigFont{12}{14.4}{bf}2}}}}}
\put(275,105){\makebox(0,0)[lb]{\smash{{{\SetFigFont{12}{14.4}{bf}2}}}}}
\put(235,425){\makebox(0,0)[lb]{\smash{{{\SetFigFont{12}{14.4}{bf}2}}}}}
\put(315,325){\makebox(0,0)[lb]{\smash{{{\SetFigFont{12}{14.4}{bf}2}}}}}
\put(420,310){\makebox(0,0)[lb]{\smash{{{\SetFigFont{12}{14.4}{bf}2}}}}}
\put(515,265){\makebox(0,0)[lb]{\smash{{{\SetFigFont{12}{14.4}{bf}2}}}}}
\put(475,585){\makebox(0,0)[lb]{\smash{{{\SetFigFont{12}{14.4}{bf}2}}}}}
\put(555,505){\makebox(0,0)[lb]{\smash{{{\SetFigFont{12}{14.4}{bf}2}}}}}
\put(675,505){\makebox(0,0)[lb]{\smash{{{\SetFigFont{12}{14.4}{bf}2}}}}}
\put(180,145){\makebox(0,0)[lb]{\smash{{{\SetFigFont{12}{14.4}{bf}1}}}}}
\put(300,145){\makebox(0,0)[lb]{\smash{{{\SetFigFont{12}{14.4}{bf}1}}}}}
\put(380,65){\makebox(0,0)[lb]{\smash{{{\SetFigFont{12}{14.4}{bf}1}}}}}
\put(380,625){\makebox(0,0)[lb]{\smash{{{\SetFigFont{12}{14.4}{bf}1}}}}}
\put(460,545){\makebox(0,0)[lb]{\smash{{{\SetFigFont{12}{14.4}{bf}1}}}}}
\put(580,545){\makebox(0,0)[lb]{\smash{{{\SetFigFont{12}{14.4}{bf}1}}}}}
\put(220,385){\makebox(0,0)[lb]{\smash{{{\SetFigFont{12}{14.4}{bf}1}}}}}
\put(540,305){\makebox(0,0)[lb]{\smash{{{\SetFigFont{12}{14.4}{bf}1}}}}}
\put(340,420){\makebox(0,0)[lb]{\smash{{{\SetFigFont{12}{14.4}{bf}1}}}}}
\put(425,350){\makebox(0,0)[lb]{\smash{{{\SetFigFont{12}{14.4}{bf}1}}}}}
\put(145,300){\makebox(0,0)[lb]{\smash{{{\SetFigFont{12}{14.4}{bf}3}}}}}
\put(220,225){\makebox(0,0)[lb]{\smash{{{\SetFigFont{12}{14.4}{bf}3}}}}}
\put(285,270){\makebox(0,0)[lb]{\smash{{{\SetFigFont{12}{14.4}{bf}3}}}}}
\put(310,225){\makebox(0,0)[lb]{\smash{{{\SetFigFont{12}{14.4}{bf}3}}}}}
\put(430,145){\makebox(0,0)[lb]{\smash{{{\SetFigFont{12}{14.4}{bf}3}}}}}
\put(590,380){\makebox(0,0)[lb]{\smash{{{\SetFigFont{12}{14.4}{bf}3}}}}}
\put(525,445){\makebox(0,0)[lb]{\smash{{{\SetFigFont{12}{14.4}{bf}3}}}}}
\put(455,390){\makebox(0,0)[lb]{\smash{{{\SetFigFont{12}{14.4}{bf}3}}}}}
\put(305,535){\makebox(0,0)[lb]{\smash{{{\SetFigFont{12}{14.4}{bf}3}}}}}
\put(385,475){\makebox(0,0)[lb]{\smash{{{\SetFigFont{12}{14.4}{bf}3}}}}}
\put(5,185){\makebox(0,0)[lb]{\smash{{{\SetFigFont{12}{14.4}{bf}2}}}}}
\put(25,185){\makebox(0,0)[lb]{\smash{{{\SetFigFont{12}{14.4}{bf}2}}}}}
\put(5,165){\makebox(0,0)[lb]{\smash{{{\SetFigFont{12}{14.4}{bf}1}}}}}
\put(25,165){\makebox(0,0)[lb]{\smash{{{\SetFigFont{12}{14.4}{bf}1}}}}}
\put(125,185){\makebox(0,0)[lb]{\smash{{{\SetFigFont{12}{14.4}{bf}2}}}}}
\put(145,185){\makebox(0,0)[lb]{\smash{{{\SetFigFont{12}{14.4}{bf}3}}}}}
\put(125,165){\makebox(0,0)[lb]{\smash{{{\SetFigFont{12}{14.4}{bf}1}}}}}
\put(145,165){\makebox(0,0)[lb]{\smash{{{\SetFigFont{12}{14.4}{bf}1}}}}}
\put(205,105){\makebox(0,0)[lb]{\smash{{{\SetFigFont{12}{14.4}{bf}2}}}}}
\put(225,105){\makebox(0,0)[lb]{\smash{{{\SetFigFont{12}{14.4}{bf}3}}}}}
\put(205,85){\makebox(0,0)[lb]{\smash{{{\SetFigFont{12}{14.4}{bf}1}}}}}
\put(225,85){\makebox(0,0)[lb]{\smash{{{\SetFigFont{12}{14.4}{bf}2}}}}}
\put(325,105){\makebox(0,0)[lb]{\smash{{{\SetFigFont{12}{14.4}{bf}3}}}}}
\put(345,105){\makebox(0,0)[lb]{\smash{{{\SetFigFont{12}{14.4}{bf}3}}}}}
\put(325,85){\makebox(0,0)[lb]{\smash{{{\SetFigFont{12}{14.4}{bf}1}}}}}
\put(345,85){\makebox(0,0)[lb]{\smash{{{\SetFigFont{12}{14.4}{bf}2}}}}}
\put(405,25){\makebox(0,0)[lb]{\smash{{{\SetFigFont{12}{14.4}{bf}3}}}}}
\put(425,25){\makebox(0,0)[lb]{\smash{{{\SetFigFont{12}{14.4}{bf}3}}}}}
\put(405,5){\makebox(0,0)[lb]{\smash{{{\SetFigFont{12}{14.4}{bf}2}}}}}
\put(425,5){\makebox(0,0)[lb]{\smash{{{\SetFigFont{12}{14.4}{bf}2}}}}}
\put(245,185){\makebox(0,0)[lb]{\smash{{{\SetFigFont{12}{14.4}{bf}3}}}}}
\put(265,185){\makebox(0,0)[lb]{\smash{{{\SetFigFont{12}{14.4}{bf}3}}}}}
\put(245,165){\makebox(0,0)[lb]{\smash{{{\SetFigFont{12}{14.4}{bf}1}}}}}
\put(265,165){\makebox(0,0)[lb]{\smash{{{\SetFigFont{12}{14.4}{bf}1}}}}}
\put(165,425){\makebox(0,0)[lb]{\smash{{{\SetFigFont{12}{14.4}{bf}2}}}}}
\put(185,425){\makebox(0,0)[lb]{\smash{{{\SetFigFont{12}{14.4}{bf}4}}}}}
\put(165,405){\makebox(0,0)[lb]{\smash{{{\SetFigFont{12}{14.4}{bf}1}}}}}
\put(185,405){\makebox(0,0)[lb]{\smash{{{\SetFigFont{12}{14.4}{bf}1}}}}}
\put(285,425){\makebox(0,0)[lb]{\smash{{{\SetFigFont{12}{14.4}{bf}3}}}}}
\put(305,425){\makebox(0,0)[lb]{\smash{{{\SetFigFont{12}{14.4}{bf}4}}}}}
\put(285,405){\makebox(0,0)[lb]{\smash{{{\SetFigFont{12}{14.4}{bf}1}}}}}
\put(305,405){\makebox(0,0)[lb]{\smash{{{\SetFigFont{12}{14.4}{bf}1}}}}}
\put(245,345){\makebox(0,0)[lb]{\smash{{{\SetFigFont{12}{14.4}{bf}2}}}}}
\put(265,345){\makebox(0,0)[lb]{\smash{{{\SetFigFont{12}{14.4}{bf}4}}}}}
\put(245,325){\makebox(0,0)[lb]{\smash{{{\SetFigFont{12}{14.4}{bf}1}}}}}
\put(265,325){\makebox(0,0)[lb]{\smash{{{\SetFigFont{12}{14.4}{bf}2}}}}}
\put(365,385){\makebox(0,0)[lb]{\smash{{{\SetFigFont{12}{14.4}{bf}3}}}}}
\put(385,385){\makebox(0,0)[lb]{\smash{{{\SetFigFont{12}{14.4}{bf}4}}}}}
\put(365,365){\makebox(0,0)[lb]{\smash{{{\SetFigFont{12}{14.4}{bf}1}}}}}
\put(385,365){\makebox(0,0)[lb]{\smash{{{\SetFigFont{12}{14.4}{bf}2}}}}}
\put(365,305){\makebox(0,0)[lb]{\smash{{{\SetFigFont{12}{14.4}{bf}2}}}}}
\put(385,305){\makebox(0,0)[lb]{\smash{{{\SetFigFont{12}{14.4}{bf}4}}}}}
\put(365,285){\makebox(0,0)[lb]{\smash{{{\SetFigFont{12}{14.4}{bf}1}}}}}
\put(385,285){\makebox(0,0)[lb]{\smash{{{\SetFigFont{12}{14.4}{bf}3}}}}}
\put(445,265){\makebox(0,0)[lb]{\smash{{{\SetFigFont{12}{14.4}{bf}3}}}}}
\put(465,265){\makebox(0,0)[lb]{\smash{{{\SetFigFont{12}{14.4}{bf}4}}}}}
\put(445,245){\makebox(0,0)[lb]{\smash{{{\SetFigFont{12}{14.4}{bf}2}}}}}
\put(465,245){\makebox(0,0)[lb]{\smash{{{\SetFigFont{12}{14.4}{bf}2}}}}}
\put(485,345){\makebox(0,0)[lb]{\smash{{{\SetFigFont{12}{14.4}{bf}3}}}}}
\put(505,345){\makebox(0,0)[lb]{\smash{{{\SetFigFont{12}{14.4}{bf}4}}}}}
\put(485,325){\makebox(0,0)[lb]{\smash{{{\SetFigFont{12}{14.4}{bf}1}}}}}
\put(505,325){\makebox(0,0)[lb]{\smash{{{\SetFigFont{12}{14.4}{bf}3}}}}}
\put(565,265){\makebox(0,0)[lb]{\smash{{{\SetFigFont{12}{14.4}{bf}3}}}}}
\put(585,265){\makebox(0,0)[lb]{\smash{{{\SetFigFont{12}{14.4}{bf}4}}}}}
\put(565,245){\makebox(0,0)[lb]{\smash{{{\SetFigFont{12}{14.4}{bf}2}}}}}
\put(585,245){\makebox(0,0)[lb]{\smash{{{\SetFigFont{12}{14.4}{bf}3}}}}}
\put(325,665){\makebox(0,0)[lb]{\smash{{{\SetFigFont{12}{14.4}{bf}4}}}}}
\put(345,665){\makebox(0,0)[lb]{\smash{{{\SetFigFont{12}{14.4}{bf}4}}}}}
\put(325,645){\makebox(0,0)[lb]{\smash{{{\SetFigFont{12}{14.4}{bf}1}}}}}
\put(345,645){\makebox(0,0)[lb]{\smash{{{\SetFigFont{12}{14.4}{bf}1}}}}}
\put(405,585){\makebox(0,0)[lb]{\smash{{{\SetFigFont{12}{14.4}{bf}4}}}}}
\put(425,585){\makebox(0,0)[lb]{\smash{{{\SetFigFont{12}{14.4}{bf}4}}}}}
\put(405,565){\makebox(0,0)[lb]{\smash{{{\SetFigFont{12}{14.4}{bf}1}}}}}
\put(425,565){\makebox(0,0)[lb]{\smash{{{\SetFigFont{12}{14.4}{bf}2}}}}}
\put(525,585){\makebox(0,0)[lb]{\smash{{{\SetFigFont{12}{14.4}{bf}4}}}}}
\put(545,585){\makebox(0,0)[lb]{\smash{{{\SetFigFont{12}{14.4}{bf}4}}}}}
\put(525,565){\makebox(0,0)[lb]{\smash{{{\SetFigFont{12}{14.4}{bf}1}}}}}
\put(545,565){\makebox(0,0)[lb]{\smash{{{\SetFigFont{12}{14.4}{bf}3}}}}}
\put(485,505){\makebox(0,0)[lb]{\smash{{{\SetFigFont{12}{14.4}{bf}4}}}}}
\put(505,505){\makebox(0,0)[lb]{\smash{{{\SetFigFont{12}{14.4}{bf}4}}}}}
\put(505,485){\makebox(0,0)[lb]{\smash{{{\SetFigFont{12}{14.4}{bf}2}}}}}
\put(485,485){\makebox(0,0)[lb]{\smash{{{\SetFigFont{12}{14.4}{bf}2}}}}}
\put(605,505){\makebox(0,0)[lb]{\smash{{{\SetFigFont{12}{14.4}{bf}4}}}}}
\put(625,505){\makebox(0,0)[lb]{\smash{{{\SetFigFont{12}{14.4}{bf}4}}}}}
\put(605,485){\makebox(0,0)[lb]{\smash{{{\SetFigFont{12}{14.4}{bf}2}}}}}
\put(625,485){\makebox(0,0)[lb]{\smash{{{\SetFigFont{12}{14.4}{bf}3}}}}}
\put(725,505){\makebox(0,0)[lb]{\smash{{{\SetFigFont{12}{14.4}{bf}4}}}}}
\put(745,505){\makebox(0,0)[lb]{\smash{{{\SetFigFont{12}{14.4}{bf}4}}}}}
\put(725,485){\makebox(0,0)[lb]{\smash{{{\SetFigFont{12}{14.4}{bf}3}}}}}
\put(745,485){\makebox(0,0)[lb]{\smash{{{\SetFigFont{12}{14.4}{bf}3}}}}}
\end{picture}
\end{center}
\caption{Crystal graph of the $U_q(\gl_4)$-module $V_{(2,\,2)}$\label{V22}}
\end{figure}
}\end{example}

\section{$\F_q[GL_n/B]$} \label{qG/B}

The subalgebra $\F_q[GL_n/B]$ of $\F_q[\mat_n]$ generated by quantum
column-shaped
tableaux can also be defined by generators and relations \cite{TT,LR}. The
generators are denoted by columns
$$\bo{\matrix{ i_k \cr
	       \vdots \cr
	       i_2 \cr
	       i_1 }}
\quad\quad 1\le k \le n,\quad 1\le i_1, \ldots ,i_k \le n ,
$$
and their products are written by juxtaposition. The relations are
\begin{description}
\item[$(R_1)$] \quad if $i_r=i_s$ for some indices $ r,\,s$, then
$$
\bo{\matrix{ i_k \cr
\vdots \cr
i_2 \cr
i_1 }}
= 0 \,,
$$
\item[$(R_2)$] \quad for $w\in S_k$  and $i_1<i_2<\cdots <i_k$, we have
$$
\bo{\matrix{ i_{w_k} \cr
\vdots \cr
i_{w_2} \cr
i_{w_1} }}
= (-q)^{-\ell (w)}\
\bo{\matrix{ i_k \cr
\vdots \cr
i_2 \cr
i_1 }}\,,
$$
\item[$(R_3)$] \quad for $k\le l$ and $j_1<j_2<\cdots <j_l$, we have
\begin{center}
\setlength{\unitlength}{0.01in}
\begingroup\makeatletter\ifx\SetFigFont\undefined
\def\x#1#2#3#4#5#6#7\relax{\def\x{#1#2#3#4#5#6}}%
\expandafter\x\fmtname xxxxxx\relax \def\y{splain}%
\ifx\x\y   
\gdef\SetFigFont#1#2#3{%
  \ifnum #1<17\tiny\else \ifnum #1<20\small\else
  \ifnum #1<24\normalsize\else \ifnum #1<29\large\else
  \ifnum #1<34\Large\else \ifnum #1<41\LARGE\else
     \huge\fi\fi\fi\fi\fi\fi
  \csname #3\endcsname}%
\else
\gdef\SetFigFont#1#2#3{\begingroup
  \count@#1\relax \ifnum 25<\count@\count@25\fi
  \def\x{\endgroup\@setsize\SetFigFont{#2pt}}%
  \expandafter\x
    \csname \romannumeral\the\count@ pt\expandafter\endcsname
    \csname @\romannumeral\the\count@ pt\endcsname
  \csname #3\endcsname}%
\fi
\fi\endgroup
\begin{picture}(230,135)(0,-10)
\path(60,0)(60,120)(30,120)
	(30,0)(60,0)
\path(30,0)(30,60)(0,60)
	(0,0)(30,0)
\path(200,0)(200,120)(165,120)
	(165,0)(200,0)
\path(230,0)(230,60)(200,60)
	(200,0)(230,0)
\put(10,5){\makebox(0,0)[lb]{\smash{{{\SetFigFont{12}{14.4}{rm}$i_1$}}}}}
\put(10,22){\makebox(0,0)[lb]{\smash{{{\SetFigFont{12}{14.4}{rm}$\vdots$}}}}}
\put(10,45){\makebox(0,0)[lb]{\smash{{{\SetFigFont{12}{14.4}{rm}$i_k$}}}}}
\put(40,5){\makebox(0,0)[lb]{\smash{{{\SetFigFont{12}{14.4}{rm}$j_1$}}}}}
\put(40,22){\makebox(0,0)[lb]{\smash{{{\SetFigFont{12}{14.4}{rm}$\vdots$}}}}}
\put(40,45){\makebox(0,0)[lb]{\smash{{{\SetFigFont{12}{14.4}{rm}$j_k$}}}}}
\put(32,65){\makebox(0,0)[lb]{\smash{{{\SetFigFont{12}{14.4}{rm}$j_{k+1}$}}}}}
\put(40,82){\makebox(0,0)[lb]{\smash{{{\SetFigFont{12}{14.4}{rm}$\vdots$}}}}}
\put(40,105){\makebox(0,0)[lb]{\smash{{{\SetFigFont{12}{14.4}{rm}$j_l$}}}}}
\put(175,5){\makebox(0,0)[lb]{\smash{{{\SetFigFont{12}{14.4}{rm}$i_1$}}}}}
\put(177,22){\makebox(0,0)[lb]{\smash{{{\SetFigFont{12}{14.4}{rm}$\vdots$}}}}}
\put(175,45){\makebox(0,0)[lb]{\smash{{{\SetFigFont{12}{14.4}{rm}$i_k$}}}}}
\put(167,65){\makebox(0,0)[lb]{\smash{{{\SetFigFont{12}{14.4}{rm}$j_{w_{k+1}}$}}}}}
\put(177,82){\makebox(0,0)[lb]{\smash{{{\SetFigFont{12}{14.4}{rm}$\vdots$}}}}}
\put(175,105){\makebox(0,0)[lb]{\smash{{{\SetFigFont{12}{14.4}{rm}$j_{w_l}$}}}}}
\put(208,5){\makebox(0,0)[lb]{\smash{{{\SetFigFont{12}{14.4}{rm}$j_{w_1}$}}}}}
\put(215,22){\makebox(0,0)[lb]{\smash{{{\SetFigFont{12}{14.4}{rm}$\vdots$}}}}}
\put(208,45){\makebox(0,0)[lb]{\smash{{{\SetFigFont{12}{14.4}{rm}$j_{w_k}$}}}}}
\put(65,45){\makebox(0,0)[lb]{\smash{{{\SetFigFont{12}{14.4}{rm}$\displaystyle
= \sum_w (-q)^{\ell(w)}$}}}}}
\end{picture}
\end{center}
where the sum runs through all $w\in S_l$ such that
$w_1<\cdots <w_k$ and $w_{k+1}<\cdots <w_l$,
\item[$(R_4)$] \quad for $k+s \le l-1$ and  $j_1<j_2<\cdots <j_l$, we have
\begin{center}
\setlength{\unitlength}{0.0118in}
\begingroup\makeatletter\ifx\SetFigFont\undefined
\def\x#1#2#3#4#5#6#7\relax{\def\x{#1#2#3#4#5#6}}%
\expandafter\x\fmtname xxxxxx\relax \def\y{splain}%
\ifx\x\y   
\gdef\SetFigFont#1#2#3{%
  \ifnum #1<17\tiny\else \ifnum #1<20\small\else
  \ifnum #1<24\normalsize\else \ifnum #1<29\large\else
  \ifnum #1<34\Large\else \ifnum #1<41\LARGE\else
     \huge\fi\fi\fi\fi\fi\fi
  \csname #3\endcsname}%
\else
\gdef\SetFigFont#1#2#3{\begingroup
  \count@#1\relax \ifnum 25<\count@\count@25\fi
  \def\x{\endgroup\@setsize\SetFigFont{#2pt}}%
  \expandafter\x
    \csname \romannumeral\the\count@ pt\expandafter\endcsname
    \csname @\romannumeral\the\count@ pt\endcsname
  \csname #3\endcsname}%
\fi
\fi\endgroup
\begin{picture}(173,175)(0,-10)
\path(100,0)(100,160)(60,160)
	(60,0)(100,0)
\path(140,0)(140,140)(100,140)
	(100,0)(140,0)
\put(70,5){\makebox(0,0)[lb]{\smash{{{\SetFigFont{12}{14.4}{rm}$i_1$}}}}}
\put(70,25){\makebox(0,0)[lb]{\smash{{{\SetFigFont{12}{14.4}{rm}$\vdots$}}}}}
\put(70,45){\makebox(0,0)[lb]{\smash{{{\SetFigFont{12}{14.4}{rm}$i_{k-1}$}}}}}
\put(70,65){\makebox(0,0)[lb]{\smash{{{\SetFigFont{12}{14.4}{rm}$j_{w_{k+1}}$}}}}}
\put(70,85){\makebox(0,0)[lb]{\smash{{{\SetFigFont{12}{14.4}{rm}$j_{w_{k+2}}$}}}}}
\put(70,115){\makebox(0,0)[lb]{\smash{{{\SetFigFont{12}{14.4}{rm}$\vdots$}}}}}
\put(70,145){\makebox(0,0)[lb]{\smash{{{\SetFigFont{12}{14.4}{rm}$j_{w_l}$}}}}}
\put(110,5){\makebox(0,0)[lb]{\smash{{{\SetFigFont{12}{14.4}{rm}$j_{w_1}$}}}}}
\put(110,25){\makebox(0,0)[lb]{\smash{{{\SetFigFont{12}{14.4}{rm}$\vdots$}}}}}
\put(110,45){\makebox(0,0)[lb]{\smash{{{\SetFigFont{12}{14.4}{rm}$j_{w_{k-1}}$}}}}}
\put(110,65){\makebox(0,0)[lb]{\smash{{{\SetFigFont{12}{14.4}{rm}$j_{w_k}$}}}}}
\put(110,85){\makebox(0,0)[lb]{\smash{{{\SetFigFont{12}{14.4}{rm}$r_1$}}}}}
\put(110,105){\makebox(0,0)[lb]{\smash{{{\SetFigFont{12}{14.4}{rm}$\vdots$}}}}}
\put(110,125){\makebox(0,0)[lb]{\smash{{{\SetFigFont{12}{14.4}{rm}$r_s$}}}}}
\put(-5,55){\makebox(0,0)[lb]{\smash{{{\SetFigFont{12}{14.4}{rm}$\displaystyle
\sum_w (-q)^{\ell (w)}$}}}}}
\put(150,55){\makebox(0,0)[lb]{\smash{{{\SetFigFont{12}{14.4}{rm}$=0$}}}}}
\end{picture}
\end{center}
where the sum runs through all permutations $w\in S_l$ such that
$w_1<\cdots <w_k$ and $w_{k+1}<\cdots <w_l$.
\end{description}
Relations $(R_1)$, $(R_2)$ express the $q$-alternating property of quantum
minors, while relations $(R_3)$, $(R_4)$ are the natural $q$-analogues
of Sylvester's and Garnir's identities respectively (\cf \cite{Le}).
Using  $(R_1)$, $(R_2)$, $(R_3)$ one can express every element of
$\F_q[GL_n/B]$ as a linear combination of products of columns, the size
of which weakly decreases from left to right. After that, following the
same strategy as in the classical straightening algorithm, one can by
means of $(R_4)$ express all such products as linear combinations of
quantum tableaux. Note that it follows from $(R_3)$, $(R_4)$ that all the
quantum tableaux $(\tau)$ appearing in the straightening of a given
tabloid $(\delta)$ have the same shape $\lambda$ obtained by reordering
the sizes of the columns of $\delta$, which means that $(\delta)$ falls
in the irreducible component $V_\lambda$ of $\F_q[GL_n/B]$.
\begin{example}{\rm
We apply the $q$-straightening algorithm to
\begin{center}
\setlength{\unitlength}{0.01in}
\begingroup\makeatletter\ifx\SetFigFont\undefined
\def\x#1#2#3#4#5#6#7\relax{\def\x{#1#2#3#4#5#6}}%
\expandafter\x\fmtname xxxxxx\relax \def\y{splain}%
\ifx\x\y   
\gdef\SetFigFont#1#2#3{%
  \ifnum #1<17\tiny\else \ifnum #1<20\small\else
  \ifnum #1<24\normalsize\else \ifnum #1<29\large\else
  \ifnum #1<34\Large\else \ifnum #1<41\LARGE\else
     \huge\fi\fi\fi\fi\fi\fi
  \csname #3\endcsname}%
\else
\gdef\SetFigFont#1#2#3{\begingroup
  \count@#1\relax \ifnum 25<\count@\count@25\fi
  \def\x{\endgroup\@setsize\SetFigFont{#2pt}}%
  \expandafter\x
    \csname \romannumeral\the\count@ pt\expandafter\endcsname
    \csname @\romannumeral\the\count@ pt\endcsname
  \csname #3\endcsname}%
\fi
\fi\endgroup
\begin{picture}(70,75)(0,-10)
\path(50,0)(50,40)(30,40)
	(30,0)(50,0)
\path(70,0)(70,60)(50,60)
	(50,0)(70,0)
\put(-15,15){\makebox(0,0)[lb]{\smash{{{\SetFigFont{12}{14.4}{rm}$(\delta)
=$}}}}}
\put(35,5){\makebox(0,0)[lb]{\smash{{{\SetFigFont{12}{14.4}{rm}1}}}}}
\put(35,25){\makebox(0,0)[lb]{\smash{{{\SetFigFont{12}{14.4}{rm}5}}}}}
\put(55,5){\makebox(0,0)[lb]{\smash{{{\SetFigFont{12}{14.4}{rm}2}}}}}
\put(55,25){\makebox(0,0)[lb]{\smash{{{\SetFigFont{12}{14.4}{rm}3}}}}}
\put(55,45){\makebox(0,0)[lb]{\smash{{{\SetFigFont{12}{14.4}{rm}6}}}}}
\end{picture}
\end{center}
It follows from  $(R_1)$, $(R_2)$, $(R_3)$ that
\begin{center}
\setlength{\unitlength}{0.01in}
\begingroup\makeatletter\ifx\SetFigFont\undefined
\def\x#1#2#3#4#5#6#7\relax{\def\x{#1#2#3#4#5#6}}%
\expandafter\x\fmtname xxxxxx\relax \def\y{splain}%
\ifx\x\y   
\gdef\SetFigFont#1#2#3{%
  \ifnum #1<17\tiny\else \ifnum #1<20\small\else
  \ifnum #1<24\normalsize\else \ifnum #1<29\large\else
  \ifnum #1<34\Large\else \ifnum #1<41\LARGE\else
     \huge\fi\fi\fi\fi\fi\fi
  \csname #3\endcsname}%
\else
\gdef\SetFigFont#1#2#3{\begingroup
  \count@#1\relax \ifnum 25<\count@\count@25\fi
  \def\x{\endgroup\@setsize\SetFigFont{#2pt}}%
  \expandafter\x
    \csname \romannumeral\the\count@ pt\expandafter\endcsname
    \csname @\romannumeral\the\count@ pt\endcsname
  \csname #3\endcsname}%
\fi
\fi\endgroup
\begin{picture}(280,75)(0,-10)
\path(20,0)(20,40)(0,40)
	(0,0)(20,0)
\path(40,0)(40,60)(20,60)
	(20,0)(40,0)
\path(100,0)(100,60)(80,60)
	(80,0)(100,0)
\path(120,0)(120,40)(100,40)
	(100,0)(120,0)
\path(180,0)(180,60)(160,60)
	(160,0)(180,0)
\path(200,0)(200,40)(180,40)
	(180,0)(200,0)
\path(260,0)(260,60)(240,60)
	(240,0)(260,0)
\path(280,0)(280,40)(260,40)
	(260,0)(280,0)
\put(5,5){\makebox(0,0)[lb]{\smash{{{\SetFigFont{12}{14.4}{rm}1}}}}}
\put(5,25){\makebox(0,0)[lb]{\smash{{{\SetFigFont{12}{14.4}{rm}5}}}}}
\put(25,5){\makebox(0,0)[lb]{\smash{{{\SetFigFont{12}{14.4}{rm}2}}}}}
\put(25,25){\makebox(0,0)[lb]{\smash{{{\SetFigFont{12}{14.4}{rm}3}}}}}
\put(25,45){\makebox(0,0)[lb]{\smash{{{\SetFigFont{12}{14.4}{rm}6}}}}}
\put(85,5){\makebox(0,0)[lb]{\smash{{{\SetFigFont{12}{14.4}{rm}1}}}}}
\put(85,25){\makebox(0,0)[lb]{\smash{{{\SetFigFont{12}{14.4}{rm}5}}}}}
\put(85,45){\makebox(0,0)[lb]{\smash{{{\SetFigFont{12}{14.4}{rm}6}}}}}
\put(105,5){\makebox(0,0)[lb]{\smash{{{\SetFigFont{12}{14.4}{rm}2}}}}}
\put(105,25){\makebox(0,0)[lb]{\smash{{{\SetFigFont{12}{14.4}{rm}3}}}}}
\put(165,5){\makebox(0,0)[lb]{\smash{{{\SetFigFont{12}{14.4}{rm}1}}}}}
\put(165,25){\makebox(0,0)[lb]{\smash{{{\SetFigFont{12}{14.4}{rm}3}}}}}
\put(165,45){\makebox(0,0)[lb]{\smash{{{\SetFigFont{12}{14.4}{rm}5}}}}}
\put(185,5){\makebox(0,0)[lb]{\smash{{{\SetFigFont{12}{14.4}{rm}2}}}}}
\put(185,25){\makebox(0,0)[lb]{\smash{{{\SetFigFont{12}{14.4}{rm}6}}}}}
\put(245,5){\makebox(0,0)[lb]{\smash{{{\SetFigFont{12}{14.4}{rm}1}}}}}
\put(245,25){\makebox(0,0)[lb]{\smash{{{\SetFigFont{12}{14.4}{rm}2}}}}}
\put(245,45){\makebox(0,0)[lb]{\smash{{{\SetFigFont{12}{14.4}{rm}5}}}}}
\put(265,5){\makebox(0,0)[lb]{\smash{{{\SetFigFont{12}{14.4}{rm}3}}}}}
\put(265,25){\makebox(0,0)[lb]{\smash{{{\SetFigFont{12}{14.4}{rm}6}}}}}
\put(55,20){\makebox(0,0)[lb]{\smash{{{\SetFigFont{12}{14.4}{rm}$=$}}}}}
\put(135,20){\makebox(0,0)[lb]{\smash{{{\SetFigFont{12}{14.4}{rm}$+$}}}}}
\put(215,20){\makebox(0,0)[lb]{\smash{{{\SetFigFont{12}{14.4}{rm}$-q$}}}}}
\end{picture}
\end{center}
Using $(R_4)$ for straightening the first term of the right-hand side, we
obtain
\begin{center}
\setlength{\unitlength}{0.01in}
\begingroup\makeatletter\ifx\SetFigFont\undefined
\def\x#1#2#3#4#5#6#7\relax{\def\x{#1#2#3#4#5#6}}%
\expandafter\x\fmtname xxxxxx\relax \def\y{splain}%
\ifx\x\y   
\gdef\SetFigFont#1#2#3{%
  \ifnum #1<17\tiny\else \ifnum #1<20\small\else
  \ifnum #1<24\normalsize\else \ifnum #1<29\large\else
  \ifnum #1<34\Large\else \ifnum #1<41\LARGE\else
     \huge\fi\fi\fi\fi\fi\fi
  \csname #3\endcsname}%
\else
\gdef\SetFigFont#1#2#3{\begingroup
  \count@#1\relax \ifnum 25<\count@\count@25\fi
  \def\x{\endgroup\@setsize\SetFigFont{#2pt}}%
  \expandafter\x
    \csname \romannumeral\the\count@ pt\expandafter\endcsname
    \csname @\romannumeral\the\count@ pt\endcsname
  \csname #3\endcsname}%
\fi
\fi\endgroup
\begin{picture}(580,75)(0,-10)
\path(20,0)(20,40)(0,40)
	(0,0)(20,0)
\path(40,0)(40,60)(20,60)
	(20,0)(40,0)
\put(5,5){\makebox(0,0)[lb]{\smash{{{\SetFigFont{12}{14.4}{rm}1}}}}}
\put(5,25){\makebox(0,0)[lb]{\smash{{{\SetFigFont{12}{14.4}{rm}5}}}}}
\put(25,5){\makebox(0,0)[lb]{\smash{{{\SetFigFont{12}{14.4}{rm}2}}}}}
\put(25,25){\makebox(0,0)[lb]{\smash{{{\SetFigFont{12}{14.4}{rm}3}}}}}
\put(25,45){\makebox(0,0)[lb]{\smash{{{\SetFigFont{12}{14.4}{rm}6}}}}}
\put(50,20){\makebox(0,0)[lb]{\smash{{{\SetFigFont{12}{14.4}{rm}$=$}}}}}
\put(62,20){\makebox(0,0)[lb]{\smash{{{\SetFigFont{12}{14.4}{rm}$(1-q^2)$}}}}}
\path(140,0)(140,60)(120,60)
	(120,0)(140,0)
\path(160,0)(160,40)(140,40)
	(140,0)(160,0)
\put(125,5){\makebox(0,0)[lb]{\smash{{{\SetFigFont{12}{14.4}{rm}1}}}}}
\put(125,25){\makebox(0,0)[lb]{\smash{{{\SetFigFont{12}{14.4}{rm}3}}}}}
\put(125,45){\makebox(0,0)[lb]{\smash{{{\SetFigFont{12}{14.4}{rm}5}}}}}
\put(145,5){\makebox(0,0)[lb]{\smash{{{\SetFigFont{12}{14.4}{rm}2}}}}}
\put(145,25){\makebox(0,0)[lb]{\smash{{{\SetFigFont{12}{14.4}{rm}6}}}}}
\path(280,0)(280,40)(260,40)
	(260,0)(280,0)
\path(260,0)(260,60)(240,60)
	(240,0)(260,0)
\put(245,5){\makebox(0,0)[lb]{\smash{{{\SetFigFont{12}{14.4}{rm}1}}}}}
\put(245,25){\makebox(0,0)[lb]{\smash{{{\SetFigFont{12}{14.4}{rm}2}}}}}
\put(245,45){\makebox(0,0)[lb]{\smash{{{\SetFigFont{12}{14.4}{rm}5}}}}}
\put(265,5){\makebox(0,0)[lb]{\smash{{{\SetFigFont{12}{14.4}{rm}3}}}}}
\put(265,25){\makebox(0,0)[lb]{\smash{{{\SetFigFont{12}{14.4}{rm}6}}}}}
\path(360,0)(360,60)(340,60)
	(340,0)(360,0)
\path(380,0)(380,40)(360,40)
	(360,0)(380,0)
\put(345,5){\makebox(0,0)[lb]{\smash{{{\SetFigFont{12}{14.4}{rm}1}}}}}
\put(345,25){\makebox(0,0)[lb]{\smash{{{\SetFigFont{12}{14.4}{rm}3}}}}}
\put(345,45){\makebox(0,0)[lb]{\smash{{{\SetFigFont{12}{14.4}{rm}6}}}}}
\put(365,5){\makebox(0,0)[lb]{\smash{{{\SetFigFont{12}{14.4}{rm}2}}}}}
\put(365,25){\makebox(0,0)[lb]{\smash{{{\SetFigFont{12}{14.4}{rm}5}}}}}
\path(560,0)(560,60)(540,60)
	(540,0)(560,0)
\path(580,0)(580,40)(560,40)
	(560,0)(580,0)
\put(545,5){\makebox(0,0)[lb]{\smash{{{\SetFigFont{12}{14.4}{rm}1}}}}}
\put(545,25){\makebox(0,0)[lb]{\smash{{{\SetFigFont{12}{14.4}{rm}2}}}}}
\put(545,45){\makebox(0,0)[lb]{\smash{{{\SetFigFont{12}{14.4}{rm}3}}}}}
\put(565,5){\makebox(0,0)[lb]{\smash{{{\SetFigFont{12}{14.4}{rm}5}}}}}
\put(565,25){\makebox(0,0)[lb]{\smash{{{\SetFigFont{12}{14.4}{rm}6}}}}}
\path(460,0)(460,60)(440,60)
	(440,0)(460,0)
\path(480,0)(480,40)(460,40)
	(460,0)(480,0)
\put(445,5){\makebox(0,0)[lb]{\smash{{{\SetFigFont{12}{14.4}{rm}1}}}}}
\put(445,25){\makebox(0,0)[lb]{\smash{{{\SetFigFont{12}{14.4}{rm}2}}}}}
\put(445,45){\makebox(0,0)[lb]{\smash{{{\SetFigFont{12}{14.4}{rm}6}}}}}
\put(465,5){\makebox(0,0)[lb]{\smash{{{\SetFigFont{12}{14.4}{rm}3}}}}}
\put(465,25){\makebox(0,0)[lb]{\smash{{{\SetFigFont{12}{14.4}{rm}5}}}}}
\put(168,20){\makebox(0,0)[lb]{\smash{{{\SetFigFont{12}{14.4}{rm}$+(q^3-q)$}}}}}
\put(300,20){\makebox(0,0)[lb]{\smash{{{\SetFigFont{12}{14.4}{rm}$+\ q$}}}}}
\put(500,20){\makebox(0,0)[lb]{\smash{{{\SetFigFont{12}{14.4}{rm}$-\ q^4$}}}}}
\put(395,20){\makebox(0,0)[lb]{\smash{{{\SetFigFont{12}{14.4}{rm}$-\ q^2$}}}}}
\end{picture}
\end{center}
Letting $q$ tend to $0$ in this expansion, we get an illustration of the
next Theorem, which can be regarded as a version of Theorem~\ref{TH1} for
$\F_q[GL_n/B]$.
}\end{example}
\begin{theorem}\label{TH2}
Let ${\cal B}$ denote the linear basis of  $\F_q[GL_n/B]$ consisting of
quantum tableaux, and let ${\cal L}$ be the $K[q]$-lattice generated by
${\cal B}$. Let $\lambda$ be a partition, and $(\delta)$ a quantum tabloid
in $V_\lambda$. Denote by $w$ the word obtained by reading the columns
of $\delta$ from top to bottom and left to right. Then,
$$
(\delta) \equiv \left\{
\matrix{ (P(w)) & \mod \ q{\cal L} & {\it if\ P(w) \ has \ shape \ \lambda} \cr
0 &  \mod \ q{\cal L} &{\it otherwise}
}\right.\,.
$$
\end{theorem}
The proof of Theorem~\ref{TH2} is derived from the following result
which describes several interesting bases of $V_\lambda$, and states that
each of them gives rise to the same crystal basis at $q=0$. We first
introduce some notations. Let $\lambda$ be a partition,
$\lambda' = (\lambda_1',\ldots ,\lambda_r')$ its conjugate, and
$\sigma (\lambda') = (\lambda_{\sigma_1}',\ldots ,\lambda_{\sigma_r}')$
the composition obtained by permuting the parts of $\lambda'$ using
the permutation $\sigma \in S_r$. A tabloid $\delta$
is said to have {\it shape} $\sigma (\lambda')$ if it is a sequence of
column-shaped Young tableaux $(c_1,\ldots ,c_r)$, the size of the column $c_i$
being
equal to $\lambda_{\sigma_i}'$, $i=1,\ldots ,r$. The {\it column reading}
of $\delta$ is the word $u_\delta$ obtained by reading the columns
of $\delta$ from top to bottom and left to right.
\begin{theorem}\label{TH3}
For any $\sigma \in S_r$, the set of quantum tabloids
$$
B_\sigma = \{(\delta)\ |\ \delta \ {\it has \ shape}\ \sigma(\lambda')\
{\it and } \ P(u_\delta) \ {\it has \ shape} \ \lambda\}
$$
is a linear basis of $V_\lambda$. The ${\cal A}$-lattice $L$ spanned by
$B_\sigma$ is
independent of $\sigma$ as well as the projection $B$ of $B_\sigma$ in
$L/qL$, and $(L,B)$ is a crystal basis of $V_\lambda$. For $(\delta) \in
B_\sigma$ and
$(\gamma) \in B_{\zeta}$, we have
$$
(\delta) \equiv (\gamma)\ \mod qL \quad \Longleftrightarrow
\quad P(u_\delta) = P(u_{\gamma})\,.
$$
Finally, if $(\delta)$ is such that $P(u_\delta)$ has not shape $\lambda$,
then $(\delta) \equiv 0 \ \mod qL$.
\end{theorem}

\Proof We fix $\sigma \in S_r$ and we write $\mu = (\mu_1, \ldots , \mu_r) :=
\sigma (\lambda')$. Recall from Example~\ref{PEXT} that the $q$-analogue
$V_{(1^k)}$ of the $k$ th exterior power of the basic representation has
a natural basis labelled by column-shaped Young tableaux. We introduce the
$K(q)$-linear map $\varphi_\mu$:
\begin{eqnarray*}
V_{(1^{\mu_1})} \otimes V_{(1^{\mu_2})}  \otimes \cdots \otimes V_{(1^{\mu_r})}
& \longrightarrow & V_\lambda \\
\bo{\matrix{b \cr \vdots \cr a}} \otimes
\bo{\matrix{d \cr \vdots \cr \vdots \cr c}}\otimes \cdots \otimes
\bo{\matrix{f \cr \vdots \cr e}}
& \longrightarrow &
\hbox{
\setlength{\unitlength}{0.02in}\begin{picture}(50,40)(0,17)
\path(0,0)(10,0)(10,30)(0,30)(0,0)
\path(10,0)(20,0)(20,40)(10,40)(10,0)
\path(40,0)(50,0)(50,30)(40,30)(40,0)
\put(3,3){$a$}
\put(3,13){$\vdots$}
\put(3,23){$b$}
\put(13,3){$c$}
\put(13,13){$\vdots$}
\put(13,23){$\vdots$}
\put(13,33){$d$}
\put(25,15){$\cdots$}
\put(43,3){$e$}
\put(43,13){$\vdots$}
\put(43,23){$f$}
\end{picture}
}
\end{eqnarray*}
Comparing (\ref{ATP1}) (\ref{ATP2}) (\ref{ATP3}) to (\ref{L1}) (\ref{L2})
(\ref{L3})
we see immediately that $\varphi_\mu$ commutes with the action of $U_q(\gl_n)$
and
hence is in fact a morphism of $U_q(\gl_n)$-modules. The image $V_\lambda$ is
irreducible and by (\ref{LR}) it appears with multiplicity one in
$V_{(1^{\mu_1})} \otimes \cdots \otimes V_{(1^{\mu_r})}$.
Therefore, the restriction of $\varphi_\mu$ to the irreducible component
of $V_{(1^{\mu_1})} \otimes \cdots \otimes V_{(1^{\mu_r})}$ with
highest weight vector the tensor product
$y^\mu := (y_{(1^{\mu_1})}) \otimes \cdots \otimes (y_{(1^{\mu_r})})$ of all
highest weight vectors  in $V_{(1^{\mu_1})}, \ldots ,V_{(1^{\mu_r})}$
is an isomorphism, and the kernel of $\varphi_\mu$ is the sum of
all components $V_\nu,\ \nu \not = \lambda$ appearing in
$V_{(1^{\mu_1})} \otimes \cdots \otimes V_{(1^{\mu_r})}$.

We now turn to the crystal bases of these various  $U_q(\gl_n)$-modules.
It follows from Example~\ref{PEXT} and Theorem~\ref{EXIST} that, if we
choose $u_{(1^k)} = y_{(1^k)}$, $L(1^k)$
is the ${\cal A}$-lattice spanned by the basis of column-shaped quantum
tableaux,
while $B(1^k)$ is the projection of this basis in $L(1^k)/qL(1^k)$.
Using Theorem~\ref{TP} we deduce that a crystal basis of
$V_{(1^{\mu_1})} \otimes \cdots \otimes V_{(1^{\mu_r})}$ can be
constructed by tensoring the crystal bases of the factors. Denote this
basis by $(L^\mu,\, B^\mu)$. Then $\varphi_\mu(L^\mu)$ is a crystal
lattice in $V_\lambda$. According to Theorem~\ref{EXIST} such a lattice
is unique up to an overall scalar multiple which can be determined by
considering the highest weight space.
Here we have
$
(L^\mu)_\lambda = {\cal A} \, y^\mu
$
and
$
\varphi_\mu(y^\mu)
= (y_\lambda)
$.
Therefore, $\varphi_\mu(L^\mu)$ is the crystal lattice $L$ of $V_\lambda$
specified by $L_\lambda = {\cal A} \, (y_\lambda)$, and it does
not depend on $\sigma$. Now $\varphi_\mu$ induces an isomorphism from the
$K$-subspace of $L^\mu/qL^\mu$ spanned by the subset of $B^\mu$ consisting
of those $b$ in the connected component of the crystal graph with source
$y^\mu \ \mod qL^\mu$,
to the $K$-space $L/qL$. It is easy to check from
Theorem~\ref{TP} and the explicit description of the operators $\hat e_i$,
$\hat f_i$ which follows it that this connected component is labelled
by the elements $b=c_1\otimes \cdots \otimes c_r \ \mod L^\mu$ such that the
word $u_b$
obtained by reading $c_1$ from top to bottom, then $c_2$ from top to bottom,
and so on, satisfies: $P(u_b)$ has shape $\lambda$. Hence its image is
$B_\sigma \ \mod qL$. This proves that $B_\sigma$ is a basis of
$V_\lambda$ and that $L = \sum_{b\in B_\sigma} {\cal A}\, b$.
Moreover, $(L,\,B_\sigma \ \mod qL)$ is a crystal basis, which proves
by unicity that $B_\sigma \ \mod qL$
does not depend on $\sigma$. On the other hand, the elements
$b$ of $B^\mu$ which belong to the other connected components are sent to
$0$ in $L/qL$. This means that the tabloids $(\delta)$
such that $P(u_\delta)$ has not shape $\lambda$ belong to $qL$.

Finally, if $\zeta$ is another permutation and $\gamma$ a tabloid of shape
$\zeta(\lambda')$, the fact that $(\delta) \equiv (\gamma) \ \mod qL$
means exactly that the words $u_\delta$ and $u_{\gamma}$ label the same vertex
in the two copies of $\Gamma(V_\lambda)$ to which they belong
in $\Gamma(V_{(1)}^{\otimes k})$, that is, by
Theorem~\ref{TCG}, $P(u_\delta) = P(u_{\gamma})$. \cqfd

\medskip\noindent
{\it Proof of Theorem~\ref{TH2}:} Let $(\delta)$ be a quantum tabloid in
$V_\lambda$.
The $q$-straightening algorithm shows that $(\delta)$ is a
$K[q,\,q^{-1}]$-linear
combination of quantum tableaux. On the other hand, Theorem~\ref{TH3} shows
that $(\delta)$ belongs to the ${\cal A}$-lattice $L$ spanned by the set
${\cal B} = B_{\rm id}$ of quantum tableaux. Therefore, $(\delta)$ is in fact
an element of the $K[q]$-lattice ${\cal L}$ spanned by ${\cal B}$. The other
statements of Theorem~\ref{TH2} follow immediately from Theorem~\ref{TH3}.
\cqfd

\medskip
We point out that the quantum tabloids $(\delta) \in B_\sigma$ are easily
computed
using Sch\"{u}tzen\-ber\-ger's jeu de taquin \cite{S}. For example, the
following graph shows the $\delta$'s giving rise to quantum tabloids
$(\delta)$ which are congruent to the quantum tableau $(\tau)$ modulo $qL$:
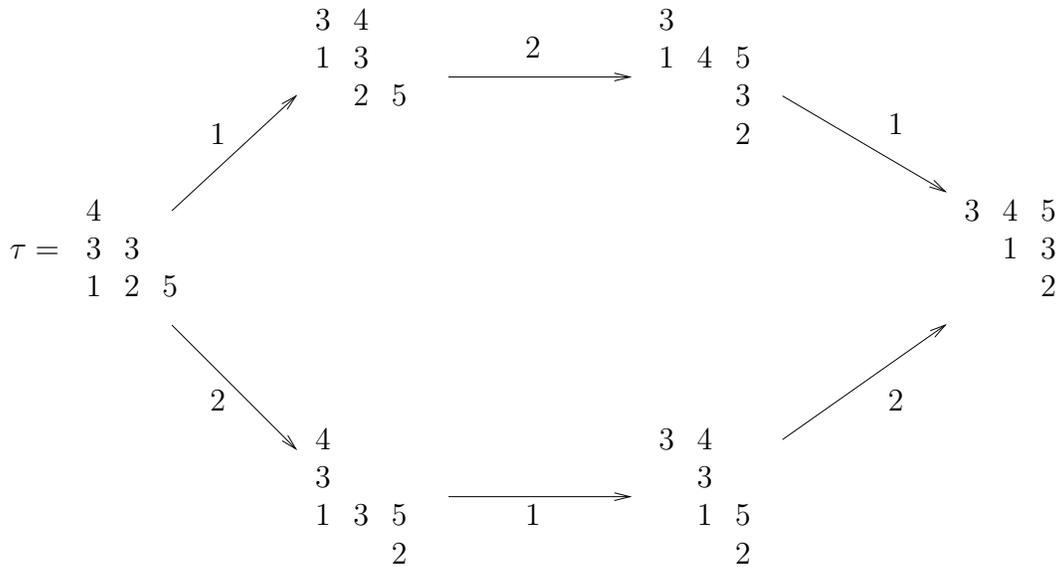
\begin{center}
\setlength{\unitlength}{0.01in}
\begingroup\makeatletter\ifx\SetFigFont\undefined
\def\x#1#2#3#4#5#6#7\relax{\def\x{#1#2#3#4#5#6}}%
\expandafter\x\fmtname xxxxxx\relax \def\y{splain}%
\ifx\x\y   
\gdef\SetFigFont#1#2#3{%
  \ifnum #1<17\tiny\else \ifnum #1<20\small\else
  \ifnum #1<24\normalsize\else \ifnum #1<29\large\else
  \ifnum #1<34\Large\else \ifnum #1<41\LARGE\else
     \huge\fi\fi\fi\fi\fi\fi
  \csname #3\endcsname}%
\else
\gdef\SetFigFont#1#2#3{\begingroup
  \count@#1\relax \ifnum 25<\count@\count@25\fi
  \def\x{\endgroup\@setsize\SetFigFont{#2pt}}%
  \expandafter\x
    \csname \romannumeral\the\count@ pt\expandafter\endcsname
    \csname @\romannumeral\the\count@ pt\endcsname
  \csname #3\endcsname}%
\fi
\fi\endgroup
\begin{picture}(546,310)(0,-10)
\path(85,185)(150,245)
\path(145.478,238.104)(150.000,245.000)(142.765,241.043)
\path(230,255)(325,255)
\path(317.000,253.000)(325.000,255.000)(317.000,257.000)
\path(405,245)(490,195)
\path(482.090,197.332)(490.000,195.000)(484.119,200.780)
\path(230,35)(325,35)
\path(317.000,33.000)(325.000,35.000)(317.000,37.000)
\path(85,125)(150,60)
\path(142.929,64.243)(150.000,60.000)(145.757,67.071)
\path(405,65)(490,125)
\path(484.618,118.753)(490.000,125.000)(482.311,122.020)
\put(160,280){\makebox(0,0)[lb]{\smash{{{\SetFigFont{12}{14.4}{rm}3}}}}}
\put(180,280){\makebox(0,0)[lb]{\smash{{{\SetFigFont{12}{14.4}{rm}4}}}}}
\put(160,260){\makebox(0,0)[lb]{\smash{{{\SetFigFont{12}{14.4}{rm}1}}}}}
\put(180,260){\makebox(0,0)[lb]{\smash{{{\SetFigFont{12}{14.4}{rm}3}}}}}
\put(180,240){\makebox(0,0)[lb]{\smash{{{\SetFigFont{12}{14.4}{rm}2}}}}}
\put(200,240){\makebox(0,0)[lb]{\smash{{{\SetFigFont{12}{14.4}{rm}5}}}}}
\put(40,180){\makebox(0,0)[lb]{\smash{{{\SetFigFont{12}{14.4}{rm}4}}}}}
\put(40,160){\makebox(0,0)[lb]{\smash{{{\SetFigFont{12}{14.4}{rm}3}}}}}
\put(40,140){\makebox(0,0)[lb]{\smash{{{\SetFigFont{12}{14.4}{rm}1}}}}}
\put(60,160){\makebox(0,0)[lb]{\smash{{{\SetFigFont{12}{14.4}{rm}3}}}}}
\put(60,140){\makebox(0,0)[lb]{\smash{{{\SetFigFont{12}{14.4}{rm}2}}}}}
\put(80,140){\makebox(0,0)[lb]{\smash{{{\SetFigFont{12}{14.4}{rm}5}}}}}
\put(160,60){\makebox(0,0)[lb]{\smash{{{\SetFigFont{12}{14.4}{rm}4}}}}}
\put(160,40){\makebox(0,0)[lb]{\smash{{{\SetFigFont{12}{14.4}{rm}3}}}}}
\put(160,20){\makebox(0,0)[lb]{\smash{{{\SetFigFont{12}{14.4}{rm}1}}}}}
\put(180,20){\makebox(0,0)[lb]{\smash{{{\SetFigFont{12}{14.4}{rm}3}}}}}
\put(200,20){\makebox(0,0)[lb]{\smash{{{\SetFigFont{12}{14.4}{rm}5}}}}}
\put(200,0){\makebox(0,0)[lb]{\smash{{{\SetFigFont{12}{14.4}{rm}2}}}}}
\put(340,280){\makebox(0,0)[lb]{\smash{{{\SetFigFont{12}{14.4}{rm}3}}}}}
\put(340,260){\makebox(0,0)[lb]{\smash{{{\SetFigFont{12}{14.4}{rm}1}}}}}
\put(360,260){\makebox(0,0)[lb]{\smash{{{\SetFigFont{12}{14.4}{rm}4}}}}}
\put(380,260){\makebox(0,0)[lb]{\smash{{{\SetFigFont{12}{14.4}{rm}5}}}}}
\put(380,240){\makebox(0,0)[lb]{\smash{{{\SetFigFont{12}{14.4}{rm}3}}}}}
\put(380,220){\makebox(0,0)[lb]{\smash{{{\SetFigFont{12}{14.4}{rm}2}}}}}
\put(340,60){\makebox(0,0)[lb]{\smash{{{\SetFigFont{12}{14.4}{rm}3}}}}}
\put(360,60){\makebox(0,0)[lb]{\smash{{{\SetFigFont{12}{14.4}{rm}4}}}}}
\put(360,40){\makebox(0,0)[lb]{\smash{{{\SetFigFont{12}{14.4}{rm}3}}}}}
\put(360,20){\makebox(0,0)[lb]{\smash{{{\SetFigFont{12}{14.4}{rm}1}}}}}
\put(380,20){\makebox(0,0)[lb]{\smash{{{\SetFigFont{12}{14.4}{rm}5}}}}}
\put(380,0){\makebox(0,0)[lb]{\smash{{{\SetFigFont{12}{14.4}{rm}2}}}}}
\put(500,180){\makebox(0,0)[lb]{\smash{{{\SetFigFont{12}{14.4}{rm}3}}}}}
\put(520,180){\makebox(0,0)[lb]{\smash{{{\SetFigFont{12}{14.4}{rm}4}}}}}
\put(520,160){\makebox(0,0)[lb]{\smash{{{\SetFigFont{12}{14.4}{rm}1}}}}}
\put(540,180){\makebox(0,0)[lb]{\smash{{{\SetFigFont{12}{14.4}{rm}5}}}}}
\put(540,160){\makebox(0,0)[lb]{\smash{{{\SetFigFont{12}{14.4}{rm}3}}}}}
\put(540,140){\makebox(0,0)[lb]{\smash{{{\SetFigFont{12}{14.4}{rm}2}}}}}
\put(105,220){\makebox(0,0)[lb]{\smash{{{\SetFigFont{12}{14.4}{rm}1}}}}}
\put(270,265){\makebox(0,0)[lb]{\smash{{{\SetFigFont{12}{14.4}{rm}2}}}}}
\put(460,225){\makebox(0,0)[lb]{\smash{{{\SetFigFont{12}{14.4}{rm}1}}}}}
\put(105,80){\makebox(0,0)[lb]{\smash{{{\SetFigFont{12}{14.4}{rm}2}}}}}
\put(270,20){\makebox(0,0)[lb]{\smash{{{\SetFigFont{12}{14.4}{rm}1}}}}}
\put(460,80){\makebox(0,0)[lb]{\smash{{{\SetFigFont{12}{14.4}{rm}2}}}}}
\put(0,160){\makebox(0,0)[lb]{\smash{{{\SetFigFont{12}{14.4}{rm}{$\tau =$}}}}}}
\end{picture}
\end{center}
The edges labelled $i$ connect tabloids $\delta ,\, \gamma$ obtained one from
the other by permuting the column lengths of the $i$ th and $(i{+}1)$ th
columns
by means of jeu de taquin. The column readings $u_\delta$ of these tabloids are
distinguished words of the plactic class of $u_\tau$, called {\it frank words}
\cite{LS2}. The combinatorics of frank words occurs in several interesting
problems
(see \cite{LS2,RS,FL}).

\section{Proof of Theorem 1.1} \label{PR}
The proof of Theorem~\ref{TH1} follows the same lines as that of
Theorem~\ref{TH2},
namely we first describe several realizations of the crystal basis of the
bimodule
${\cal F}_q[\mat_n]$ and then deduce Theorem~\ref{TH1} by comparing them, using
the unicity property stated in Theorem~\ref{EXIST}.

We retain the notations of Section~\ref{qG/B} (before Theorem~\ref{TH3}). In
particular, we fix a partition $\lambda$ of $k$, and set $\mu =
\sigma(\lambda')$
for $\sigma \in S_r$. We denote by $W_\lambda$ the subspace of ${\cal
F}_q[\mat_n]$
whose decomposition as a $U_q(\gl_n)$-bimodule is
$ W_\lambda \simeq \oplus_{\nu \le \lambda} V_\nu \otimes V_\nu $.
Thus, $W_{(k)}$ is the homogeneous component of degree $k$ of ${\cal
F}_q[\mat_n]$.
For each $\nu \le \lambda$, we choose a standard Young tableau $\tau_\nu$ of
shape $\nu$. We can now state
\begin{theorem}\label{TH4}
Let $B_{\lambda , \sigma}$ denote the following set of quantum bitabloids:
$$
B_{\lambda , \sigma} = \{(\delta|\delta')\ |\ \delta ,  \delta' \ {\it have \
shape}\
\mu \ {\it and}\ Q(u_\delta) = Q(u_\delta') = \tau_\nu \ {\it for\ some}\
\nu\le\lambda\}
\,.
$$
Then $B_{\lambda , \sigma}$ is a basis of $W_\lambda$. The ${\cal A}$-lattice
$L_\lambda$
spanned by $B_{\lambda , \sigma}$ does not depend on $\sigma$. The projection
$B_\lambda$
of $B_{\lambda , \sigma}$ in $L_\lambda/qL_\lambda$ is also independent of
$\sigma$,
and $(L_\lambda, B_\lambda)$ is a crystal basis of $W_\lambda$. Moreover for
$(\delta | \delta') \in B_{\lambda , \sigma}$ and
$(\gamma | \gamma') \in B_{\lambda , \zeta}$, we have
$$
(\delta | \delta') \equiv (\gamma | \gamma') \ \mod qL_\lambda \quad
\Longleftrightarrow \quad
P(u_\delta) = P(u_\gamma) \ {\it and} \ P(u_{\delta'}) = P(u_{\gamma'})\,.
$$
Finally, the crystal graph attached to $(L_\lambda , B_\lambda)$ is the
bi-coloured
graph given by
$$
(\delta | \delta')  \stackrel{i_l}{\longrightarrow} (\gamma | \delta')
\quad \Longleftrightarrow \quad
f_i^\dagger (\delta | \delta') = (\gamma | \delta')
\quad \Longleftrightarrow \quad
\hat f_i u_\delta = u_\gamma \,,
$$
$$
(\delta | \delta')  \stackrel{i_r}{\longrightarrow} (\delta | \gamma')
\quad \Longleftrightarrow \quad
f_i (\delta | \delta') = (\delta | \gamma')
\quad \Longleftrightarrow \quad
\hat f_i u_{\delta'} = u_{\gamma'} \,.
$$
\end{theorem}
\Proof We imitate the proof of Theorem~\ref{TH3} and introduce the
$K(q)$-linear
map $\Phi_\mu$:
$$
\matrix{
\left(V_{(1^{\mu_1})}\otimes \cdots \otimes V_{(1^{\mu_r})}\right)
\otimes
\left(V_{(1^{\mu_1})}\otimes \cdots \otimes V_{(1^{\mu_r})}\right)
&
\longrightarrow
&
W_\lambda
\cr
(c_1\otimes \cdots \otimes c_r)\otimes (d_1\otimes \cdots \otimes d_r)
&
\longrightarrow
&
(c_1\cdots c_r | d_1\cdots d_r)
}
$$
Here $c_i,\,d_j$ denote elements of the canonical bases of $V_{(1^{\mu_i})}$,
$V_{(1^{\mu_j})}$, {\it i.e.} columns of size $\mu_i$, $\mu_j$, and
$(c_1\cdots c_r | d_1\cdots d_r)$ is the quantum bitabloid formed on these
columns.
The map $\Phi_\mu$ is in fact a homomorphism of $U_q(\gl_n)$-bimodules, and
therefore
sends the crystal lattice ${\cal L}_{\lambda , \sigma}$ spanned by the tensors
$(c_1\otimes \cdots \otimes c_r)\otimes (d_1\otimes \cdots \otimes d_r)$ onto a
crystal lattice $L_{\lambda , \sigma}$ in $W_\lambda$. To describe precisely
$L_{\lambda , \sigma}$, it is enough to determine its submodule $L_{\lambda ,
\sigma}^+$
of highest weight vectors. We shall prove that
$L_{\lambda , \sigma}^+ = \oplus_{\nu \le \lambda} {\cal A} \, (y_\nu |
y_\nu)$.

To do this we have to introduce some notations. We write for short
$V^\mu := V_{(1^{\mu_1})}\otimes \cdots \otimes V_{(1^{\mu_r})}$, and we
consider a source vertex
$$
y =
\bo{\matrix{y_1 \cr \vdots \cr y_{\mu_1}}} \otimes \cdots \otimes
\bo{\matrix{y_{k-\mu_r+1} \cr \vdots \cr y_k}}
$$
of the crystal graph of $V^\mu$. This means that $y_1\cdots y_k$ is a
Yamanouchi word on $\{1,\ldots ,n\}$ such that the columns of $y$ are
increasing
from bottom to top. By definition of a crystal basis, there exists a highest
weight vector $T_y$ in $V^\mu$ such that
$T_y \equiv y \ \mod q{\cal L}_{\lambda , \sigma}$. Now, let
$$
w =
\bo{\matrix{w_1 \cr \vdots \cr w_{\mu_1}}} \otimes \cdots \otimes
\bo{\matrix{w_{k-\mu_r+1} \cr \vdots \cr w_k}}
$$
be any monomial tensor of the same weight $\nu$ as $y$. Then
$\Phi_\mu(T_y \otimes w) = \alpha_{y,w}(q)\, (y_\nu | y_\nu)$ for some scalar
$\alpha_{y,w}(q) \in K(q)$. Indeed,
$$
e_i^\dagger \Phi_\mu(T_y \otimes w) = \Phi_\mu( e_i^\dagger T_y \otimes w) =
\Phi_\mu( (e_i T_y) \otimes w) = 0 \,,
$$
and  on the other hand $\Phi_\mu(T_y \otimes w)$ has weight $(\nu ,\,\nu)$.
Therefore $\Phi_\mu(T_y \otimes w)$ is a scalar multiple of the highest weight
vector $(y_\nu | y_\nu)$. But from the definition of $\Phi_\mu$, we have
\begin{equation}\label{EXP}
\Phi_\mu(T_y \otimes w) = \sum_{\bf u} \kappa_{\bf u}(q)\, ({\bf u} | {\bf w})
\,.
\end{equation}
Here {\bf w} is the tabloid made of the columns of $w$, {\bf u} runs through
the tabloids of the same shape and weight as {\bf w}, and
$ \kappa_{\bf u}(q) \in {\cal A}$. Moreover, $ \kappa_{\bf y}(0) = 1$ and
$ \kappa_{\bf u}(0) = 0$ for ${\bf u} \not = {\bf y}$. Expanding (\ref{EXP})
on the basis
$$\beta_w=\{t_{u_1w_1}\cdots t_{u_kw_k} \ | \ u_1\cdots u_k \ {\rm has \
weight}\ \nu
\ {\rm and}\ w_i = w_j \Longrightarrow u_i \le u_j\}$$
and comparing to the similar expansion of $(y_\nu | y_\nu)$, we obtain by
checking
the coefficient of $t_{w_1w_1}\cdots t_{w_kw_k}$ that
$\alpha_{y,w}(q) =  \kappa_{\bf w}(q)$. Therefore,
$\Phi_\mu(T_y \otimes w) \equiv \delta_{w\,y}\,(y_\nu | y_\nu) \ \mod q(y_\nu |
y_\nu)$,
and if $y'$ is another source vertex of $\Gamma(V^\mu)$, we have
\begin{equation}\label{YY'}
\Phi_\mu(T_y \otimes T_{y'}) \equiv
\delta_{y\,y'}\,(y_\nu | y_\nu) \ \mod q(y_\nu | y_\nu)\,.
\end{equation}
This proves that
$L_{\lambda , \sigma}^+ = \oplus_{\nu \le \lambda} {\cal A} \, (y_\nu |
y_\nu)$,
as stated.

It follows that $L_{\lambda ,\sigma}$ does not in fact depend on
$\sigma$, and we may write $L_{\lambda ,\sigma} = L_\lambda$. We can now
consider
the map induced by $\Phi_\mu$ from
${\cal L}_{\lambda , \sigma}/q{\cal L}_{\lambda , \sigma}$
to $L_\lambda /qL_\lambda$ (we still denote it by $\Phi_\mu$). By (\ref{YY'}),
it maps to zero all the connected components of the crystal graph of
$V^\mu \otimes V^\mu$ with source vertex $y\otimes y'$ such that $y \not = y'$.
This means precisely that
$
(\delta|\delta') \equiv 0 \ \mod qL_\lambda
$
whenever $Q(u_\delta) \not = Q(u_{\delta'})$. On the other hand, to each
$\tau_\nu$ corresponds a unique source vertex $y = y(\tau_\nu)$ such that
$Q(u_y)=\tau_\nu$. The vertices $w\otimes w'$ of the connected components
of $\Gamma(V^\mu \otimes V^\mu)$ with origin $y(\tau_\nu) \otimes y(\tau_\nu)$,
$\nu \le \lambda$, span a subspace of
${\cal L}_{\lambda , \sigma}/q{\cal L}_{\lambda , \sigma}$ isomorphic under
$\Phi_\mu$ to $L_\lambda /qL_\lambda$. Therefore
$(L_\lambda,\, B_{\lambda,\sigma}\ \mod qL_\lambda)$ is a crystal basis of
$W_\lambda$, and this implies that $B_{\lambda,\sigma}\ \mod qL_\lambda$
is independent of $\sigma$. The other statements follow now easily from
Theorem~\ref{TCG}. \cqfd

\medskip\noindent
{\it Proof of Theorem~\ref{TH1}:} Denote by $L$ the crystal lattice of
${\cal F}_q[\mat_n]$ whose submodule of highest weight vectors is equal to
$ \oplus_\nu\, {\cal A} \, (y_\nu | y_\nu)$. It follows from the proof of
Theorem~\ref{TH4} that the crystal lattice $L_\lambda$ of $W_\lambda$ is
nothing but $L\cap W_\lambda$. We also see that $L$ is spanned over ${\cal A}$
by the set of bitableaux $(\tau\,|\,\tau')$. Indeed, if we choose
$\tau_\lambda$
to be the standard Young tableau whose first column contains $1,\,2,\,\ldots
,\,
\lambda_1'$, whose second column contains
$\lambda_1'+1,\,\ldots ,\,\lambda_1'+\lambda_2'$, and so on, we have
$$
BT_\lambda := \{(\tau\,|\,\tau') \ | \ \tau ,\, \tau'\ {\rm are\ Young\
tableaux\ of\ shape}\ \lambda\} \subset B_{\lambda ,{\rm id}}
$$
and $BT_\lambda \ \mod L_\lambda$ is the part of $B_\lambda$ which labels the
connected component of $\Gamma(W_\lambda)$ corresponding to $V_\lambda \otimes
V_\lambda$.
Therefore, $(L, \sqcup_\lambda BT_\lambda \ \mod L)$ is a crystal basis of
${\cal F}_q[\mat_n]$ and $L= \oplus_\lambda \,{\cal A}\,  BT_\lambda$.
On the other hand, taking $\lambda = (k)$ in Theorem~\ref{TH4}, we find that
the set of monomials
$$
B_k=\{t_{i_1j_1}\cdots t_{i_kj_k}\ | \ Q(i_1\cdots i_k) = Q(j_1\cdots j_k) =
\tau_\nu \ {\rm for \ some}\ \nu\}
$$
gives rise to the same crystal basis under the guise
$(L, \sqcup_k B_k \ \mod L)$. Hence, applying again Theorem~\ref{EXIST},
it follows from the description of the crystal graph that
\begin{equation}\label{DER}
t_{i_1j_1}\cdots t_{i_kj_k} \equiv
\left\{\matrix{
(P(i_1\cdots i_k)|P(j_1\cdots j_k)) &\mod qL &{\rm if} \
Q(i_1\cdots i_k)=Q(j_1\cdots j_k) \cr
0           &\mod qL &{\rm otherwise}} \right.
\end{equation}
Finally, it is known that the coefficients of the straightening of a bitabloid
belong
to $K[q,q^{-1}]$. Since we have just proved that they also belong to ${\cal
A}$,
we can replace in (\ref{DER}) the crystal lattice $L$ by the $K[q]$-lattice
$\L$
spanned by the quantum bitableaux. \cqfd


\footnotesize

\begin{thebibliography}{ABCD}

\bibitem{BZ} {\sc  A. D. Berenstein, A. V. Zelevinsky}, {\it
Canonical bases for the quantum group of type $A_r$ and piecewise-linear
combinatorics},
Preprint Institut Gaspard Monge, 1994.

\bibitem{BKW}{\sc C. Burd\'ik, R. C. King, T. A. Welsh}, {\it The explicit
construction of irreducible representations of the quantum algebras
$U_q(sl_n)$},
preprint 1993.

\bibitem{ChPr}{\sc V. Chari, A. Pressley}, {\it A guide to
quantum groups}, Cambridge Universty Press, 1994.

\bibitem{DJM}{\sc E. Date, M. Jimbo, T. Miwa},  Representations
of $U_q(gl_n)$ at $q=0$ and the Robinson-Schensted correspondence,
in {\it Physics and mathematics of strings, Memorial volume of
V. Knizhnik}, L. Brink, D. Friedan, A.M. Polyakov eds., World
Scientific, 1990.

\bibitem{D}{\sc J. D\'esarm\'enien}, {\it An algorithm for the
Rota straightening formula}, Discrete Math. {\bf 30} (1980), 51-68.

\bibitem{DKR}{\sc J. D\'esarm\'enien, J.P.S. Kung, G.C. Rota},
{\it Invariant theory, Young bitableaux and combinatorics}, Adv. Math.
{\bf 27} (1978), 63-92.

\bibitem{Dr}{\sc V. G. Drinfeld}, {\it Hopf algebras and the quantum
Yang-Baxter equation},
Soviet Math. Dokl. {\bf 32} (1985), 254-258.

\bibitem{FL} {\sc W. Fulton, A. Lascoux}, {\it A Pieri formula in the
Grothendieck ring of a flag bundle}, Duke Math. J., {\bf 76} (1994), 711-729.

\bibitem{Gr}{\sc J.A. Green}, Classical invariants and the general linear
group, in {\it Representation theory of finite groups and finite dimensional
algebras}, G.O. Michler, C.M. Ringel eds., Progress in Mathematics 95,
Birkh\"auser,
1991.

\bibitem{HZ}{\sc R.Q. Huang, J.J. Zhang}, {\it Standard basis theorem
for quantum linear groups}, Adv. Math. {\bf 102} (1993), 202-229.

\bibitem{Ji1}{\sc M. Jimbo}, {\it A q-difference analogue of $U(g)$ and the
Yang-Baxter equation}, Lett. Math. Phys. {\bf 10} (1985), 63-69.

\bibitem{Ji2}{\sc M. Jimbo}, {\it A q-analogue of $U(gl(N+1))$, Hecke algebra,
and the Yang-Baxter equation}, Lett. Math. Phys. {\bf 11} (1986), 247-252.

\bibitem{Ka1}{\sc M. Kashiwara}, {\it Crystalizing the $q$-analogue of
universal enveloping algebras}, Commun. Math. Phys. {\bf 133} (1990), 249-260.

\bibitem{Ka2}{\sc M. Kashiwara}, {\it On crystal bases of the $q$-analogue
of universal enveloping algebras},  Duke Math. J. {\bf 63} (1991), 465-516.

\bibitem{KN}{\sc M. Kashiwara, T. Nakashima}, {\it Crystal graphs
for representations of the $q$-analogue of classical Lie algebras},
J. Algebra, {\bf 165} (1994), 295-345.


\bibitem{Kn}{\sc D.E. Knuth}, {\it Permutations, matrices and generalized Young
tableaux}, Pacific. J. Math. {\bf 34} (1970), 709-727.

\bibitem{LLT}{\sc A. Lascoux, B. Leclerc, J.Y. Thibon}, {\it Crystal graphs
and $q$-analogues of weight multiplicities for root systems of type $A_n$},
Lett. Math. Phys. {\it to appear}.

\bibitem{LR} {\sc V. Lakshmibai,  N. Yu. Reshetikhin},
Quantum deformations of $SL_n/B$ and its Schubert varieties,
in {\it Special Functions}, ICM-90 Satellite Conference Proceedings,
M. Kashiwara, T. Miwa eds., Springer.

\bibitem{LS1} {\sc A. Lascoux, M. P. Sch\"{u}tzenberger},  Le mono\"{\i}de
plaxique,
in {\it Noncommutative structures in algebra and geometric combinatorics }
A. de Luca Ed., Quaderni della Ricerca Scientifica del C. N. R., Roma, 1981

\bibitem{LS2}{\sc A. Lascoux, M.P. Sch\"utzenberger}, Keys and standard
bases, in {\it Invariant theory and tableaux}, D. Stanton ed., Springer,
1990.

\bibitem{Le} {\sc B. Leclerc}, {\it On identities satisfied by
minors of a matrix}, Adv. in Math., {\bf 100} (1993), 101-132.

\bibitem{NYM} {\sc M. Noumi, H. Yamada, K. Mimachi}, {\it Finite dimensional
representations of the quantum group $GL_q(n;C)$ and the zonal spherical
functions on $U_q(n-1)\backslash U_q(n)$}, Japan. J. Math., {\bf 19} (1993),
31-80.

\bibitem{RS}{\sc V. Reiner, M. Shimozono}, {\it Key polynomials and a flagged
Littlewood-Richardson rule}, J. Comb. Theory A, {\bf 70} (1995), 107-143.

\bibitem{RTF}{\sc N.Y. Reshetikhin, L.A. Takhtadzhyan, L.D.
Faddeev}, {\it Quantization of Lie groups and Lie algebras}, Leningrad Math.
J.,{\bf 1} (1990), 193--225.

\bibitem{Ro}{\sc G. de B. Robinson}, {\it On the representation theory
of the symmetric group}, Amer. J. Math. {\bf 60} (1938), 745-760.

\bibitem{Sch}{\sc C. Schensted}, {\it Longest increasing and decreasing
subsequences}, Canad. J. Math. {\bf 13} (1961), 179-191.

\bibitem{S}{\sc  M. P. Sch\"{u}tzenberger}, La correspondance de Robinson,
in {\it Combinatoire et repr\'esentation du groupe sym\'etrique}, Strasbourg
1976, D. Foata ed., Springer Lect. Notes in Math. 579.

\bibitem{TT} {\sc E. Taft, J. Towber}, {\it Quantum deformation of
flag schemes and Grassman schemes I -- A $q$-deformation of the
shape-algebra for $GL(n)$}, J. of Algebra, {\bf 142} (1991), 1--36.

\end{thebibliography}
\end{document}